\documentclass[11pt,reqno]{amsart}
\usepackage[T1]{fontenc}
\usepackage{lmodern,amsmath,amssymb,mathrsfs,mathtools,booktabs,enumitem}
\usepackage[margin=1.02in]{geometry}
\usepackage{microtype}
\usepackage{xcolor}
\usepackage[colorlinks=true,linkcolor=blue!45!black,citecolor=blue!45!black,urlcolor=blue!45!black]{hyperref}
\newtheorem{theorem}{Theorem}[section]
\newtheorem{proposition}[theorem]{Proposition}
\newtheorem{lemma}[theorem]{Lemma}
\newtheorem{corollary}[theorem]{Corollary}
\theoremstyle{remark}
\newtheorem{remark}[theorem]{Remark}
\newtheorem*{mainremark}{Remark}
\numberwithin{equation}{section}
\newcommand{\R}{\mathbb R}
\newcommand{\dd}{\,\mathrm d}
\newcommand{\eps}{\varepsilon}

\allowdisplaybreaks[3]
\AtBeginDocument{%
  \addtolength{\abovedisplayshortskip}{0pt plus 1.5pt}%
  \addtolength{\belowdisplayshortskip}{0pt plus 1.5pt}%
}
\hypersetup{pdftitle={Time-Decay Rates toward Viscous Shock Waves for the One-Dimensional Compressible Navier-Stokes-Fourier System},pdfauthor={Le Chen and Yujie Wang},pdfsubject={Non-isentropic shock decay under small H2 perturbations}}
\title[Time-decay rates toward non-isentropic viscous shocks]{Time-Decay Rates toward Viscous Shock Waves for the\
One-Dimensional Compressible Navier--Stokes--Fourier System}
\author{Le Chen}
\author{Yujie Wang}
\thanks{Le Chen. School of Mathematics, Statistics and Mechanics,
Beijing University of Technology, Beijing 100124, P. R. China
(\texttt{chenle@emails.bjut.edu.cn}).}
\thanks{Yujie Wang (corresponding author). School of Mathematics,
Statistics and Mechanics, Beijing University of Technology, Beijing
100124, P. R. China (\texttt{15316072697@163.com}).}
\date{}
\begin{document}
\begin{abstract}
We study small $H^2(\R)$ perturbations of a fixed weak viscous
$3$-shock for the one-dimensional compressible Navier--Stokes--Fourier
system. For each fixed initial perturbation, the error from the
shifted profile in the uniform norm and the shift velocity are
$o(t^{-\frac{1}{4}})$. Nearby stationary profiles give spatial
estimates in terms of the time derivative. Estimates along the
characteristics bound the perturbation at time $t$ by the
dissipation on $[\frac{t}{2},t]$. The exponent $\frac{1}{4}$ is
sharp for a uniform bound on the small $H^2$ ball.
\end{abstract}
\maketitle
\vspace{-0.6\baselineskip}
\begin{quote}

\end{quote}
\setcounter{tocdepth}{1}
\tableofcontents

\section{Introduction}

We consider the one-dimensional compressible Navier--Stokes--Fourier
system in Lagrangian coordinates
\begin{equation}\label{eq:NSF}
 \begin{cases}
 v_t-u_x=0,\\
 u_t+p(v,\theta)_x=\mu\left(\frac{u_x}{v}\right)_x,\\
 c_v\theta_t+p(v,\theta)u_x
       =\kappa\left(\frac{\theta_x}{v}\right)_x+\frac{\mu u_x^2}{v},
 \end{cases}
 \qquad (t,x)\in(0,\infty)\times\R ,
\end{equation}
with initial data $(v,u,\theta)(0,x)=(v_0,u_0,\theta_0)(x)$.
Here $v>0$, $u$ and $\theta>0$ are the specific volume, velocity and
absolute temperature. The pressure is $p(v,\theta)=\frac{R\theta}{v}$.
The constants $R,c_v,\mu,\kappa>0$ are the gas constant, the specific
heat at constant volume, the viscosity and the heat conductivity,
respectively. We write $\gamma=1+\frac{R}{c_v}>1$.

Let $\bar q=(V,U,\Theta)$ be the profile of a weak viscous $3$-shock
of speed $\sigma$. We study the decay of small $H^2(\R)$ perturbations
toward $\bar q(x-\sigma t-X(t))$, with the shift $X$ specified by
\eqref{eq:moving-frame}.

\section{Viscous shock profile and main results}\label{sec:profiles}

\subsection{Viscous shock profile and shift}

Let $q=(v,u,\theta)$, $q_0=(v_0,u_0,\theta_0)$ and
$q_\pm=(v_\pm,u_\pm,\theta_\pm)$. The end-state pressures are
$p_\pm=\frac{R\theta_\pm}{v_\pm}$. The total energy is
$E=c_v\theta+\frac{u^2}{2}$. Multiply the momentum equation by $u$
and add the temperature equation to get
\[
 E_t+(pu)_x=\left(\frac{\kappa\theta_x+\mu uu_x}{v}\right)_x.
\]
The end states satisfy the Rankine--Hugoniot relations
\begin{equation}\label{eq:RH}
 -\sigma[v]-[u]=0,\qquad
 -\sigma[u]+[p]=0,\qquad
 -\sigma[E]+[pu]=0,
\end{equation}
where $[f]=f_+-f_-$. For the $3$-shock under consideration,
\begin{equation}\label{eq:lax}
 \delta=v_+-v_->0,\qquad c_->\sigma>c_+>0,
 \qquad c(v,\theta)=\frac{\sqrt{\gamma R\theta}}{v},
 \quad c_\pm=c(v_\pm,\theta_\pm).
\end{equation}
By \eqref{eq:RH},
\[
 \begin{aligned}
 [u]&=-\sigma\delta,\qquad [p]=\sigma[u]=-\sigma^2\delta,\\
 c_v[\theta]
 &=\frac{[pu]}{\sigma}-\frac{u_++u_-}{2}[u]\\
 &=\frac{p_++p_-}{2\sigma}[u]
       +\frac{u_++u_-}{2}\left(\frac{[p]}{\sigma}-[u]\right)
 =-\frac{p_++p_-}{2}\delta .
 \end{aligned}
\]

Fix a translation of the profile $\bar q=(V,U,\Theta)$. The profile
equations are
\begin{equation}\label{eq:profileODE}
 \begin{gathered}
 -\sigma V'-U'=0,\qquad
 -\sigma U'+P'=\mu\left(\frac{U'}{V}\right)',\\
 -\sigma c_v\Theta'+PU'
 =\kappa\left(\frac{\Theta'}{V}\right)'+\frac{\mu(U')^2}{V},
 \qquad P=\frac{R\Theta}{V}.
 \end{gathered}
\end{equation}
For a sufficiently weak shock, the profile exists and stays in a compact
set of positive states. It satisfies
\begin{equation}\label{eq:profilebounds}
 \begin{gathered}
 V'>0,\qquad U'=-\sigma V',\qquad |\Theta'|\le CV',
 \qquad |\bar q'|\le C\delta^2 e^{-c\delta|\xi|},\\
 |\bar q''|\le C\delta|\bar q'|,\qquad
 |\partial_\xi^j\bar q|\le C_j\delta^3e^{-c\delta|\xi|},\quad j\ge3.
 \end{gathered}
\end{equation}
The bounds through order two follow from~\cite[Lemma~2.3]{KVW}
and the mass equation in \eqref{eq:profileODE}.

For a profile of speed $s$, the three conserved fluxes are
\begin{equation}\label{eq:fluxes}
 \begin{aligned}
 J_0&=-sV-U,\qquad J_1=-sU+P-\frac{\mu U'}{V},\\
 J_2&=-s\left(c_v\Theta+\frac{U^2}{2}\right)
             +PU-\frac{\kappa\Theta'+\mu UU'}{V}.
 \end{aligned}
\end{equation}
The first flux relation determines $U$ and gives
\begin{equation}\label{eq:reducedODE}
 \begin{aligned}
 U&=-sV-J_0,\qquad V'=\frac{V}{\mu s}\left(sU-\frac{R\Theta}{V}+J_1\right),\\
 \Theta'&=\frac{V}{\kappa}
          \left(-sc_v\Theta+\frac{sU^2}{2}+J_1U-J_2\right).
 \end{aligned}
\end{equation}
Substitute $U=-sV-J_0$ in the equations for $V'$ and $\Theta'$.
Their right-hand sides are smooth in $(V,\Theta)$ on the range of
the profile. At $s=\sigma$, differentiation gives
\[
 |(V,\Theta)'''|
 \le C\bigl(|(V,\Theta)''|+|(V,\Theta)'|^2\bigr)
 \le C\delta^3e^{-c\delta|\xi|}.
\]
At each higher order, the linear term uses the preceding derivative.
Every product contains at least two lower-order derivatives.
This gives \eqref{eq:profilebounds} by induction, with
$U=-\sigma V-J_0$.

The weight and the coefficient in the shift equation are
\begin{equation}\label{eq:weight}
 a(\xi)=1+\delta^{-\frac{1}{4}}(V(\xi)-v_-),
 \qquad
 M=\frac{3\gamma(\gamma+1)p_-}{4v_-^2c_-^3}
       \left(1+\frac{2\kappa(\gamma-1)^2}{\mu R\gamma}\right).
\end{equation}
For the weight in \eqref{eq:weight}, we use~\cite[(4.17)]{KVW}
with parameter $\delta^{\frac{3}{4}}$. We take $M$ from the choice
after (4.62) in the proof of~\cite[Lemma~4.4]{KVW} and evaluate it
at the left state.
For $(f,g,h)$, the shift functional is
\begin{equation}\label{eq:feedback}
 \ell(f,g,h)=-\frac{M}{\delta}\int_\R a
 \left(\frac{PV'}{V}f+U'g+\frac{c_v\Theta'}{\Theta}h\right)\dd\xi.
\end{equation}
The shift and perturbation satisfy
\begin{equation}\label{eq:moving-frame}
 \begin{aligned}
 \xi&=x-\sigma t-X(t),\qquad X'=\ell(Z),\qquad X(0)=0,\\
 Z(t,\xi)&=q(t,\xi+\sigma t+X(t))-\bar q(\xi).
 \end{aligned}
\end{equation}
In the moving coordinate $\xi$, we use $q=\bar q+Z$ with
$Z=(\phi,\psi,\zeta)$.

By \eqref{eq:lax} and \eqref{eq:profilebounds},
$\int_\R V'(\xi)\dd\xi=\delta$ and $v_-\le V\le v_+$.
Then $1\le a\le 1+\delta^{\frac{3}{4}}$ by \eqref{eq:weight}. Using
\eqref{eq:feedback} and the Cauchy--Schwarz inequality, we have
\begin{equation}\label{eq:feedbackbounds}
 \begin{aligned}
 |\ell(Z)|&\le \frac{C}{\delta}\int_\R V'|Z|\dd\xi
 \le \min\left\{C_\delta\|Z\|_{L^\infty},
       \left(\frac{C}{\delta}\int_\R V'|Z|^2\dd\xi\right)^{\frac{1}{2}}\right\},\\
 \ell(\bar q')&=-\frac{M}{\delta}\int_\R a
 \left(\frac{P}{V}(V')^2+(U')^2+\frac{c_v}{\Theta}(\Theta')^2\right)\dd\xi<0.
 \end{aligned}
\end{equation}
The weighted norm is $m_\omega(f)=\|\sqrt{\omega}f\|_{L^2}$, where
$\omega(\xi)=\exp(-\eta_0\sqrt{1+\xi^2})$ and $\eta_0>0$ is
sufficiently small and fixed.
For $f\in H^1(\R)$, the fundamental theorem of calculus and
Cauchy--Schwarz give
\begin{equation}\label{eq:anchored-weight}
 \begin{aligned}
 m_\omega(f)^2
 &=\int_\R\omega(\xi)|f(\xi)|^2\dd\xi
 \le \int_\R\omega(\xi)\int_{-1}^{1}
   \bigl(|f(y)|^2+|\xi-y|\|f_\xi\|_{L^2}^2\bigr)\dd y\dd\xi\\
 &\le\int_\R\omega(\xi)
   \left(\int_{-1}^{1}|f(y)|^2\dd y
              +2(1+|\xi|)\|f_\xi\|_{L^2}^2\right)\dd\xi\\
 &\le C_\delta\left(\int_\R V'|f|^2\dd\xi+\|f_\xi\|_{L^2}^2\right).
 \end{aligned}
\end{equation}
For the last inequality, $V'$ is bounded below by a positive constant
on $[-1,1]$, and $\int_\R (1+|\xi|)\omega(\xi)\dd\xi<\infty$.

\subsection{Main results}

Matsumura and Nishihara~\cite{MN} used an energy estimate for the
integrated perturbation to study an isentropic shock. Kawashima and
Matsumura~\cite{KM} treated a heat-conducting gas under the
Nishida--Smoller condition with small zero-integral data. For the full
equations, Liu~\cite{Liu} used mass, momentum and energy to track the
shock translation and the outgoing diffusion waves. Tang and
Zhang~\cite{TZ} later proved stability for a class of large initial
perturbations of weak non-isentropic shocks.

The pointwise semigroup approach to shock stability was developed
in~\cite{ZH}. For dissipative hyperbolic-parabolic systems, see
also~\cite{SK} for the dissipation condition and~\cite{MZG} for
pointwise Green function bounds with real viscosity. These bounds
give sharp $L^p$ decay for weak shocks with small $L^1\cap H^3$ data,
for $2\le p\le\infty$~\cite{MZ}. For localized perturbations,
\cite{Raoofi,HRZ,LZ} give more detailed descriptions of the solution
near a shifted shock and along the outgoing characteristics.

Kang and Vasseur~\cite{KV} proved contraction up to a shift for
barotropic shocks. Kang, Vasseur and Wang applied the method to a
shock and a rarefaction in the barotropic system~\cite{KVWbar}, and
to a shock, a contact wave and a rarefaction in the full
system~\cite{KVW}. We use the single-shock form of the relative
entropy estimate in~\cite{KVW}. For two shocks in the full system,
see~\cite{HM}. Huang and Wang~\cite{HW} use a time-dependent shift
in the vanishing-dissipation limit for single-shock data.

In three dimensions, Wang and Wang~\cite{WW} use a time-dependent
shift to prove stability of a planar barotropic shock under small
$H^2$ perturbations on $\R\times\mathbb T^2$. For the full equations
on $\R^3$, Huang and Xu~\cite{HX} obtain stability and decay with
initial perturbations in $H^3\cap\dot H^{-1}$. In one dimension,
\cite{HHLX} gives decay for an isentropic rarefaction wave with
small $H^2$ data. The decay theory for a viscous contact wave in the
full equations uses an integrated perturbation~\cite{LWX}.

Here the initial perturbation is small in $H^2$. The proof proceeds
as follows.
\begin{enumerate}[leftmargin=*,label=(\roman*),itemsep=0.35em]
\item The energy estimate gives a time between $\frac{t}{2}$ and
      $\frac{3t}{4}$ when the dissipation is small. At that time
      we use four localized moments to choose a stationary profile
      with the current frame speed. Comparing with this profile
      gives \eqref{eq:st-nash}.
\item The acoustic equations still contain the entropy gradient.
      We compare the frozen equations with a parabolic system
      whose kernel is nonnegative. Matching the slow eigenvalue
      through fourth order gives the kernel difference estimate
      \eqref{eq:lift-kernel-bound}. The entropy and acoustic
      components travel at different speeds. This lets us estimate
      the entropy gradient along the acoustic kernels and use
      \eqref{eq:comparison-final} to carry the bound up to $t$.
\item To test the uniform rate, we place an initial perturbation far
      to the right of the shock and let its amplitude tend to zero.
      Diffusion in the contact mode of the right end state gives
      the lower bound in Proposition~\ref{prop:uniform-sharpness}.
\end{enumerate}

\begin{theorem}\label{thm:main}
Fix $R,c_v,\mu,\kappa>0$ and a sufficiently weak nonzero $3$-shock
as above. There are $\eps_\delta,C_\delta>0$ such that, if
\[
 Z_0=q_0-\bar q\in H^2(\R;\R^3),\qquad
 \eps=\|Z_0\|_{H^2}\le\eps_\delta,
\]
then \eqref{eq:NSF}, with the phase \eqref{eq:moving-frame}, has a
unique global strong solution in the small perturbation class. It satisfies
\[
 Z\in C([0,\infty);H^2),\qquad
 (\psi_\xi,\zeta_\xi)\in L^2(0,\infty;H^2),\qquad
 \phi_\xi\in L^2(0,\infty;H^1),
\]
and $v,\theta$ stay uniformly positive. Moreover,
\begin{equation}\label{eq:main-uniform}
 \sup_{t\ge0}(1+t)^{\frac{1}{4}}
       (\|Z(t)\|_\infty+|X'(t)|)\le C_\delta\eps.
\end{equation}
The dissipation is
\begin{equation}\label{eq:def-B}
 B(t)=\|Z_\xi(t)\|_{H^1}^2+m_\omega(Z(t))^2+|X'(t)|^2.
\end{equation}
Then $\int_0^\infty B\le C_\delta\eps^2$. For $t\ge1$,
\begin{equation}\label{eq:main-tail}
 (1+t)(\|Z(t)\|_\infty+|X'(t)|)^4
 \le C_\delta\eps^2\int_{\frac{t}{2}}^{t}B(\tau)\dd\tau.
\end{equation}
In particular, for each fixed initial datum,
\begin{equation}\label{eq:main-little-o}
 \|q(t,\cdot)-\bar q(\cdot-\sigma t-X(t))\|_\infty
        +|X'(t)|=o(t^{-\frac{1}{4}}).
\end{equation}
\end{theorem}

\begin{mainremark}
The rate in \eqref{eq:main-little-o} is for each fixed initial
perturbation. The initial data used for the lower bound in
Proposition~\ref{prop:uniform-sharpness} can depend on $t$.
That lower bound makes $\frac{1}{4}$ optimal for the uniform
estimate \eqref{eq:main-uniform}.
\end{mainremark}

Unless we specify a domain, the norms are over $\R$.
A prime on a profile denotes differentiation in $\xi$.
The constants $C_\delta,c_\delta>0$ may change from line to line.
They may depend on the fixed shock, the physical parameters, the weight
and the spatial cutoffs.

\section{Global energy estimates}\label{sec:energy}

In \eqref{eq:moving-frame}, write $d=X'=\ell(Z)$ for the shift velocity and
$b=\sigma+d$ for the frame speed. The pressure perturbation is $\pi=p-P$,
where $P$ is given in \eqref{eq:profileODE}. The localized dissipation is
\begin{equation}\label{eq:GS-energy}
 \mathcal G^S(Z)=\int_{\mathbb R}V'|Z|^2\,d\xi .
\end{equation}
\begin{lemma}[Weighted relative entropy]\label{lemma:entropy}
Fix a sufficiently weak shock. There are $\varepsilon_1,c_0,c_1>0$ such that
$\sup_{0\le \tau\le T}\|Z(\tau)\|_{H^1}\le\varepsilon_1$ implies
\begin{equation}\label{eq:pure-entropy-differential}
 \frac{\mathrm d}{\mathrm dt}\mathcal E_a(Z)
 +c_0\mathcal G^S(Z)+c_1\mathcal D_0(Z)
 +\frac{\delta}{4M}|d|^2\le0,
\end{equation}
where
\begin{align}
 \mathcal E_a(Z)
 &=\int_{\mathbb R}a\left[
   \frac{\psi^2}{2}
   +R\Theta\Phi\left(\frac{v}{V}\right)
   +c_v\Theta\Phi\left(\frac{\theta}{\Theta}\right)
   \right]\dd\xi,\qquad \Phi(z)=z-1-\log z,
   \label{eq:pure-entropy-functional}\\
 \mathcal D_0(Z)
 &=\int_{\mathbb R}a\left[
       \frac{\mu}{v}\psi_\xi^2
       +\frac{\kappa}{v\theta}\zeta_\xi^2
       \right]\dd\xi .\label{eq:pure-D0}
\end{align}
Moreover,
\[
 c\|Z\|_{L^2}^2\le \mathcal E_a(Z)\le C\|Z\|_{L^2}^2,\qquad
 c\|(\psi_\xi,\zeta_\xi)\|_{L^2}^2
 \le \mathcal D_0(Z)\le C\|(\psi_\xi,\zeta_\xi)\|_{L^2}^2 .
\]
\end{lemma}

\begin{proof}
For the single shock in \eqref{eq:moving-frame}, the energy and diffusion
in~\cite[(4.36) and (4.38)]{KVW} reduce to
\eqref{eq:pure-entropy-functional} and \eqref{eq:pure-D0}. The weight in
\eqref{eq:weight} satisfies~\cite[(4.18)]{KVW} because
$\frac{\delta}{\delta^{\frac{3}{4}}}=\delta^{\frac{1}{4}}$ and
$\delta^{\frac{3}{4}}\le\sqrt{\delta}$ for $0<\delta\le1$. Constant
rarefaction and contact profiles make the interaction terms and profile
residuals vanish. The inequality following (4.71) in~\cite{KVW} gives
\eqref{eq:pure-entropy-differential} with $M$ from \eqref{eq:weight}.

For $z>0$, integrate $\Phi'(z)=\frac{z-1}{z}$ to obtain
\[
 \frac{(z-1)^2}{2\max\{1,z\}}
 \le \Phi(z)=(z-1)^2\int_0^1\frac{\tau}{1+\tau(z-1)}\dd\tau
 \le \frac{(z-1)^2}{2\min\{1,z\}} .
\]
Sobolev embedding and small $\|Z\|_{H^1}$ bound $v,\theta,V,\Theta$ away
from zero and above. The estimate for $\Phi$, applied to $\frac{v}{V}$ and
$\frac{\theta}{\Theta}$, gives the bounds for $\mathcal E_a$ and $\mathcal
D_0$ since $a$ is bounded above and below.
\end{proof}

\begin{proposition}[Global $H^2$ estimate]\label{prop:global}
Fix a sufficiently weak nonzero shock and the shift specified above. There
is $\varepsilon_\delta>0$ such that each $Z_0\in H^2(\mathbb R)$ with
$\|Z_0\|_{H^2}\le\varepsilon_\delta$ gives a unique global strong solution
with positive specific volume and temperature. In the moving coordinate,
\begin{equation}\label{eq:global-regularity}
 Z\in C([0,\infty);H^2),\qquad
 \phi_\xi\in L^2(0,\infty;H^1),\qquad
 (\psi_\xi,\zeta_\xi)\in L^2(0,\infty;H^2).
\end{equation}
For every $0\le s<t$,
\begin{align}
 &\sup_{s\le\tau\le t}\|Z(\tau)\|_{H^2}^2
 +\int_s^t\left[
   \|\phi_\xi\|_{H^1}^2
   +\|(\psi_\xi,\zeta_\xi)\|_{H^2}^2
   +\mathcal G^S(Z)+|d|^2
   \right]d\tau
 \le C_\delta\|Z(s)\|_{H^2}^2 .\label{eq:global-H2}
\end{align}
\end{proposition}

\begin{proof}
We derive the estimates for a smooth solution in a small $H^2$ neighborhood
of the shock and then justify them with the regularity in
\eqref{eq:global-regularity}.

\paragraph{First-order estimates.}
Since the shock is fixed and $\delta>0$, integrating
\eqref{eq:pure-entropy-differential} gives
\[
 \begin{aligned}
 &\|Z(t)\|_{L^2}^2+\int_s^t\Bigl[
 \mathcal G^S(Z)+\|(\psi_\xi,\zeta_\xi)\|_{L^2}^2+|d|^2\Bigr]\dd\tau\\
 &\quad\le C_\delta\mathcal E_a(Z(s))
 \le C_\delta\|Z(s)\|_{L^2}^2 .
 \end{aligned}
\]
The profile residuals and the contact-profile derivatives vanish for the
single shock. Thus \cite[(5.8) and (5.11)]{KVW} and the preceding $L^2$
bound give, on $[s,t]$,
\begin{align}
 &\|\phi_\xi(t)\|_{L^2}^2+\int_s^t\|\phi_\xi\|_{L^2}^2\dd\tau\notag\\
 &\quad\le C\bigl(\|Z(s)\|_{H^1}^2+\|Z(t)\|_{L^2}^2\bigr)
 +C\int_s^t\Bigl[\|(\psi_\xi,\zeta_\xi)\|_{L^2}^2
 +\delta\bigl(|d|^2+\mathcal G^S(Z)\bigr)
 +\varepsilon_1^2\|\psi_{\xi\xi}\|_{L^2}^2\Bigr]\dd\tau\notag\\
 &\quad\le C_\delta\|Z(s)\|_{H^1}^2+
 C\varepsilon_1^2\int_s^t\|\psi_{\xi\xi}\|_{L^2}^2\dd\tau,
 \label{eq:H1-density}
\end{align}
\begin{equation}\label{eq:H1-parabolic}
 \begin{aligned}
 &\|(\psi_\xi,\zeta_\xi)(t)\|_{L^2}^2+
 \int_s^t\|(\psi_{\xi\xi},\zeta_{\xi\xi})\|_{L^2}^2\dd\tau\\
 &\quad\le C\bigl(\|Z(s)\|_{H^1}^2+\delta\|Z(t)\|_{L^2}^2\bigr)
 +C_*\int_s^t\|\phi_\xi\|_{L^2}^2\dd\tau\\
 &\qquad+C\int_s^t\Bigl[\|(\psi_\xi,\zeta_\xi)\|_{L^2}^2
              +\delta\bigl(|d|^2+\mathcal G^S(Z)\bigr)\Bigr]\dd\tau\\
 &\quad\le C_\delta\|Z(s)\|_{H^1}^2+C_*\int_s^t\|\phi_\xi\|_{L^2}^2\dd\tau.
 \end{aligned}
\end{equation}
Adding \eqref{eq:H1-density} and $\frac{1}{2C_*}$ times
\eqref{eq:H1-parabolic} and using
$C\varepsilon_1^2\le\frac{1}{4C_*}$, we obtain
\[
 \begin{aligned}
 &\|\phi_\xi(t)\|_2^2+\frac{1}{2C_*}\|(\psi_\xi,\zeta_\xi)(t)\|_2^2\\
 &\qquad+\int_s^t\left[\frac{1}{2}\|\phi_\xi\|_2^2+
              \frac{1}{4C_*}\|(\psi_{\xi\xi},\zeta_{\xi\xi})\|_2^2
                   \right]\dd\tau
 \le C_\delta\|Z(s)\|_{H^1}^2 .
 \end{aligned}
\]
For
\begin{equation}\label{eq:B1-energy-definition}
 \mathcal B_1=
 \|\phi_\xi\|_2^2+\|(\psi_\xi,\zeta_\xi)\|_{H^1}^2
 +\mathcal G^S+|d|^2,
\end{equation}
this and the $L^2$ bound give
\begin{equation}\label{eq:homogeneous-H1}
 \sup_{s\le\tau\le t}\|Z(\tau)\|_{H^1}^2
 +\int_s^t\mathcal B_1(\tau)\,d\tau
 \le C_\delta\|Z(s)\|_{H^1}^2 .
\end{equation}

\paragraph{Second-order energy and the pressure cancellation.}
The material derivative is $D_t=\partial_t-b\partial_\xi$.
The viscous and heat flux differences are
\begin{equation}\label{eq:energy-fluxes}
 J_v=\frac{u_\xi}{v}-\frac{U'}{V},\qquad
 J_\theta=\frac{\theta_\xi}{v}-\frac{\Theta'}{V}.
\end{equation}
The perturbation equations are
\begin{equation}\label{eq:energy-perturbation-system}
 \begin{aligned}
 D_t\phi&=\psi_\xi+dV',\\
 D_t\psi&=-\pi_\xi+\mu(J_v)_\xi+dU',\\
 c_vD_t\zeta&=-p\psi_\xi-\pi U'
      +\kappa(J_\theta)_\xi
      +\mu\left(\frac{u_\xi^2}{v}-\frac{(U')^2}{V}\right)+c_vd\Theta'.
 \end{aligned}
\end{equation}
On an interval where $\|Z\|_{H^2}\le\varepsilon$, the positive bounds
for $v$ and $\theta$ give
\begin{equation}\label{eq:second-order-weight}
 c\|Z_{\xi\xi}\|_2^2\le\frac{1}{2}\int_\R
 \left(\frac{R\theta}{v^2}\phi_{\xi\xi}^2+\psi_{\xi\xi}^2
       +\frac{c_v}{\theta}\zeta_{\xi\xi}^2\right)\dd\xi
 \le C\|Z_{\xi\xi}\|_2^2 .
\end{equation}
By \eqref{eq:profilebounds} and the Sobolev inequality,
\begin{equation}\label{eq:first-derivative-smallness}
 \|q_\xi\|_\infty\le\|\bar q'\|_\infty+C\|Z\|_{H^2}
 \le C(\delta^2+\varepsilon).
\end{equation}

Differentiate \eqref{eq:energy-perturbation-system} twice in $\xi$. Test the
three equations with $\frac{R\theta}{v^2}\phi_{\xi\xi}$, $\psi_{\xi\xi}$ and
$\frac{\zeta_{\xi\xi}}{\theta}$, respectively. Integrate the pressure term
by parts and use $\frac{p}{\theta}=\frac{R}{v}$ to obtain
\begin{equation}\label{eq:pressure-weak-pairing}
 \begin{aligned}
 &\int\left(\frac{R\theta}{v^2}\phi_{\xi\xi}\psi_{\xi\xi\xi}
       -\psi_{\xi\xi}\pi_{\xi\xi\xi}
       -\frac{p}{\theta}\zeta_{\xi\xi}\psi_{\xi\xi\xi}\right)\dd\xi\\
 &\quad=\int\psi_{\xi\xi\xi}\left(\pi_{\xi\xi}
       +\frac{R\theta}{v^2}\phi_{\xi\xi}
       -\frac{R}{v}\zeta_{\xi\xi}\right)\dd\xi .
 \end{aligned}
\end{equation}
Differentiating $p=\frac{R\theta}{v}$ and $P=\frac{R\Theta}{V}$ gives
\begin{equation}\label{eq:pressure-second-decomposition}
 \begin{aligned}
 &\pi_{\xi\xi}+\frac{R\theta}{v^2}\phi_{\xi\xi}
              -\frac{R}{v}\zeta_{\xi\xi}\\
 &\quad=-\left(\frac{R\theta}{v^2}-\frac{R\Theta}{V^2}\right)V''
         +\left(\frac{R}{v}-\frac{R}{V}\right)\Theta''\\
 &\qquad+2R\left(\frac{\theta v_\xi^2}{v^3}
              -\frac{\Theta(V')^2}{V^3}
              -\frac{v_\xi\theta_\xi}{v^2}
              +\frac{V'\Theta'}{V^2}\right).
 \end{aligned}
\end{equation}
The product rule gives the viscous and heat flux derivatives
\begin{align}
 &(J_v)_{\xi\xi}-\frac{\psi_{\xi\xi\xi}}{v}
                  +\frac{u_\xi}{v^2}\phi_{\xi\xi}\notag\\
 &\quad=\left(\frac{1}{v}-\frac{1}{V}\right)U'''
        -\left(\frac{u_\xi}{v^2}-\frac{U'}{V^2}\right)V''
        -2\left(\frac{v_\xi u_{\xi\xi}}{v^2}
                         -\frac{V'U''}{V^2}\right)\notag\\
 &\qquad+2\left(\frac{u_\xi v_\xi^2}{v^3}
                  -\frac{U'(V')^2}{V^3}\right),
 \label{eq:momentum-flux-second}\displaybreak[3]\\
 &(J_\theta)_{\xi\xi}-\frac{\zeta_{\xi\xi\xi}}{v}
                  +\frac{\theta_\xi}{v^2}\phi_{\xi\xi}\notag\\
 &\quad=\left(\frac{1}{v}-\frac{1}{V}\right)\Theta'''
        -\left(\frac{\theta_\xi}{v^2}-\frac{\Theta'}{V^2}\right)V''
        -2\left(\frac{v_\xi\theta_{\xi\xi}}{v^2}
                         -\frac{V'\Theta''}{V^2}\right)\notag\\
 &\qquad+2\left(\frac{\theta_\xi v_\xi^2}{v^3}
                  -\frac{\Theta'(V')^2}{V^3}\right).
 \label{eq:heat-flux-second}
\end{align}
The bounds \eqref{eq:profilebounds}, \eqref{eq:first-derivative-smallness}
and \eqref{eq:anchored-weight} give
\begin{equation}\label{eq:second-order-remainders}
 \begin{aligned}
 &\left\|\pi_{\xi\xi}+\frac{R\theta}{v^2}\phi_{\xi\xi}
                     -\frac{R}{v}\zeta_{\xi\xi}\right\|_2\\
 &\quad+\left\|(J_v)_{\xi\xi}-\frac{\psi_{\xi\xi\xi}}{v}
                  +\frac{u_\xi}{v^2}\phi_{\xi\xi}\right\|_2
       +\left\|(J_\theta)_{\xi\xi}-\frac{\zeta_{\xi\xi\xi}}{v}
                  +\frac{\theta_\xi}{v^2}\phi_{\xi\xi}\right\|_2\\
 &\quad\le C(\delta^2+\varepsilon)
        \bigl(\|Z_\xi\|_2+\|(\psi_{\xi\xi},\zeta_{\xi\xi})\|_2\bigr)
      +C\bigl\|(|\bar q'|+|\bar q''|+|\bar q'''|)Z\bigr\|_2\\
 &\quad\le C_\delta\bigl(\|Z_\xi\|_2
          +\|(\psi_{\xi\xi},\zeta_{\xi\xi})\|_2+m_\omega(Z)\bigr)
 \le C_\delta\mathcal B_1^{\frac{1}{2}} .
 \end{aligned}
\end{equation}

Integrate the viscous and heat conduction terms by parts and add
\eqref{eq:pressure-weak-pairing}. Using
\eqref{eq:first-derivative-smallness}, \eqref{eq:second-order-remainders}
and Young's inequality gives
\begin{align*}
 &\int\psi_{\xi\xi\xi}\left(\pi_{\xi\xi}
       +\frac{R\theta}{v^2}\phi_{\xi\xi}
       -\frac{R}{v}\zeta_{\xi\xi}\right)\dd\xi\\*
 &\quad-\mu\int\psi_{\xi\xi\xi}(J_v)_{\xi\xi}\dd\xi
       -\kappa\int\left(\frac{\zeta_{\xi\xi\xi}}{\theta}
             -\frac{\theta_\xi\zeta_{\xi\xi}}{\theta^2}\right)
                 (J_\theta)_{\xi\xi}\dd\xi\\*
 &\quad+\int\left(\frac{\mu}{v}\psi_{\xi\xi\xi}^2
                  +\frac{\kappa}{v\theta}\zeta_{\xi\xi\xi}^2\right)\dd\xi\\
 &\quad\le C\left(\|(\psi_{\xi\xi\xi},\zeta_{\xi\xi\xi})\|_2
                        +\mathcal B_1^{\frac{1}{2}}\right)
          \left((\delta^2+\varepsilon)\|\phi_{\xi\xi}\|_2
                        +C_\delta\mathcal B_1^{\frac{1}{2}}\right)\\
 &\quad\le\frac{1}{2}\int\left(\frac{\mu}{v}\psi_{\xi\xi\xi}^2
                  +\frac{\kappa}{v\theta}\zeta_{\xi\xi\xi}^2\right)\dd\xi
     +C(\delta^2+\varepsilon)^2\|\phi_{\xi\xi}\|_2^2
     +C_\delta\mathcal B_1 .
\end{align*}

\paragraph{Coefficient derivatives and viscous heating.}
The equations for $v$ and $\theta$ give
\begin{equation}\label{eq:actual-material-equations}
 D_tv=u_\xi,\qquad
 c_vD_t\theta=\kappa\left(\frac{\theta_\xi}{v}\right)_\xi
                   -pu_\xi+\frac{\mu u_\xi^2}{v}.
\end{equation}
Because $b$ is independent of $\xi$, the time and transport terms give the
time derivative of the energy in \eqref{eq:second-order-weight} and terms
involving the material derivatives of its coefficients. By
\eqref{eq:actual-material-equations}, \eqref{eq:profilebounds},
\eqref{eq:first-derivative-smallness} and Sobolev embedding,
\begin{equation}\label{eq:temperature-weight-estimate}
 \begin{aligned}
 &\frac{1}{2}\int\left(
 \left|D_t\left(\frac{R\theta}{v^2}\right)\right|\phi_{\xi\xi}^2
 +\left|D_t\left(\frac{c_v}{\theta}\right)\right|\zeta_{\xi\xi}^2
 \right)\dd\xi\\
 &\quad\le C\int\left(
 |\theta_{\xi\xi}|+|\theta_\xi v_\xi|+|u_\xi|+u_\xi^2
 \right)(\phi_{\xi\xi}^2+\zeta_{\xi\xi}^2)\dd\xi\\
 &\quad\le C\left(
 \|\zeta_{\xi\xi}\|_2^{\frac{1}{2}}
 \|\zeta_{\xi\xi\xi}\|_2^{\frac{1}{2}}+\delta^2+\varepsilon
 \right)\|Z_{\xi\xi}\|_2^2\\
 &\quad\le \eta\|\zeta_{\xi\xi\xi}\|_2^2
 +C_\eta\left(\|\zeta_{\xi\xi}\|_2^{\frac{2}{3}}
             \|Z_{\xi\xi}\|_2^{\frac{2}{3}}+\delta^2+\varepsilon
       \right)\|Z_{\xi\xi}\|_2^2\\
 &\quad\le \eta\|\zeta_{\xi\xi\xi}\|_2^2
 +C_\eta(\delta^2+\varepsilon)\|\phi_{\xi\xi}\|_2^2
 +C_\eta\mathcal B_1
 \end{aligned}
\end{equation}
for any $\eta>0$. The last two bounds use Young's inequality and
$\|Z\|_{H^2}\le\varepsilon\le1$.

The product rule for the viscous heating term gives
\[
 \begin{aligned}
 \left(\frac{u_\xi^2}{v}-\frac{(U')^2}{V}\right)_{\xi\xi}
 &=\frac{2u_\xi}{v}\psi_{\xi\xi\xi}
   -\frac{u_\xi^2}{v^2}\phi_{\xi\xi}
   +\frac{2}{v}\psi_{\xi\xi}^2+\frac{4U''}{v}\psi_{\xi\xi}\\
 &\quad+2\left(\frac{u_\xi}{v}-\frac{U'}{V}\right)U'''
   -\left(\frac{u_\xi^2}{v^2}-\frac{(U')^2}{V^2}\right)V''
   +2\left(\frac{1}{v}-\frac{1}{V}\right)(U'')^2\\
 &\quad-4\left(\frac{u_\xi u_{\xi\xi}v_\xi}{v^2}
               -\frac{U'U''V'}{V^2}\right)
   +2\left(\frac{u_\xi^2v_\xi^2}{v^3}-\frac{(U')^2(V')^2}{V^3}\right).
 \end{aligned}
\]
Using \eqref{eq:profilebounds}, \eqref{eq:anchored-weight} and
\eqref{eq:first-derivative-smallness}, we obtain
\begin{align}
 &\left|\mu\int_\R\frac{\zeta_{\xi\xi}}{\theta}
             \left(\frac{u_\xi^2}{v}-\frac{(U')^2}{V}\right)_{\xi\xi}
              \dd\xi\right|\notag\\*
 &\quad\le C\|\zeta_{\xi\xi}\|_2\left(
 (\delta^2+\varepsilon)\|\psi_{\xi\xi\xi}\|_2
 +(\delta^2+\varepsilon)^2\|\phi_{\xi\xi}\|_2
 +C_\delta\mathcal B_1^{\frac{1}{2}}\right)
 +C\int|\psi_{\xi\xi}|^2|\zeta_{\xi\xi}|\dd\xi\notag\\
 &\quad\le C\|\zeta_{\xi\xi}\|_2\left(
 (\delta^2+\varepsilon)\|\psi_{\xi\xi\xi}\|_2
 +(\delta^2+\varepsilon)^2\|\phi_{\xi\xi}\|_2
 +C_\delta\mathcal B_1^{\frac{1}{2}}\right)\notag\\*
 &\qquad+C\|\zeta_{\xi\xi}\|_2
              \|\psi_{\xi\xi}\|_2^{\frac{3}{2}}
              \|\psi_{\xi\xi\xi}\|_2^{\frac{1}{2}}\notag\\
 &\quad\le \eta\|\psi_{\xi\xi\xi}\|_2^2
          +C(\delta^2+\varepsilon)^2\|\phi_{\xi\xi}\|_2^2
          +C_{\eta,\delta}\mathcal B_1 .
 \label{eq:heating-product-estimate}
\end{align}
This uses the Sobolev inequality for $\psi_{\xi\xi}$ and Young's inequality.

Use $|\pi|+|\pi_\xi|\le C(|Z|+|Z_\xi|)$ with \eqref{eq:profilebounds},
\eqref{eq:anchored-weight}, \eqref{eq:first-derivative-smallness},
\eqref{eq:pressure-second-decomposition} and
\eqref{eq:second-order-remainders} to obtain
\[
 \begin{aligned}
 &\left|\int\frac{\zeta_{\xi\xi}}{\theta}
 \left(2p_\xi\psi_{\xi\xi}+p_{\xi\xi}\psi_\xi+(\pi U')_{\xi\xi}\right)
 \dd\xi\right|\\
 &\quad+|d|\left|\int\left(
 \frac{R\theta}{v^2}\phi_{\xi\xi}V'''+\psi_{\xi\xi}U'''
 +\frac{c_v}{\theta}\zeta_{\xi\xi}\Theta'''
 \right)\dd\xi\right|\\
 &\quad\le C\|\zeta_{\xi\xi}\|_2\left(
 (\|\psi_\xi\|_\infty+\|U'\|_\infty)\|\pi_{\xi\xi}\|_2
 + C_\delta\mathcal B_1^{\frac{1}{2}}\right)
 +C|d|\|\bar q'''\|_2\|Z_{\xi\xi}\|_2\\
 &\quad\le C(\delta^2+\varepsilon)
       (\|\zeta_{\xi\xi}\|_2+|d|)\|\phi_{\xi\xi}\|_2+C_\delta\mathcal B_1\\
 &\quad\le C(\delta^2+\varepsilon)\|\phi_{\xi\xi}\|_2^2+C_\delta\mathcal B_1 .
 \end{aligned}
\]
Add the bounds for the pressure, diffusion and the terms above. Choose
$\eta$ small enough to absorb the third derivatives in
\eqref{eq:temperature-weight-estimate} and
\eqref{eq:heating-product-estimate}. For sufficiently small $\delta$ and
$\varepsilon$, this gives
\begin{equation}\label{eq:second-order-physical-energy}
 \begin{aligned}
 &\frac{1}{2}\frac{\mathrm d}{\mathrm dt}\int_\R
    \left(\frac{R\theta}{v^2}\phi_{\xi\xi}^2+\psi_{\xi\xi}^2
                       +\frac{c_v}{\theta}\zeta_{\xi\xi}^2\right)\dd\xi
     +c\|(\psi_{\xi\xi\xi},\zeta_{\xi\xi\xi})\|_2^2\\
 &\quad\le C_\eta(\delta^2+\varepsilon)\|\phi_{\xi\xi}\|_2^2
       +C_{\eta,\delta}\mathcal B_1
 \le \eta\|\phi_{\xi\xi}\|_2^2+C_{\eta,\delta}\mathcal B_1 .
 \end{aligned}
\end{equation}

\paragraph{Estimate for the second derivative of the specific volume.}
Since $b$ is independent of $\xi$, integration by parts in
\eqref{eq:energy-perturbation-system} gives
\begin{align}
 &\frac{\mathrm d}{\mathrm dt}\int_{\mathbb R}\psi_\xi\phi_{\xi\xi}\dd\xi
  =\int_{\mathbb R}\left(\phi_{\xi\xi}D_t\psi_\xi
             +\psi_\xi D_t\phi_{\xi\xi}\right)\dd\xi
  \label{eq:density-cross-definition}\\
 &\quad=-\int_{\mathbb R}\left(\pi_{\xi\xi}-\mu(J_v)_{\xi\xi}\right)
       \phi_{\xi\xi}\dd\xi-\|\psi_{\xi\xi}\|_2^2
  +d\int_{\mathbb R}\left(U''\phi_{\xi\xi}+\psi_\xi V'''\right)\dd\xi .
  \label{eq:density-cross-identity}
\end{align}
For sufficiently small $\delta$ and $\varepsilon$,
\eqref{eq:first-derivative-smallness} and the positive bounds for $v$ and
$\theta$ give $\frac{R\theta-\mu u_\xi}{v^2}\ge 2c_1>0$. Use
\eqref{eq:pressure-second-decomposition} and \eqref{eq:momentum-flux-second}
in \eqref{eq:density-cross-identity}, then apply
\eqref{eq:second-order-remainders} to obtain
\begin{equation}\label{eq:density-cross-inequality}
 \begin{aligned}
 &-\frac{\mathrm d}{\mathrm dt}\int_{\mathbb R}\psi_\xi\phi_{\xi\xi}\dd\xi
     +c_1\|\phi_{\xi\xi}\|_2^2\\
 &\quad\le \int_{\mathbb R}
        \left(\frac{R}{v}\zeta_{\xi\xi}-\frac{\mu}{v}\psi_{\xi\xi\xi}\right)
        \phi_{\xi\xi}\dd\xi
     +C_\delta\mathcal B_1^{\frac{1}{2}}\|\phi_{\xi\xi}\|_2
     +C_\delta\mathcal B_1\\
 &\qquad-\int_{\mathbb R}
          \left(\frac{R\theta-\mu u_\xi}{v^2}-c_1\right)
          \phi_{\xi\xi}^2\dd\xi\\
 &\quad\le \left(C\|\psi_{\xi\xi\xi}\|_2+C_\delta\mathcal B_1^{\frac{1}{2}}\right)
          \|\phi_{\xi\xi}\|_2+C_\delta\mathcal B_1-c_1\|\phi_{\xi\xi}\|_2^2\\
 &\quad\le C_1\|\psi_{\xi\xi\xi}\|_2^2+C_\delta\mathcal B_1 .
 \end{aligned}
\end{equation}
The constants $c_1$ and $C_1$ are uniform in $\delta$ and $\varepsilon$.

Choose $\rho>0$ with $\rho C_1\le\frac{c}{2}$ and take $\eta\le\frac{\rho
c_1}{2}$ in \eqref{eq:second-order-physical-energy}. Add $\rho$ times
\eqref{eq:density-cross-inequality} to
\eqref{eq:second-order-physical-energy} to obtain
\begin{equation}\label{eq:closed-second-order-energy}
 \begin{aligned}
 &\frac{\mathrm d}{\mathrm dt}\int_\R
   \left[\frac{1}{2}\left(\frac{R\theta}{v^2}\phi_{\xi\xi}^2
                     +\psi_{\xi\xi}^2+\frac{c_v}{\theta}\zeta_{\xi\xi}^2
                        \right)-\rho\psi_\xi\phi_{\xi\xi}\right]\dd\xi\\
 &\qquad+\frac{\rho c_1}{2}\|\phi_{\xi\xi}\|_2^2
          +\frac{c}{2}\|(\psi_{\xi\xi\xi},\zeta_{\xi\xi\xi})\|_2^2
 \le C_\delta\mathcal B_1 .
 \end{aligned}
\end{equation}
By Young's inequality,
\[
 \begin{aligned}
 \rho\left|\int_\R\psi_\xi\phi_{\xi\xi}\dd\xi\right|
 &\le\frac{1}{4}\int_\R\frac{R\theta}{v^2}\phi_{\xi\xi}^2\dd\xi
     +\rho^2\int_\R\frac{v^2}{R\theta}\psi_\xi^2\dd\xi\\
 &\le\frac{1}{4}\int_\R\left(\frac{R\theta}{v^2}\phi_{\xi\xi}^2
                   +\psi_{\xi\xi}^2+\frac{c_v}{\theta}\zeta_{\xi\xi}^2
                    \right)\dd\xi+C_\rho\|\psi_\xi\|_2^2 .
 \end{aligned}
\]
Integrate \eqref{eq:closed-second-order-energy} on $[s,\tau]$ for
$s\le\tau\le t$. Apply the preceding bound at both endpoints and use
\eqref{eq:second-order-weight} and \eqref{eq:homogeneous-H1} to obtain
\[
 \begin{aligned}
 &\sup_{s\le\tau\le t}\|Z(\tau)\|_{H^2}^2
 +\int_s^t\left(\|\phi_{\xi\xi}\|_2^2
       +\|(\psi_{\xi\xi\xi},\zeta_{\xi\xi\xi})\|_2^2+\mathcal B_1\right)\dd\tau\\
 &\quad\le C\left(\|Z_{\xi\xi}(s)\|_2^2+
                        \sup_{s\le\tau\le t}\|Z(\tau)\|_{H^1}^2\right)
     +C_\delta\int_s^t\mathcal B_1(\tau)\dd\tau\\
 &\quad\le C_\delta\left(\|Z_{\xi\xi}(s)\|_2^2+\|Z(s)\|_{H^1}^2\right)
 \le C_\delta\|Z(s)\|_{H^2}^2 .
 \end{aligned}
\]
This proves \eqref{eq:global-H2} on the bootstrap interval.

\paragraph{Regularity and continuation.}
\cite[Proposition~3.1]{KVW} gives local existence and uniqueness in $H^1$.
Choose its reference function with the same end states as $\bar q$. The
difference from $\bar q$ lies in $H^2$ by \eqref{eq:profilebounds}.

For local existence in $H^2$, start in the frame $x-\sigma t$ with $d=0$ in
\eqref{eq:energy-perturbation-system}. Use transport for the mass equation
and the heat operators with coefficients $\frac{\mu}{v_-}$ and
$\frac{\kappa}{c_v v_-}$ for the other equations. The fluxes satisfy
\[
 (J_v,J_\theta)-\frac{1}{v_-}(\psi_\xi,\zeta_\xi)
 =\left(\frac{1}{v}-\frac{1}{v_-}\right)(\psi_\xi,\zeta_\xi)
   -\frac{\phi}{vV}(U',\Theta').
\]
On the ball $\|Z\|_{C([0,T];H^2)}+\|(\psi,\zeta)\|_{L^2(0,T;H^3)}
\le\varepsilon$, the product estimate and \eqref{eq:profilebounds} give
\[
 \begin{aligned}
 &\left\|\left[\left(\frac{1}{v}-\frac{1}{v_-}\right)
                       (\psi_\xi,\zeta_\xi)\right]_\xi
                         \right\|_{L^2(0,T;H^1)}\\
 &\quad\le C\left(
        \left\|\frac{1}{v}-\frac{1}{v_-}\right\|_{L^\infty([0,T]\times\mathbb R)}
       +\left\|\left(\frac{1}{v}\right)_\xi\right\|_{L^\infty(0,T;H^1)}
            \right)\|(\psi,\zeta)\|_{L^2(0,T;H^3)}\\
 &\quad\le C(\delta+\varepsilon)\|(\psi,\zeta)\|_{L^2(0,T;H^3)} .
 \end{aligned}
\]
The coefficient estimate
$\left\|\frac{1}{v_1}-\frac{1}{v_2}\right\|_{H^2}\le C\|v_1-v_2\|_{H^2}$
gives the flux term a Lipschitz constant $C(\delta+\varepsilon)$ in the
local solution norm. The other heat sources are locally Lipschitz from $H^2$
to $H^1$, with a bound $C\sqrt T$ in $L^2(0,T;H^1)$. The transport bound has
the same factor by $\|\psi_\xi\|_{L^1(0,T;H^2)}\le\sqrt
T\|\psi\|_{L^2(0,T;H^3)}$. For $T\le1$, the homogeneous heat and transport
solutions have norm at most $C\|Z_0\|_{H^2}$. Thus
$C(\delta+\varepsilon+\sqrt T)\le\frac{1}{2}$ and
$\|Z_0\|_{H^2}\le\frac{\varepsilon}{2C}$ give a contraction on the ball and
a common local existence time.

The boundedness of $\ell$ on $L^2$ in \eqref{eq:feedback} and
$\|q(t,\cdot+x)-q(t,\cdot+y)\|_2\le |x-y|\|q_x(t)\|_2$ make the right-hand
side of \eqref{eq:moving-frame} uniformly Lipschitz in $X$ on the local
interval. The scalar equation therefore has a unique solution there. Since
$\bar q'\in H^3$, translation by $X(t)$ preserves the local Sobolev bounds
on the perturbation.

In the shifted frame, the mass equation in
\eqref{eq:energy-perturbation-system} gives, for $0\le s<t\le T$,
\[
 \begin{aligned}
 &\left\|\phi\left(t,\cdot-\int_s^t b(\tau)\dd\tau\right)
                 -\phi(s)\right\|_{H^2}\\
 &\quad\le \int_s^t
       \left(\|\psi_\xi(\tau)\|_{H^2}+|d(\tau)|\|V'\|_{H^2}\right)\dd\tau\\
 &\quad\le C_\delta\sqrt{t-s}\left(\int_s^t
             \left(\|\psi_\xi(\tau)\|_{H^2}^2+|d(\tau)|^2\right)
             \dd\tau\right)^{\frac{1}{2}} .
 \end{aligned}
\]
Apply the product rule in the other two equations, using the positive lower
bounds for $v$ and $\theta$, to obtain
\[
 \begin{aligned}
 \|(\psi_t,\zeta_t)\|_{L^2(s,t;H^1)}
 &\le C_\delta\left(\sqrt{t-s}\sup_{s\le\tau\le t}\|Z(\tau)\|_{H^2}
         +\|(\psi_\xi,\zeta_\xi)\|_{L^2(s,t;H^2)}
         +\|d\|_{L^2(s,t)}\right).
 \end{aligned}
\]
Translation is continuous in $H^2$. The $L^2(0,T;H^3)$ and time-derivative
bounds for $(\psi,\zeta)$ give continuity in $H^2$ by the Hilbert-space
energy identity. Thus $Z\in C([0,T];H^2)$ in the shifted frame.

The transport formula gives $D_t\phi_{\xi\xi}\in L^2(0,T;L^2)$ and
$(D_t\psi_{\xi\xi},D_t\zeta_{\xi\xi})\in L^2(0,T;H^{-1})$. The coefficients
in \eqref{eq:second-order-weight} have bounded spatial first derivatives and
material derivatives in $L^2(0,T;L^\infty)$ by
\eqref{eq:actual-material-equations}. Time regularization in the translated
frame gives the energy identity with the parabolic terms paired in $H^{-1}$
and $H^1$. In particular, \eqref{eq:pressure-weak-pairing} and
\eqref{eq:density-cross-identity} hold at this regularity and justify
\eqref{eq:global-H2}.

Take $\varepsilon_\delta$ with
$C_\delta^{\frac{1}{2}}\varepsilon_\delta\le\frac{\varepsilon}{4C}$, where
$C\ge1$ is the constant in the local bound on the data. On every interval
with $\sup_{0\le\tau\le t}\|Z(\tau)\|_{H^2}\le\varepsilon$,
\eqref{eq:global-H2} gives
\[
 \sup_{0\le\tau\le t}\|Z(\tau)\|_{H^2}
 \le C_\delta^{\frac{1}{2}}\|Z_0\|_{H^2}
 \le\frac{\varepsilon}{4C}
 <\frac{\varepsilon}{2C}.
\]
Continuity preserves this bound. Sobolev embedding then keeps $v$ and
$\theta$ uniformly positive. At any $s$, translate the data by $\sigma
s+X(s)$ to use $Z(s)$ in the same local construction. The common existence
time gives a global solution. Integrate the estimates from any $s\ge0$ to
get \eqref{eq:global-H2} with the same constant. Set $s=0$ and let $t$ tend
to infinity to obtain the integral bounds in \eqref{eq:global-regularity}.
\end{proof}

\subsection{An energy estimate for the time derivative}\label{sec:time}

Use the coordinate \eqref{eq:moving-frame}, the parameters $a$ and $M$ in
\eqref{eq:weight}, and the functional $\ell$ in \eqref{eq:feedback}. The
linear relation $d=\ell(Z)$ gives
\begin{equation}\label{eq:time-derivatives}
 F=Z_t=(W,Y,Q),\qquad e=d_t=\ell(F).
\end{equation}
With $B$ as in \eqref{eq:def-B}, \eqref{eq:anchored-weight} and
\eqref{eq:global-H2} give, for $0\le s<t$,
\begin{equation}\label{eq:time-B-integral}
 \begin{aligned}
 \int_s^t B(\rho)\dd\rho
 &\le C_\delta\int_s^t
    \left(\|\phi_{\xi\xi}\|_2^2+\mathcal B_1\right)\dd\rho
 \le C_\delta\|Z(s)\|_{H^2}^2 .
 \end{aligned}
\end{equation}
We write
\begin{equation}\label{eq:time-basic-quantities}
 S(t)=\|Z(t)\|_\infty,\qquad
 \Gamma(t)^2=\|F_\xi(t)\|_2^2+m_\omega(F(t))^2 .
\end{equation}

\begin{proposition}[Time-derivative energy]\label{prop:time}
Let $Z$ be the solution in Proposition~\ref{prop:global}. Assume
$\sup_{t\ge0}\|Z(t)\|_{H^2}\le\eps_\delta$ for sufficiently small
$\eps_\delta>0$. For some fixed $\tau>0$, set
\begin{equation}\label{eq:time-energy-definition}
 \begin{aligned}
 E_F&=\frac{1}{2}\int_{\R}a\left(
       \frac{R\theta}{v^2}W^2+Y^2+\frac{c_v}{\theta}Q^2
                               \right)\dd\xi,\\
 \mathcal T&=E_F+\tau\int_\R
             \left(\frac{\mu}{2v}W_\xi^2-YW_\xi\right)\dd\xi .
 \end{aligned}
\end{equation}
The function $\mathcal T$ is absolutely continuous on each finite time
interval and satisfies
\begin{equation}\label{eq:time-energy-equivalence}
 c_\delta\left(\|F\|_2^2+\|W_\xi\|_2^2\right)
 \le\mathcal T
 \le C_\delta\left(\|F\|_2^2+\|W_\xi\|_2^2\right)
 \le C_\delta B .
\end{equation}
For almost every $t>0$,
\begin{equation}\label{eq:time-energy-inequality}
 \mathcal T'+c_\delta\Gamma^2
 \le C_\delta B\mathcal T+C_\delta S^2B .
\end{equation}
For $I=[s,t]$ with $0\le s<t$, put
\begin{equation}\label{eq:time-window-notation}
 \begin{aligned}
 \beta_I^2&=\int_s^t B(\rho)\dd\rho,
 &Q_I^2&=\int_s^t\Gamma(\rho)^2\dd\rho,\\
 S_I&=\sup_{s\le\rho\le t}S(\rho),
 &\mathfrak b_s&=\sqrt{\mathcal T(s)} .
 \end{aligned}
\end{equation}
Then
\begin{equation}\label{eq:time-window-estimate}
 \begin{aligned}
 \sup_{s\le\rho\le t}\sqrt{\mathcal T(\rho)}+Q_I
 &\le C_\delta\left[\sqrt{\mathcal T(s)}+
       \left(\int_s^t S(\rho)^2 B(\rho)\dd\rho\right)^{\frac{1}{2}}\right]\\
 &\le C_\delta\bigl(\mathfrak b_s+S_I\beta_I\bigr).
 \end{aligned}
\end{equation}
The same constant applies to every finite interval $[s,t]$.
\end{proposition}

\begin{proof}
We derive the estimates for smooth solutions and justify them for the
solution of Proposition~\ref{prop:global} at the end. The subscripts $2$ and
$\infty$ denote norms over $\R$.

\medskip
\noindent\emph{The differentiated equations.}
The profile $\bar q$ is fixed in $\xi$, so $q_t=F$. With $b_t=e$ from
\eqref{eq:time-derivatives},
\begin{equation}\label{eq:time-chain-rule}
 \partial_t(D_t q)=D_tF-eq_\xi,\qquad
 p_t=-\frac{R\theta}{v^2}W+\frac{R}{v}Q .
\end{equation}
Differentiating \eqref{eq:energy-perturbation-system} in time and
using \eqref{eq:time-chain-rule} gives
\begin{equation}\label{eq:time-differentiated-system}
 \begin{aligned}
 D_tW-Y_\xi&=e v_\xi,\\
 D_tY+(p_t)_\xi&=\mu\left(\frac{Y_\xi}{v}
                           -\frac{u_\xi W}{v^2}\right)_\xi+e u_\xi,\\
 c_vD_tQ+pY_\xi+u_\xi p_t
 &=\kappa\left(\frac{Q_\xi}{v}-\frac{\theta_\xi W}{v^2}\right)_\xi
   +\mu\left(\frac{2u_\xi Y_\xi}{v}-\frac{u_\xi^2W}{v^2}\right)
   +c_ve\theta_\xi .
 \end{aligned}
\end{equation}
The bounds \eqref{eq:profilebounds} and \eqref{eq:feedback} and the identity
$\int_{\R} V'\dd\xi=\delta$ give, for $f\in L^2(\R;\R^3)$,
\begin{equation}\label{eq:time-phase-functional-bounds}
 \begin{aligned}
 |\ell(f)|
 &\le\frac{C_\delta}{\delta}\int_{\R}V'|f|\dd\xi
 \le\frac{C_\delta}{\delta}
     \left(\|V'\|_\infty\int_{\R} V'\dd\xi\right)^{\frac{1}{2}}\|f\|_2
 \le C_\delta\sqrt{\delta}\,\|f\|_2 .
 \end{aligned}
\end{equation}
Taking $f=Z$ and $f=F$ gives
$|d|\le C_\delta\sqrt{\delta}\,\|Z\|_2$ and
$|e|^2\le C_\delta\delta\|F\|_2^2$.
By Sobolev embedding and the $H^2$ bound,
\begin{equation}\label{eq:time-gradient-basic-bounds}
 \|Z_\xi\|_\infty^2\le C\|Z_\xi\|_{H^1}^2\le CB,\qquad
 \|\psi_{\xi\xi}\|_2^4
 \le\eps_\delta^2\|\psi_{\xi\xi}\|_2^2\le\eps_\delta^2 B .
\end{equation}

\medskip
\noindent\emph{The physical quadratic identity.}
Test the three equations in \eqref{eq:time-differentiated-system} with
$\frac{aR\theta}{v^2}W$, $aY$ and $\frac{aQ}{\theta}$, respectively. Since
$a$ is independent of $t$ and $b$ is independent of $\xi$,
\[
 \begin{aligned}
 &\int_{\R}a\left(
    \frac{R\theta}{v^2}W D_tW+Y D_tY+\frac{c_v}{\theta}Q D_tQ
         \right)\dd\xi\\
 &\quad=E_F'+\frac{b}{2}\int_{\R}a'
         \left(\frac{R\theta}{v^2}W^2+Y^2+\frac{c_v}{\theta}Q^2\right)\dd\xi\\
 &\qquad-\frac{1}{2}\int_{\R}a\left(
       W^2 D_t\left(\frac{R\theta}{v^2}\right)
       +Q^2 D_t\left(\frac{c_v}{\theta}\right)\right)\dd\xi .
 \end{aligned}
\]
Using \eqref{eq:time-chain-rule} and $\frac{p}{\theta}=\frac{R}{v}$,
we obtain
\begin{equation}\label{eq:time-pressure-flux}
 \begin{aligned}
 &\frac{R\theta}{v^2}WY_\xi-Y(p_t)_\xi-\frac{p}{\theta}QY_\xi
      =-p_tY_\xi-Y(p_t)_\xi=-(Yp_t)_\xi,\\
 &-\int_{\R}a(Yp_t)_\xi\dd\xi=\int_{\R}a'Yp_t\dd\xi .
 \end{aligned}
\end{equation}
By \eqref{eq:actual-material-equations} and \eqref{eq:time-chain-rule},
\begin{equation}\label{eq:time-euler-square}
 \begin{aligned}
 &\frac{1}{2}W^2D_t\left(\frac{R\theta}{v^2}\right)
  +\frac{1}{2}Q^2D_t\left(\frac{c_v}{\theta}\right)
  -\frac{u_\xi Qp_t}{\theta}\\
 &\quad-\frac{1}{2}\left(\frac{R}{c_v v^2}W^2-\frac{Q^2}{\theta^2}\right)
     \left[\kappa\left(\frac{\theta_\xi}{v}\right)_\xi
             +\frac{\mu u_\xi^2}{v}\right]\\
 &\quad=-\frac{R\theta u_\xi}{v^3}W^2
      -\frac{p u_\xi}{2}\left(\frac{R}{c_v v^2}W^2-\frac{Q^2}{\theta^2}\right)
      -\frac{u_\xi Qp_t}{\theta}\\
 &\quad=-\frac{Ru_\xi}{2v\theta}
     \left[(\gamma+1)\frac{\theta^2}{v^2}W^2
                  -\frac{2\theta}{v}WQ+Q^2\right]\\
 &\quad=-\frac{Ru_\xi}{2v\theta}
     \left[\left(Q-\frac{\theta}{v}W\right)^2
                  +\gamma\frac{\theta^2}{v^2}W^2\right] .
 \end{aligned}
\end{equation}
Integration by parts in the viscosity and heat conduction terms gives
\[
 -\mu\int_{\R}(aY)_\xi
         \left(\frac{Y_\xi}{v}-\frac{u_\xi W}{v^2}\right)\dd\xi
 -\kappa\int_{\R}\left(\frac{aQ}{\theta}\right)_\xi
         \left(\frac{Q_\xi}{v}-\frac{\theta_\xi W}{v^2}\right)\dd\xi .
\]
Adding the phase and viscous heating terms, we have
\begin{equation}\label{eq:time-exact-EF}
 \begin{aligned}
 &E_F'+\int_{\R}a\left(\frac{\mu Y_\xi^2}{v}
                     +\frac{\kappa Q_\xi^2}{v\theta}\right)\dd\xi\\
 &\quad=-\frac{b}{2}\int_{\R}a'
       \left(\frac{R\theta}{v^2}W^2+Y^2+\frac{c_v}{\theta}Q^2\right)\dd\xi
       +\int_{\R}a'Yp_t\dd\xi\\
 &\qquad-\int_{\R}a\frac{Ru_\xi}{2v\theta}
       \left[\left(Q-\frac{\theta}{v}W\right)^2
                        +\gamma\frac{\theta^2}{v^2}W^2\right]\dd\xi\\
 &\qquad+e\mathcal Y_q(F)+\mathcal V_q(F)+\mathcal W_q(F),
 \end{aligned}
\end{equation}
where
\begin{equation}\label{eq:time-Y-functional}
 \mathcal Y_q(F)=\int_{\R}a
       \left(\frac{R\theta}{v^2}v_\xi W+u_\xi Y
                      +\frac{c_v\theta_\xi}{\theta}Q\right)\dd\xi,
\end{equation}
\begin{equation}\label{eq:time-V-functional}
 \begin{aligned}
 \mathcal V_q(F)
 &=\int_{\R}a\left[\frac{\mu u_\xi}{v^2}WY_\xi
      +\frac{\kappa\theta_\xi}{v^2\theta}WQ_\xi
      +\frac{\mu Q}{\theta}
         \left(\frac{2u_\xi Y_\xi}{v}-\frac{u_\xi^2W}{v^2}\right)\right]\dd\xi\\
 &\quad-\mu\int_{\R}a'Y
                \left(\frac{Y_\xi}{v}-\frac{u_\xi W}{v^2}\right)\dd\xi
       -\kappa\int_{\R}\left(\frac{a}{\theta}\right)_\xi Q
                \left(\frac{Q_\xi}{v}-\frac{\theta_\xi W}{v^2}\right)\dd\xi,
 \end{aligned}
\end{equation}
and
\begin{equation}\label{eq:time-W-functional}
 \mathcal W_q(F)=\frac{1}{2}\int_{\R}a
   \left(\frac{R}{c_v v^2}W^2-\frac{Q^2}{\theta^2}\right)
   \left[\kappa\left(\frac{\theta_\xi}{v}\right)_\xi
                      +\frac{\mu u_\xi^2}{v}\right]\dd\xi .
\end{equation}

\medskip
\noindent\emph{Estimate at the viscous shock profile.}
At the profile, use the dissipation
\begin{equation}\label{eq:time-frozen-dissipation}
 D_F^0=\int_{\R}V'|F|^2\dd\xi
             +\|(Y_\xi,Q_\xi)\|_2^2+\delta e^2 .
\end{equation}
For some $c_0>0$, we have
\begin{equation}\label{eq:time-frozen-coercivity}
 \int_{\R}a\left(\frac{R\Theta}{V^2}W,Y,\frac{c_v}{\Theta}Q\right)
     \cdot\bigl[D_q\mathscr N_\sigma(\bar q)F+e\bar q'\bigr]\dd\xi
 \le-c_0D_F^0 .
\end{equation}

Apply Proposition~\ref{prop:global} to $Z^\rho(0)=\rho F$ with $F=(W,Y,Q)\in
C_c^\infty(\R;\R^3)$ and sufficiently small $\rho>0$. At $t=0$ the shift
velocity is $\rho e=\rho\ell(F)$. Taylor's formula in
\eqref{eq:st-spatial-operator}, with $\mathscr N_\sigma(\bar q)=0$, gives
\[
 \begin{aligned}
 &\left\|\frac{Z_t^\rho(0)}{\rho}
                   -D_q\mathscr N_\sigma(\bar q)F-e\bar q'\right\|_2\\
 &\quad\le
   \left\|\frac{\mathscr N_\sigma(\bar q+\rho F)-\mathscr N_\sigma(\bar q)}
                 {\rho}-D_q\mathscr N_\sigma(\bar q)F\right\|_2
                      +\rho|e|\|F_\xi\|_2
 \le C_\delta\rho\|F\|_{H^2}^2 .
 \end{aligned}
\]
The bound follows from the product rule in \eqref{eq:st-spatial-operator},
$H^2\hookrightarrow W^{1,\infty}$ and
\eqref{eq:time-phase-functional-bounds}. Since $Z^\rho$ is continuous in
$H^2$, \eqref{eq:energy-perturbation-system} gives continuity of $Z_t^\rho$
in $L^2$. The integrated inequality \eqref{eq:pure-entropy-differential}
therefore bounds the right derivative at $t=0$. Differentiating
\eqref{eq:pure-entropy-functional} gives
\[
 \left.\frac{1}{\rho^2}\frac{\mathrm d}{\mathrm dt}
                  \mathcal E_a(Z^\rho(t))\right|_{t=0}
 =\int_{\R}a\left(\frac{R\Theta W}{V(V+\rho W)},Y,
                    \frac{c_v Q}{\Theta+\rho Q}\right)
                  \cdot\frac{Z_t^\rho(0)}{\rho}\dd\xi .
\]
Divide the bound from \eqref{eq:pure-entropy-differential} at $t=0$ by
$\rho^2$ and use the preceding estimates to pass to the limit.
\begin{equation}\label{eq:time-entropy-limit}
 \begin{aligned}
 &\int_{\R}a\left(\frac{R\Theta}{V^2}W,Y,\frac{c_v}{\Theta}Q\right)
        \cdot\bigl[D_q\mathscr N_\sigma(\bar q)F+e\bar q'\bigr]\dd\xi\\
 &\quad=\lim_{\rho\downarrow0}
          \left.\frac{1}{\rho^2}\frac{\mathrm d}{\mathrm dt}
                    \mathcal E_a(Z^\rho(t))\right|_{t=0}\\
 &\quad\le -c\int_{\R}V'|F|^2\dd\xi
         -c\int_{\R}a\left(
                \frac{\mu Y_\xi^2}{V}+\frac{\kappa Q_\xi^2}{V\Theta}
                         \right)\dd\xi-\frac{\delta}{4M}e^2
 \le -c_0D_F^0 .
 \end{aligned}
\end{equation}
At the profile, \eqref{eq:profileODE} gives
$D_t(V,\Theta)=-\sigma(V',\Theta')$, in agreement with
\eqref{eq:actual-material-equations}. The pairing in
\eqref{eq:time-frozen-coercivity} is therefore the right-hand side of
\eqref{eq:time-exact-EF} evaluated at the profile, minus the diffusion
integral. By \eqref{eq:time-phase-functional-bounds} and the bounded profile
coefficients, this expression is continuous in
$\|F\|_2+\|(Y_\xi,Q_\xi)\|_2$. Density gives
\eqref{eq:time-frozen-coercivity} for $W\in L^2$ and $Y,Q\in H^1$.

By \eqref{eq:feedback} and \eqref{eq:time-Y-functional},
\begin{equation}\label{eq:time-exact-phase-sign}
 \begin{aligned}
 e\mathcal Y_{\bar q}(F)
 &=e\int_{\R}a\left(
        \frac{P}{V}V'W+U'Y+\frac{c_v\Theta'}{\Theta}Q\right)\dd\xi\\
 &=-\frac{\delta}{M}e\ell(F)=-\frac{\delta}{M}e^2 .
 \end{aligned}
\end{equation}

\medskip
\noindent\emph{Comparison with the current coefficients.}
For each fixed $\eta>0$, choose $\eps_\delta$ small enough to obtain
\begin{equation}\label{eq:time-EF-estimate}
 E_F'+\frac{c_0}{2}D_F^0
 \le\eta\|W_\xi\|_2^2+C_{\eta,\delta}(B+S^2)\|F\|_2^2.
\end{equation}
By \eqref{eq:time-frozen-coercivity},
\[
 E_F'+c_0D_F^0\le E_F'-
 \int_{\R}a\left(\frac{R\Theta}{V^2}W,Y,\frac{c_v}{\Theta}Q\right)
       \cdot\bigl[D_q\mathscr N_\sigma(\bar q)F+e\bar q'\bigr]\dd\xi .
\]
By \eqref{eq:time-exact-EF}, the right-hand side contains the differences
between the current and profile coefficients, evaluated on the same $F$. The
ratios in \eqref{eq:time-exact-EF} and \eqref{eq:time-Y-functional} are
Lipschitz on the positive state range. The bounds \eqref{eq:profilebounds}
and $\int_{\R}V'\dd\xi=\delta$ give
\begin{align}
 &|e(\mathcal Y_q-\mathcal Y_{\bar q})(F)|\notag\\*
 &\quad=|e|\Bigg|\int_{\R}a\Bigg[
      \frac{R\theta}{v^2}\phi_\xi W+\psi_\xi Y
                   +\frac{c_v}{\theta}\zeta_\xi Q\notag\\*
 &\qquad\qquad+\left(\frac{R\theta}{v^2}-\frac{R\Theta}{V^2}\right)V'W
       +\left(\frac{c_v}{\theta}-\frac{c_v}{\Theta}\right)\Theta'Q
                                               \Bigg]\dd\xi\Bigg|\notag\\
 &\quad\le C|e|\left[\int_\R|Z_\xi||F|\dd\xi+
                              S\int_\R V'|F|\dd\xi\right]\notag\\
 &\quad\le C|e|\left[\|Z_\xi\|_2\|F\|_2+
                 S\sqrt{\delta}\left(\int_\R V'|F|^2\dd\xi\right)^{\frac{1}{2}}\right]\notag\\
 &\quad\le\eta\delta e^2+C_{\eta,\delta}
             \left[B\|F\|_2^2+S^2\int_\R V'|F|^2\dd\xi\right].
 \label{eq:time-phase-error}
\end{align}
For the two terms with $a'$ in \eqref{eq:time-exact-EF}, use $b=\sigma+d$
and \eqref{eq:time-chain-rule}. With $a'=\delta^{-\frac{1}{4}}V'$ and
$|d|\le C_\delta S$ from \eqref{eq:feedbackbounds}, we have
\begin{equation}\label{eq:time-moving-weight-error}
 \begin{aligned}
 &\Bigg|\int_{\R}a'\Bigg[
     -\frac{d}{2}\left(\frac{R\theta}{v^2}W^2+Y^2
                                +\frac{c_v}{\theta}Q^2\right)\\
 &\qquad-\frac{\sigma}{2}\left[
      \left(\frac{R\theta}{v^2}-\frac{R\Theta}{V^2}\right)W^2
      +\left(\frac{c_v}{\theta}-\frac{c_v}{\Theta}\right)Q^2\right]\\
 &\qquad+Y\left(p_t+\frac{R\Theta}{V^2}W-\frac{R}{V}Q\right)
                                       \Bigg]\dd\xi\Bigg|\\
 &\quad\le C\delta^{-\frac{1}{4}}(|d|+S)\int_\R V'|F|^2\dd\xi
 \le C_\delta S\int_\R V'|F|^2\dd\xi .
 \end{aligned}
\end{equation}
In \eqref{eq:time-euler-square}, substitute $u_\xi=U'+\psi_\xi$, expand the
square and integrate the term with $\psi_\xi$ by parts.
\begin{equation}\label{eq:time-compression-parts}
 \begin{aligned}
 &-\frac{R}{2}\int_{\R}a\psi_\xi
      \left(\frac{(\gamma+1)\theta}{v^3}W^2
               -\frac{2}{v^2}WQ+\frac{1}{v\theta}Q^2\right)\dd\xi\\
 &\quad=\frac{R}{2}\int_{\R}\psi\left[
       (\gamma+1)\left(\frac{a\theta}{v^3}\right)_\xi W^2
       -2\left(\frac{a}{v^2}\right)_\xi WQ
       +\left(\frac{a}{v\theta}\right)_\xi Q^2\right]\dd\xi\\
 &\qquad+R\int_{\R}a\psi\left[
       \frac{(\gamma+1)\theta}{v^3}WW_\xi
       -\frac{W_\xi Q+WQ_\xi}{v^2}+\frac{QQ_\xi}{v\theta}\right]\dd\xi .
 \end{aligned}
\end{equation}
By \eqref{eq:profilebounds} and \eqref{eq:weight}, $|a'|+|\bar q'|\le
C_\delta V'$. The coefficient derivatives in
\eqref{eq:time-compression-parts} are bounded by $C(|a'|+|\bar
q'|+|Z_\xi|)$. With the Lipschitz bounds on the ratios, Young's inequality
and \eqref{eq:time-gradient-basic-bounds}, we obtain
\begin{align}
 &\Bigg|-\frac{R}{2}\int_{\R}a\Bigg[
 (\gamma+1)\left(\frac{\theta u_\xi}{v^3}
                      -\frac{\Theta U'}{V^3}\right)W^2\notag\\*
 &\qquad-2\left(\frac{u_\xi}{v^2}-\frac{U'}{V^2}\right)WQ
       +\left(\frac{u_\xi}{v\theta}-\frac{U'}{V\Theta}\right)Q^2
                                   \Bigg]\dd\xi\Bigg|\notag\\
 &\quad\le CS\int_{\R}(|a'|+|\bar q'|)|F|^2\dd\xi
       +CS\int_{\R}(|Z_\xi||F|^2+|F||F_\xi|)\dd\xi\notag\\
 &\quad\le C_\delta S\int_\R V'|F|^2\dd\xi
      +CS\|Z_\xi\|_\infty\|F\|_2^2+CS\|F\|_2\|F_\xi\|_2\notag\\*
 &\quad\le C_\delta S\int_\R V'|F|^2\dd\xi+\eta\|F_\xi\|_2^2
                 +C_{\eta,\delta}(B+S^2)\|F\|_2^2 .
 \label{eq:time-compression-bounds}
\end{align}

Subtract the profile contributions from the diffusion terms in
\eqref{eq:time-exact-EF} and the terms in \eqref{eq:time-V-functional}. For
the fluxes in \eqref{eq:time-differentiated-system}, the difference is
\[
 \begin{aligned}
 &\frac{(Y_\xi,Q_\xi)}{v}-\frac{W(u_\xi,\theta_\xi)}{v^2}
      -\frac{(Y_\xi,Q_\xi)}{V}+\frac{W(U',\Theta')}{V^2}\\
 &\quad=\left(\frac{1}{v}-\frac{1}{V}\right)(Y_\xi,Q_\xi)
      -\frac{W}{v^2}(\psi_\xi,\zeta_\xi)
      -W\left(\frac{1}{v^2}-\frac{1}{V^2}\right)(U',\Theta').
 \end{aligned}
\]
For the differentiated weight and the squared derivatives, we use
\[
 \begin{aligned}
 \left(\frac{a}{\theta}\right)_\xi-\left(\frac{a}{\Theta}\right)_\xi
 &=a'\left(\frac{1}{\theta}-\frac{1}{\Theta}\right)
     -\frac{a\zeta_\xi}{\theta^2}
     -a\left(\frac{1}{\theta^2}-\frac{1}{\Theta^2}\right)\Theta',\\
 u_\xi^2-(U')^2&=2U'\psi_\xi+\psi_\xi^2,\qquad
 \theta_\xi^2-(\Theta')^2=2\Theta'\zeta_\xi+\zeta_\xi^2 .
 \end{aligned}
\]
The terms contain $Y_\xi$ and $Q_\xi$ as derivatives of $F$. Use the
preceding identities, $|a'|+|\bar q'|\le C_\delta V'$ and
\eqref{eq:time-gradient-basic-bounds} to obtain
\begin{align}
 &\Bigg|\mathcal V_q(F)-\mathcal V_{\bar q}(F)
  -\int_{\R}a\left[
     \mu\left(\frac{1}{v}-\frac{1}{V}\right)Y_\xi^2
    +\kappa\left(\frac{1}{v\theta}-\frac{1}{V\Theta}\right)Q_\xi^2
                                   \right]\dd\xi\Bigg|\notag\\*
 &\quad\le CS\|(Y_\xi,Q_\xi)\|_2^2\notag\\*
 &\qquad+C\int_\R
       \left[|Z_\xi|+S(|a'|+|\bar q'|)\right]|F|\,|(Y_\xi,Q_\xi)|\dd\xi\notag\\*
 &\qquad+C\int_\R\left[
       |Z_\xi|^2+(|a'|+|\bar q'|)|Z_\xi|
             +S(|a'|+|\bar q'|)^2\right]|F|^2\dd\xi\notag\\
 &\quad\le CS\|(Y_\xi,Q_\xi)\|_2^2
       +C\|Z_\xi\|_\infty\|F\|_2\|(Y_\xi,Q_\xi)\|_2\notag\\*
 &\qquad+C_\delta S\left(\int_\R V'|F|^2\dd\xi\right)^{\frac{1}{2}}
                                  \|(Y_\xi,Q_\xi)\|_2\notag\\*
 &\qquad+C_\delta(S+\|Z_\xi\|_\infty)\int_\R V'|F|^2\dd\xi
                          +C\|Z_\xi\|_\infty^2\|F\|_2^2\notag\\
 &\quad\le (CS+\eta)\|(Y_\xi,Q_\xi)\|_2^2\notag\\*
 &\qquad+C_{\eta,\delta}(S+S^2+\|Z_\xi\|_\infty)\int_\R V'|F|^2\dd\xi
                        +C_{\eta,\delta}B\|F\|_2^2 .
 \label{eq:time-gradient-coefficient-error}
\end{align}
By \eqref{eq:profilebounds} and \eqref{eq:weight},
\begin{equation}\label{eq:time-profile-products}
 |\bar q'|^2+|a'|^2+|\bar q''|
 \le C_\delta\left((V')^2+\delta V'\right)
 \le C_\delta V'.
\end{equation}

In \eqref{eq:time-W-functional}, subtract the value at the profile and
integrate the heat-flux difference by parts. The profile equation
\eqref{eq:profileODE} gives
\begin{align*}
 &\mathcal W_q(F)-\mathcal W_{\bar q}(F)\\*
 &\quad=-\frac{\kappa}{2}\int_{\R}
 \left(\frac{\zeta_\xi}{v}-\frac{\Theta'\phi}{vV}\right)
 \Bigg[\frac{R}{c_v}\left(\frac{a}{v^2}\right)_\xi W^2
                          -\left(\frac{a}{\theta^2}\right)_\xi Q^2\\*
 &\hspace{11em}+2a\left(\frac{RWW_\xi}{c_v v^2}
                                   -\frac{QQ_\xi}{\theta^2}\right)\Bigg]\dd\xi\\
 &\qquad+\frac{1}{2}\int_{\R}a\left[
   \frac{R}{c_v}\left(\frac{1}{v^2}-\frac{1}{V^2}\right)W^2
    -\left(\frac{1}{\theta^2}-\frac{1}{\Theta^2}\right)Q^2\right]
                 (PU'-\sigma c_v\Theta')\dd\xi\\*
 &\qquad+\frac{\mu}{2}\int_{\R}a\left(\frac{R W^2}{c_v v^2}
      -\frac{Q^2}{\theta^2}\right)
   \left(\frac{2U'\psi_\xi+\psi_\xi^2}{v}
                 -\frac{(U')^2\phi}{vV}\right)\dd\xi .
\end{align*}
The coefficient derivatives in the first integral are bounded by
$C(|a'|+|\bar q'|+|Z_\xi|)$. By \eqref{eq:profilebounds}, $|PU'-\sigma
c_v\Theta'|\le C_\delta V'$. Since $S\le1$, apply
\eqref{eq:time-profile-products}, \eqref{eq:time-gradient-basic-bounds} and
Young's inequality to obtain
\[
 \begin{aligned}
 &|\mathcal W_q(F)-\mathcal W_{\bar q}(F)|\\
 &\quad\le C\int_{\R} (|Z_\xi|+S|\bar q'|)|F||F_\xi|\dd\xi
    +C_\delta\int_{\R}\left[|Z_\xi|^2+(S+|Z_\xi|)V'\right]|F|^2\dd\xi\\
 &\quad\le C\|Z_\xi\|_\infty\|F\|_2\|F_\xi\|_2
    +C_\delta S \left(\int_\R V'|F|^2\dd\xi\right)^{\frac{1}{2}}\|F_\xi\|_2\\
 &\qquad+C_\delta(S+\|Z_\xi\|_\infty)\int_\R V'|F|^2\dd\xi
                      +C_\delta\|Z_\xi\|_\infty^2\|F\|_2^2\\
 &\quad\le\eta\|F_\xi\|_2^2
      +C_{\eta,\delta}(S+S^2+\|Z_\xi\|_\infty)\int_\R V'|F|^2\dd\xi
      +C_{\eta,\delta}B\|F\|_2^2 .
 \end{aligned}
\]
Use $\frac{c_0}{8}$ as the Young parameter in \eqref{eq:time-phase-error}
and \eqref{eq:time-gradient-coefficient-error}. Use
$\min\{\frac{\eta}{2},\frac{c_0}{16}\}$ in
\eqref{eq:time-compression-bounds} and in the heat coefficient bound. Since
$S+\|Z_\xi\|_\infty\le C\eps_\delta$, choose $\eps_\delta$ with
$C_{\eta,\delta}(S+S^2+\|Z_\xi\|_\infty)\le\frac{c_0}{8}$. Add the bounds
and use \eqref{eq:time-frozen-dissipation} and
$\|F_\xi\|_2^2=\|W_\xi\|_2^2+\|(Y_\xi,Q_\xi)\|_2^2$ to obtain
\[
 \begin{aligned}
 E_F'+c_0D_F^0
 &\le\frac{c_0}{8}\left(\delta e^2+\|(Y_\xi,Q_\xi)\|_2^2\right)
      +\min\left\{\eta,\frac{c_0}{8}\right\}\|F_\xi\|_2^2\\
 &\quad+C_{\eta,\delta}(S+S^2+\|Z_\xi\|_\infty)
             \left(\int_\R V'|F|^2\dd\xi+\|(Y_\xi,Q_\xi)\|_2^2\right)\\
 &\quad+C_{\eta,\delta}(B+S^2)\|F\|_2^2\\
 &\le\frac{3c_0}{8}D_F^0+\eta\|W_\xi\|_2^2+
                        C_{\eta,\delta}(B+S^2)\|F\|_2^2 .
 \end{aligned}
\]
This proves \eqref{eq:time-EF-estimate}.

\medskip
\noindent\emph{Gradient estimate for $W$.}
\eqref{eq:actual-material-equations} and
\eqref{eq:time-differentiated-system} give $D_tv=u_\xi$ and
$D_tW_\xi=Y_{\xi\xi}+e v_{\xi\xi}$. By the product rule,
\[
 \begin{aligned}
 D_t\left(\frac{\mu W_\xi}{v}-Y\right)
 &=\frac{\mu}{v}(Y_{\xi\xi}+e v_{\xi\xi})
       -\frac{\mu u_\xi}{v^2}W_\xi-D_tY\\
 &=(p_t)_\xi+\frac{\mu v_\xi}{v^2}Y_\xi
       +\mu\left(\frac{u_\xi}{v^2}\right)_\xi W
       +e\left(\frac{\mu v_{\xi\xi}}{v}-u_\xi\right).
 \end{aligned}
\]
The spatial integral of the transport term vanishes because $b$ is
independent of $\xi$. Integration by parts gives
\[
 \begin{aligned}
 &\frac{\mathrm d}{\mathrm dt}\int_\R
           \left(\frac{\mu}{2v}W_\xi^2-YW_\xi\right)\dd\xi\\
 &\quad=\int_\R\left[
 W_\xi D_t\left(\frac{\mu W_\xi}{v}-Y\right)
 -Y(Y_{\xi\xi}+e v_{\xi\xi})
 +\frac{\mu u_\xi}{2v^2}W_\xi^2\right]\dd\xi\\
 &\quad=\int_\R W_\xi D_t\left(\frac{\mu W_\xi}{v}-Y\right)\dd\xi
 +\|Y_\xi\|_2^2-e\int_\R Yv_{\xi\xi}\dd\xi
 +\frac{\mu}{2}\int_\R\frac{u_\xi}{v^2}W_\xi^2\dd\xi .
 \end{aligned}
\]
Expanding $(p_t)_\xi$ by \eqref{eq:time-chain-rule} yields
\begin{equation}\label{eq:time-exact-KF}
 \begin{aligned}
 &\frac{\mathrm d}{\mathrm dt}\int_\R
       \left(\frac{\mu}{2v}W_\xi^2-YW_\xi\right)\dd\xi
    +\int_\R\frac{R\theta}{v^2}W_\xi^2\dd\xi\\
 &\quad=\|Y_\xi\|_2^2+\int_\R\frac{R}{v}Q_\xi W_\xi\dd\xi\\
 &\qquad+\int_\R\left[-\left(\frac{R\theta}{v^2}\right)_\xi W
       +\left(\frac{R}{v}\right)_\xi Q
       +\frac{\mu v_\xi}{v^2}Y_\xi
       +\mu\left(\frac{u_\xi}{v^2}\right)_\xi W\right]
       W_\xi\dd\xi\\
 &\qquad+\frac{\mu}{2}\int_\R\frac{u_\xi}{v^2}W_\xi^2\dd\xi
       +e\int_\R\left[
       \left(\frac{\mu v_{\xi\xi}}{v}-u_\xi\right)W_\xi
                         -Yv_{\xi\xi}\right]\dd\xi .
 \end{aligned}
\end{equation}
On the positive state range, fix $c_\delta>0$ with $\frac{R\theta}{v^2}\ge
c_\delta$. The derivatives of $\frac{R\theta}{v^2}$ and $\frac{R}{v}$ are
bounded by $C(|\bar q'|+|Z_\xi|)$ and $\|v_\xi\|_\infty\le C_\delta$. Apply
\eqref{eq:time-gradient-basic-bounds}, \eqref{eq:time-profile-products} and
Young's inequality to obtain
\begin{align*}
 &\left|\int_{\R}\left[\frac{R}{v}Q_\xi
       -\left(\frac{R\theta}{v^2}\right)_\xi W
       +\left(\frac{R}{v}\right)_\xi Q
       +\frac{\mu v_\xi}{v^2}Y_\xi\right]W_\xi\dd\xi\right|\\*
 &\quad\le C_\delta\|(Y_\xi,Q_\xi)\|_2\|W_\xi\|_2
       +C\int_{\R}(|\bar q'|+|Z_\xi|)|F||W_\xi|\dd\xi\\
 &\quad\le C_\delta\left(\|(Y_\xi,Q_\xi)\|_2+\left(\int_\R V'|F|^2\dd\xi\right)^{\frac{1}{2}}
                  +\|Z_\xi\|_\infty\|F\|_2\right)\|W_\xi\|_2\\
 &\quad\le\eta\|W_\xi\|_2^2
       +C_{\eta,\delta}\left(\int_\R V'|F|^2\dd\xi+\|(Y_\xi,Q_\xi)\|_2^2\right)
       +C_{\eta,\delta}B\|F\|_2^2 .
\end{align*}

Expanding the viscous coefficient in \eqref{eq:time-exact-KF} gives
\[
 \left(\frac{u_\xi}{v^2}\right)_\xi
 =\frac{\psi_{\xi\xi}}{v^2}+\frac{U''}{v^2}
   -\frac{2}{v^3}\left(U'V'+U'\phi_\xi+V'\psi_\xi
                                    +\psi_\xi\phi_\xi\right).
\]
Since $\|Z_\xi\|_\infty\le C\eps_\delta$, \eqref{eq:profilebounds} and
\eqref{eq:time-profile-products} bound the terms with $U''$, $U'$ or $V'$ by
$C_\delta V'$. Integration of $(W^2)_\xi$ gives $\|W\|_\infty^2\le
2\|W\|_2\|W_\xi\|_2$. Also, \eqref{eq:time-gradient-basic-bounds} gives
$\|Z_\xi\|_\infty^4\le C\eps_\delta^2 B$. Apply Young's inequality with
exponents $\frac{4}{3}$ and $4$ to the term with $\|W_\xi\|_2^{\frac{3}{2}}$
to obtain
\begin{align}
 &\left|\mu\int_{\R}
       \left(\frac{u_\xi}{v^2}\right)_\xi WW_\xi\dd\xi\right|\notag\\*
 &\quad\le C_\delta\int_{\R}V'|W||W_\xi|\dd\xi
    +C\|\psi_{\xi\xi}\|_2\|W\|_\infty\|W_\xi\|_2
    +C\|Z_\xi\|_\infty^2\|W\|_2\|W_\xi\|_2\notag\\
 &\quad\le C_\delta\left(\left(\int_\R V'|F|^2\dd\xi\right)^{\frac{1}{2}}
             +\|Z_\xi\|_\infty^2\|W\|_2\right)\|W_\xi\|_2
     +C\|\psi_{\xi\xi}\|_2\|W\|_2^{\frac{1}{2}}
                            \|W_\xi\|_2^{\frac{3}{2}}\notag\\
 &\quad\le\eta\|W_\xi\|_2^2+C_{\eta,\delta}\int_\R V'|F|^2\dd\xi
    +C_{\eta,\delta}\left(\|\psi_{\xi\xi}\|_2^4
                         +\|Z_\xi\|_\infty^4\right)\|W\|_2^2\notag\\
 &\quad\le\eta\|W_\xi\|_2^2+C_{\eta,\delta}\int_\R V'|F|^2\dd\xi
                    +C_{\eta,\delta}\eps_\delta^2 B\|F\|_2^2 .
 \label{eq:time-critical-KF-product}
\end{align}
Since $U'=-\sigma V'<0$, the next term in \eqref{eq:time-exact-KF} satisfies
\[
 \frac{\mu}{2}\int_{\R}\frac{u_\xi}{v^2}W_\xi^2\dd\xi
 \le\frac{\mu}{2}\int_{\R}\frac{\psi_\xi}{v^2}W_\xi^2\dd\xi
 \le C\|\psi_\xi\|_\infty\|W_\xi\|_2^2
 \le C\eps_\delta\|W_\xi\|_2^2 .
\]

For the phase terms in \eqref{eq:time-exact-KF}, \eqref{eq:profilebounds}
gives
\[
 \|(V',V'',U')\|_2\le C_\delta\|V'\|_2
 \le C_\delta\left(\|V'\|_\infty\int_\R V'\dd\xi\right)^{\frac{1}{2}}
 \le C_\delta\delta^{\frac{3}{2}} .
\]
Using $v_{\xi\xi}=V''+\phi_{\xi\xi}$, $u_\xi=U'+\psi_\xi$,
\eqref{eq:time-phase-functional-bounds} and \eqref{eq:def-B}, we have
\begin{equation}\label{eq:time-KF-phase-first}
 \begin{aligned}
 &|e|\left|\int_{\R}
       \left(\frac{\mu v_{\xi\xi}}{v}-u_\xi\right)W_\xi\dd\xi\right|\\
 &\quad\le C_\delta|e|\left(\|(U',V'')\|_2+\|Z_\xi\|_{H^1}\right)
                  \|W_\xi\|_2
 \le C_\delta|e|\left(\delta^{\frac{3}{2}}+\|Z_\xi\|_{H^1}\right)\|W_\xi\|_2\\
 &\quad\le\eta\|W_\xi\|_2^2+C_{\eta,\delta}
                 \left(\delta^3 e^2+|e|^2\|Z_\xi\|_{H^1}^2\right)\\
 &\quad\le\eta\|W_\xi\|_2^2+
            C_{\eta,\delta}\left(\delta^3 e^2+B\|F\|_2^2\right).
 \end{aligned}
\end{equation}
Integrating the remaining phase term by parts gives
\begin{equation}\label{eq:time-KF-phase-last}
 \begin{aligned}
 \left|-e\int_{\R}Yv_{\xi\xi}\dd\xi\right|
 &=\left|e\int_{\R}Y_\xi(V'+\phi_\xi)\dd\xi\right|\\
 &\le |e|\|Y_\xi\|_2\bigl(\|V'\|_2+\|\phi_\xi\|_2\bigr)\\
 &\le C\|Y_\xi\|_2^2+C_\delta\left(\delta^3 e^2+B\|F\|_2^2\right).
 \end{aligned}
\end{equation}

Use $\frac{c_\delta}{8}$ as the Young parameter in the three bounds with
$\eta\|W_\xi\|_2^2$. If $C\eps_\delta\le\frac{c_\delta}{8}$, the coefficient
of $\|W_\xi\|_2^2$ on the right of \eqref{eq:time-exact-KF} is at most
$\frac{3c_\delta}{8}+C\eps_\delta \le\frac{c_\delta}{2}$. Use
$\delta^3\le\delta$ and $c_1=\frac{c_\delta}{2}$ to obtain
\begin{equation}\label{eq:time-KF-estimate}
 \begin{aligned}
 &\frac{\mathrm d}{\mathrm dt}\int_\R
       \left(\frac{\mu}{2v}W_\xi^2-YW_\xi\right)\dd\xi
       +c_1\|W_\xi\|_2^2\\
 &\quad\le C_1\left(\int_\R V'|F|^2\dd\xi+
                      \|(Y_\xi,Q_\xi)\|_2^2+\delta e^2\right)
          +C_\delta B\|F\|_2^2
 =C_1D_F^0+C_\delta B\|F\|_2^2 .
 \end{aligned}
\end{equation}

\medskip
\noindent\emph{Combination of the estimates.}
Fix $C_1$ as in \eqref{eq:time-KF-estimate}. Choose $\tau>0$ with $\tau
C_1\le\frac{c_0}{4}$ and $\frac{\tau v}{\mu}\le\frac{a}{4}$, which is
possible since $v$ is bounded above and $a\ge1$. Take $0<\eta\le\frac{\tau
c_1}{2}$ in \eqref{eq:time-EF-estimate} and shrink $\eps_\delta$ for that
estimate. Add \eqref{eq:time-EF-estimate} and $\tau$ times
\eqref{eq:time-KF-estimate} to obtain
\begin{equation}\label{eq:time-combined-before-B}
 \begin{aligned}
 &\mathcal T'+\frac{c_0}{4}D_F^0+\frac{\tau c_1}{2}\|W_\xi\|_2^2\\
 &\quad\le\mathcal T'+\left(\frac{c_0}{2}-\tau C_1\right)D_F^0
                       +(\tau c_1-\eta)\|W_\xi\|_2^2\\
 &\quad\le C_{\eta,\delta}(B+S^2)\|F\|_2^2+\tau C_\delta B\|F\|_2^2
       \le C_\delta(B+S^2)\|F\|_2^2 .
 \end{aligned}
\end{equation}
Young's inequality gives
$|YW_\xi|\le\frac{\mu}{4v}W_\xi^2+\frac{v}{\mu}Y^2$. Apply $\frac{\tau
v}{\mu}\le\frac{a}{4}$ in \eqref{eq:time-energy-definition} to obtain
\[
 \begin{aligned}
 c_\delta\bigl(\|F\|_2^2+\|W_\xi\|_2^2\bigr)
 &\le\frac{1}{2}E_F+\tau\int_{\R}\frac{\mu}{4v}W_\xi^2\dd\xi\\
 &\le E_F+\tau\int_{\R}
      \left(\frac{\mu}{4v}W_\xi^2-\frac{v}{\mu}Y^2\right)\dd\xi
       \le\mathcal T,\\
 \mathcal T
 &\le E_F+\tau\int_{\R}
      \left(\frac{3\mu}{4v}W_\xi^2+\frac{v}{\mu}Y^2\right)\dd\xi\\
 &\le C_\delta\bigl(\|F\|_2^2+\|W_\xi\|_2^2\bigr).
 \end{aligned}
\]
This proves the first two inequalities in \eqref{eq:time-energy-equivalence}.

By \eqref{eq:anchored-weight} and \eqref{eq:time-frozen-dissipation},
\begin{equation}\label{eq:time-local-anchor}
 \begin{aligned}
 \Gamma^2
 &=\|F_\xi\|_2^2+m_\omega(F)^2
 \le C_\delta\left(\int_\R V'|F|^2\dd\xi+\|F_\xi\|_2^2\right)\\
 &\le C_\delta\left(D_F^0+\|W_\xi\|_2^2\right)
 \le C_\delta\left(\frac{c_0}{4}D_F^0+
                          \frac{\tau c_1}{2}\|W_\xi\|_2^2\right).
 \end{aligned}
\end{equation}

The mass equation in \eqref{eq:energy-perturbation-system},
$\pi=\frac{R\zeta}{v}-\frac{R\Theta\phi}{vV}$ and
\eqref{eq:energy-fluxes} give
\[
 \begin{aligned}
 W_\xi&=b\phi_{\xi\xi}+\psi_{\xi\xi}+dV'',\\
 \pi_\xi&=\frac{R\zeta_\xi}{v}-\frac{R\Theta\phi_\xi}{vV}
   -\frac{Rv_\xi\zeta}{v^2}-\left(\frac{R\Theta}{vV}\right)_\xi\phi,\\
 (J_v,J_\theta)_\xi
 &=\frac{(\psi_{\xi\xi},\zeta_{\xi\xi})}{v}
    -\frac{v_\xi(\psi_\xi,\zeta_\xi)}{v^2}
    -\frac{\phi_\xi(U',\Theta')}{vV}
    -\left(\frac{(U',\Theta')}{vV}\right)_\xi\phi,\\
 \frac{u_\xi^2}{v}-\frac{(U')^2}{V}
 &=\frac{2U'\psi_\xi+\psi_\xi^2}{v}-\frac{(U')^2\phi}{vV}.
 \end{aligned}
\]
By \eqref{eq:profilebounds} and the choice of $\omega$, $|\bar q'|+|\bar
q''|\le C_\delta\sqrt{\omega}$. Sobolev embedding and
\eqref{eq:feedbackbounds} bound $\|Z\|_\infty+\|Z_\xi\|_\infty+|b|$
uniformly. Use $v_\xi=V'+\phi_\xi$ in the preceding formulas and take norms
in \eqref{eq:energy-perturbation-system} to obtain
\begin{equation}\label{eq:time-F-controlled-by-B}
 \begin{aligned}
 &\|F\|_2+\|W_\xi\|_2\\
 &\quad\le C_\delta\left[
     \|Z_\xi\|_{H^1}+|d|
     +\left\|\left(\pi_\xi,\pi U',(J_v)_\xi,(J_\theta)_\xi,
                      \frac{u_\xi^2}{v}-\frac{(U')^2}{V}\right)\right\|_2\right]\\
 &\quad\le C_\delta\left[
    \|Z_\xi\|_{H^1}+|d|+(\|Z\|_\infty+\|Z_\xi\|_\infty)\|Z_\xi\|_2
                      +\big\|(|\bar q'|+|\bar q''|)|Z|\big\|_2\right]\\
 &\quad\le C_\delta\left(\|Z_\xi\|_{H^1}+m_\omega(Z)+|d|\right)
 \le C_\delta\sqrt{B}.
 \end{aligned}
\end{equation}
This completes \eqref{eq:time-energy-equivalence}.
Using \eqref{eq:time-local-anchor} in \eqref{eq:time-combined-before-B},
we have
\[
 \begin{aligned}
 \mathcal T'+c_\delta\Gamma^2
 &\le C_\delta B\|F\|_2^2+C_\delta S^2\|F\|_2^2\\
 &\le C_\delta B\mathcal T+C_\delta S^2B .
 \end{aligned}
\]
\eqref{eq:time-energy-equivalence} and \eqref{eq:time-F-controlled-by-B}
bound the two terms on the last line, respectively, giving
\eqref{eq:time-energy-inequality}.

Use the coefficient $C_\delta$ of $B\mathcal T$ in
\eqref{eq:time-energy-inequality}. Multiply that inequality by
$\exp(-C_\delta\int_s^\rho B(\lambda)\dd\lambda)$ and integrate. Since
$B,\mathcal T\ge0$, for $s\le\rho\le t$,
\[
 \begin{aligned}
 &\exp\left(-C_\delta\int_s^t B(\lambda)\dd\lambda\right)
    \left[\mathcal T(\rho)+c_\delta\int_s^\rho
                                    \Gamma(\lambda)^2\dd\lambda\right]\\
 &\quad\le
   \exp\left(-C_\delta\int_s^\rho B(\lambda)\dd\lambda\right)
       \mathcal T(\rho)
   +c_\delta\int_s^\rho
      \exp\left(-C_\delta\int_s^\lambda B(r)\dd r\right)
                                     \Gamma(\lambda)^2\dd\lambda\\
 &\quad\le \mathcal T(s)+C_\delta\int_s^\rho
     \exp\left(-C_\delta\int_s^\lambda B(r)\dd r\right)
                                     S(\lambda)^2 B(\lambda)\dd\lambda\\
 &\quad\le \mathcal T(s)+C_\delta\int_s^t
                                     S(\lambda)^2 B(\lambda)\dd\lambda .
 \end{aligned}
\]
Take the supremum in the energy term and set $\rho=t$ in the dissipation
term. Then \eqref{eq:time-B-integral} gives
\begin{equation}\label{eq:time-window-squared}
 \begin{aligned}
 &\sup_{s\le\rho\le t}\mathcal T(\rho)+\int_s^t\Gamma(\rho)^2\dd\rho\\
 &\quad\le C_\delta\exp\left(C_\delta\int_s^t B(\lambda)\dd\lambda\right)
             \left[\mathcal T(s)+\int_s^t S(\lambda)^2B(\lambda)\dd\lambda\right]\\
 &\quad\le C_\delta\exp\bigl(C_\delta\|Z(s)\|_{H^2}^2\bigr)
             \left[\mathcal T(s)+\int_s^t S(\lambda)^2B(\lambda)\dd\lambda\right]\\
 &\quad\le C_\delta
             \left[\mathcal T(s)+\int_s^t S(\lambda)^2B(\lambda)\dd\lambda\right].
 \end{aligned}
\end{equation}
The bound $\|Z(s)\|_{H^2}\le\eps_\delta$ makes the constant uniform in
$s<t$. By \eqref{eq:time-window-notation},
\[
 \begin{aligned}
 \sup_{s\le\rho\le t}\sqrt{\mathcal T(\rho)}+Q_I
 &\le \sqrt{2}\left(\sup_{s\le\rho\le t}\mathcal T(\rho)
                           +\int_s^t\Gamma(\lambda)^2\dd\lambda\right)^{\frac{1}{2}}\\
 &\le C_\delta\left[\sqrt{\mathcal T(s)}
                +\left(\int_s^tS(\lambda)^2B(\lambda)\dd\lambda\right)^{\frac{1}{2}}\right]
 \le C_\delta(\mathfrak b_s+S_I\beta_I).
 \end{aligned}
\]
This gives \eqref{eq:time-window-estimate}.

For the solution of Proposition~\ref{prop:global}, fix $0<T<\infty$. By
\eqref{eq:global-regularity}, $b$ and $d=\ell(Z)$ are continuous in time,
and $Z$ and its spatial derivatives through order two are continuous in
$L^2$. By Sobolev embedding, the coefficients used in
\eqref{eq:time-F-controlled-by-B} are continuous in $L^\infty$. Equation
\eqref{eq:energy-perturbation-system} thus gives $F\in C([0,T];L^2)$ and
$W_\xi\in C([0,T];L^2)$. The mass equation and the bound for
$(\psi_t,\zeta_t)$ in Proposition~\ref{prop:global} give
\[
 \begin{aligned}
 \|F\|_{L^2(0,T;H^1)}
 &\le \|b\phi_\xi+\psi_\xi+dV'\|_{L^2(0,T;H^1)}
       +\|(\psi_t,\zeta_t)\|_{L^2(0,T;H^1)}\\
 &\le C_\delta\left(\sqrt T\sup_{0\le t\le T}\|Z(t)\|_{H^2}
         +\|(\psi_\xi,\zeta_\xi)\|_{L^2(0,T;H^2)}
         +\|d\|_{L^2(0,T)}\right)<\infty .
 \end{aligned}
\]
Differentiate \eqref{eq:energy-perturbation-system} in distributions to get
\eqref{eq:time-differentiated-system}. With $\|\partial_\xi
f\|_{H^{-1}}\le\|f\|_2$,
\eqref{eq:time-chain-rule}--\eqref{eq:time-phase-functional-bounds} give
\[
 \begin{aligned}
 \|F_t\|_{H^{-1}}
 &\le C_\delta\left(\|F_\xi\|_2+\|F\|_2+\|p_t\|_2
                         +|e|\|q_\xi\|_2\right)
 \le C_\delta\|F\|_{H^1}.
 \end{aligned}
\]
In particular, $F_t\in L^2(0,T;H^{-1})$. The square roots of the weights in
$2E_F$ are uniformly bounded in $W^{1,\infty}$. Since $a$ is fixed and
$(v_t,\theta_t)=(W,Q)$, their time derivatives are bounded in $L^\infty$ by
$C_\delta(\|W\|_\infty+\|Q\|_\infty)$. Multiplication by the weights is
bounded on $H^1$ and $H^{-1}$. The product rule and Sobolev embedding give
\[
 \begin{aligned}
 &\left\|\partial_t\left(\frac{\sqrt{aR\theta}}{v}W,\sqrt a\,Y,
             \sqrt{\frac{ac_v}{\theta}}\,Q\right)\right\|_{H^{-1}}\\
 &\quad\le C_\delta\left[\|F_t\|_{H^{-1}}+
                   (\|W\|_\infty+\|Q\|_\infty)\|F\|_2\right]
 \le C_\delta(1+\|F\|_2)\|F\|_{H^1}.
 \end{aligned}
\]
The weighted products thus lie in $L^2(0,T;H^1)$ with time derivatives in
$L^2(0,T;H^{-1})$. The $H^1$--$H^{-1}$ energy identity makes $E_F$
absolutely continuous and gives, for almost every $t$,
\[
 \begin{aligned}
 E_F'
 &=\left\langle W_t,\frac{aR\theta}{v^2}W\right\rangle
    +\langle Y_t,aY\rangle
    +\left\langle Q_t,\frac{ac_v}{\theta}Q\right\rangle\\
 &\quad+\frac{1}{2}\int_{\R}a\left[
       W^2\partial_t\left(\frac{R\theta}{v^2}\right)
       +Q^2\partial_t\left(\frac{c_v}{\theta}\right)\right]\dd\xi .
 \end{aligned}
\]
Substituting \eqref{eq:time-differentiated-system} and integrating
by parts in the $H^{-1}$--$H^1$ pairing give \eqref{eq:time-exact-EF}.

The identity for $D_t(\frac{\mu W_\xi}{v}-Y)$ in \eqref{eq:time-exact-KF}
holds in $L^2(0,T;H^{-1})$ because $D_tW_\xi=Y_{\xi\xi}+e v_{\xi\xi}$ holds
there. By \eqref{eq:global-regularity}, $(\frac{u_\xi}{v^2})_\xi$ is
uniformly bounded in $L^2$. Apply \eqref{eq:time-chain-rule},
\eqref{eq:time-phase-functional-bounds} and Sobolev embedding to obtain
\[
 \begin{aligned}
 \left\|D_t\left(\frac{\mu W_\xi}{v}-Y\right)\right\|_2
 &\le\|(p_t)_\xi\|_2+C_\delta\|Y_\xi\|_2
       +\mu\left\|\left(\frac{u_\xi}{v^2}\right)_\xi\right\|_2
           \|W\|_\infty+C_\delta|e|\|(v_{\xi\xi},u_\xi)\|_2\\
 &\le C_\delta\|F\|_{H^1}.
 \end{aligned}
\]
Translation by $\xi\mapsto\xi-\int_0^t b(\rho)\dd\rho$ turns the material
derivative into the time derivative of the translated field. Thus the
translated $\frac{\mu W_\xi}{v}-Y$ lies in $H^1(0,T;L^2)$. The translated
$Y$ lies in $L^2(0,T;H^1)$ with time derivative in $L^2(0,T;H^{-1})$. Since
$D_tv=u_\xi$ is bounded, the translated weight $\sqrt{\frac{v}{\mu}}$ has
bounded time derivative. The product rule and the $H^1$--$H^{-1}$ energy
identity make the integral absolutely continuous and give, for almost every
$t$,
\begin{align*}
 &\frac{\mathrm d}{\mathrm dt}\int_\R
             \left(\frac{\mu}{2v}W_\xi^2-YW_\xi\right)\dd\xi
 =\frac{\mathrm d}{\mathrm dt}\int_\R\frac{v}{2\mu}
              \left[\left(\frac{\mu W_\xi}{v}-Y\right)^2-Y^2\right]\dd\xi\\*
 &\quad=\int_\R\frac{v}{\mu}
        \left(\frac{\mu W_\xi}{v}-Y\right)
             D_t\left(\frac{\mu W_\xi}{v}-Y\right)\dd\xi
        -\left\langle D_tY,\frac{vY}{\mu}\right\rangle\\*
 &\qquad+\frac{1}{2\mu}\int_\R u_\xi
                \left[\left(\frac{\mu W_\xi}{v}-Y\right)^2-Y^2\right]\dd\xi\\
 &\quad=\int_\R W_\xi D_t\left(\frac{\mu W_\xi}{v}-Y\right)\dd\xi
        -\left\langle D_t\left(\frac{\mu W_\xi}{v}\right),
                                      \frac{vY}{\mu}\right\rangle
        +\int_\R u_\xi\left(\frac{\mu}{2v^2}W_\xi^2
                                     -\frac{YW_\xi}{v}\right)\dd\xi\\
 &\quad=\int_\R W_\xi D_t\left(\frac{\mu W_\xi}{v}-Y\right)\dd\xi
        +\|Y_\xi\|_2^2-e\int_\R Yv_{\xi\xi}\dd\xi
        +\frac{\mu}{2}\int_\R\frac{u_\xi}{v^2}W_\xi^2\dd\xi .
\end{align*}
The last equality uses $D_tv=u_\xi$ and $D_tW_\xi=Y_{\xi\xi}+e v_{\xi\xi}$
in the pairing. Substitute the material derivative and
\eqref{eq:time-chain-rule} to get \eqref{eq:time-exact-KF}. Both terms in
\eqref{eq:time-energy-definition} are absolutely continuous on $[0,T]$, so
$\mathcal T$ is also absolutely continuous there. By
\eqref{eq:time-energy-equivalence}, \eqref{eq:def-B} and
\eqref{eq:time-phase-functional-bounds},
\[
 \mathcal T(0)
 \le C_\delta\bigl(\|F(0)\|_2^2+\|W_\xi(0)\|_2^2\bigr)
 \le C_\delta B(0)\le C_\delta\|Z_0\|_{H^2}^2 .
\]
Integrating \eqref{eq:time-energy-inequality}, we have, for
$0\le s<t\le T$,
\[
 \begin{aligned}
 &\mathcal T(t)+c_\delta\int_s^t\Gamma(\rho)^2\dd\rho
 =\mathcal T(s)+\int_s^t
      \left[\mathcal T'(\rho)+c_\delta\Gamma(\rho)^2\right]\dd\rho\\
 &\quad\le\mathcal T(s)+C_\delta\int_s^t
               B(\rho)\left(\mathcal T(\rho)+S(\rho)^2\right)\dd\rho.
 \end{aligned}
\]
The integrating factor argument gives \eqref{eq:time-window-estimate} for
every such $s,t$. The time $T$ is arbitrary.
\end{proof}

\begin{remark}
On each interval $I$, the contribution of $S^2B$ in
\eqref{eq:time-window-estimate} is controlled by $S_I\beta_I$.
\end{remark}

\section{Stationary comparison profiles and spatial estimates}
\label{sec:stationary}

The comparison profiles have variable end states. Four fixed linear
functionals of the perturbation determine their parameters. This
construction removes the bounded stationary solutions of the spatial
operator without an integral condition on the initial perturbation.

\subsection{The stationary family and four fixed moments}

For a profile of speed $s>0$ in the coordinate $\xi=x-st$,
let $J=(J_0,J_1,J_2)$ be the fluxes in \eqref{eq:fluxes}.
The mass relation gives $U=-sV-J_0$ and $U'=-sV'$, so
\begin{equation}\label{eq:st-fluxes}
 \begin{aligned}
 \frac{\mu s}{V}V'
 &=-\frac{\mu U'}{V}=sU-P+J_1,\\
 J_2-UJ_1
 &=-s\left(c_v\Theta+\frac{U^2}{2}\right)+sU^2
                         -\frac{\kappa\Theta'}{V}
 =-sc_v\Theta+\frac{sU^2}{2}-\frac{\kappa\Theta'}{V}.
 \end{aligned}
\end{equation}
Solving \eqref{eq:st-fluxes} for $V'$ and $\Theta'$ with
$P=\frac{R\Theta}{V}$ gives
\begin{equation}\label{eq:st-reduced}
 \begin{aligned}
 U&=-sV-J_0,\\
 V'&=\frac{V}{\mu s}
       \left(-s^2V-sJ_0-\frac{R\Theta}{V}+J_1\right),\\
 \Theta'&=\frac{V}{\kappa}
       \left(-sc_v\Theta+\frac{s}{2}(sV+J_0)^2
                            -J_1(sV+J_0)-J_2\right).
 \end{aligned}
\end{equation}
The last two equations form a closed system for $(V,\Theta)$.
Write $\bar J$ for the fluxes of the reference profile $\bar q$.
Choose the exponent $\eta_0$ in the weight $\omega$ of
\eqref{eq:anchored-weight} sufficiently small for the tail estimates.
For an integer $k\ge0$, $E^k_{2\eta_0}$ consists of $C^k$ functions
with limits $f_\pm=\lim_{\xi\to\pm\infty}f(\xi)$ and finite norm
\begin{equation}\label{eq:st-profile-space}
 \|f\|_{E^k_{2\eta_0}}
 =|f_-|+|f_+|
  +\sum_{\pm}\sum_{j=0}^k
    \sup_{\pm\xi\ge0}e^{2\eta_0|\xi|}
       |\partial_\xi^j(f(\xi)-f_\pm)|.
\end{equation}
By \eqref{eq:profilebounds}, for $\pm\xi\ge0$,
\[
 \begin{aligned}
 \sum_{j=0}^k|\partial_\xi^j(\bar q(\xi)-q_\pm)|
 &\le \int_0^\infty|\bar q'(\xi\pm y)|\dd y
           +\sum_{j=1}^k|\partial_\xi^j\bar q(\xi)|\\
 &\le C_{\delta,k}e^{-c\delta|\xi|}.
 \end{aligned}
\]
Thus $\bar q\in E^k_{2\eta_0}$ when $2\eta_0<c\delta$.

\begin{lemma}\label{lem:st-raw-family}
For $(s,J)$ near $(\sigma,\bar J)$, \eqref{eq:st-reduced} has a family
of connecting profiles $\widehat q_{s,J}$ with
$\widehat V_{s,J}(0)=V(0)$. Each connecting profile sufficiently close
to $\bar q$ in $L^\infty$ is a translate of a member of the family.
The family is smooth into $E^k_{2\eta_0}$ for every fixed finite $k$.
\end{lemma}

\begin{proof}
At an equilibrium $(v,\vartheta)$, the brackets in
\eqref{eq:st-reduced} vanish and
$p=\frac{R\vartheta}{v}=J_1-s(sv+J_0)$. Differentiating the last
two right-hand sides in $(V,\Theta)$ gives
\begin{equation}\label{eq:st-spatial-matrix}
 \begin{aligned}
 K(v,\vartheta;s)
 &=\begin{pmatrix}
    \frac{v}{\mu s}\left(\frac{p}{v}-s^2\right)
       &-\frac{R}{\mu s}\\[2mm]
    -\frac{svp}{\kappa}&-\frac{sc_vv}{\kappa}
   \end{pmatrix},\\
 \det K
 &=\frac{c_vv^2}{\mu\kappa}\left(s^2-\frac{p}{v}\right)
      -\frac{Rvp}{\mu\kappa}
  =\frac{c_vv^2}{\mu\kappa}(s^2-c^2),
 \end{aligned}
\end{equation}
where $c^2=\frac{\gamma p}{v}$ by \eqref{eq:lax}.
At the ends of the reference profile, \eqref{eq:lax} gives
\[
 \begin{aligned}
 \det K(v_-,\theta_-;\sigma)
 &=\frac{c_vv_-^2}{\mu\kappa}(\sigma^2-c_-^2)
 <0<\frac{c_vv_+^2}{\mu\kappa}(\sigma^2-c_+^2)
 =\det K(v_+,\theta_+;\sigma),\\
 \operatorname{tr}K(v_+,\theta_+;\sigma)
 &=\frac{v_+}{\mu\sigma}\left(\frac{c_+^2}{\gamma}-\sigma^2\right)
      -\frac{\sigma c_v v_+}{\kappa}<0.
 \end{aligned}
\]
Since $\det K\ne0$ at both ends, the implicit function theorem gives
unique nearby equilibria
$(\widehat V_{s,J,\pm},\widehat\Theta_{s,J,\pm})$ that depend smoothly
on $(s,J)$. The left equilibrium is a saddle with a one-dimensional
unstable manifold. The right equilibrium is a sink. The velocity
relation $U=-sV-J_0$ and the mean value theorem give
\[
 \begin{aligned}
 |\widehat q_{s,J,\pm}-q_\pm|
 &\le (1+s)|\widehat V_{s,J,\pm}-v_\pm|
                +|\widehat\Theta_{s,J,\pm}-\theta_\pm|
                +v_\pm|s-\sigma|+|J_0-\bar J_0|\\
 &\le C_\delta\bigl(|s-\sigma|+|J-\bar J|\bigr).
 \end{aligned}
\]

Since $K_{12}=-\frac{R}{\mu s}\ne0$ in \eqref{eq:st-spatial-matrix},
the left unstable manifold is locally a graph over $V$ that depends
smoothly on $(s,J)$. Fix $\xi_0>0$ sufficiently large. For nearby
parameters the graph meets $V=V(-\xi_0)$ on the branch
$V>\widehat V_{s,J,-}$ because
\[
 V(-\xi_0)-\widehat V_{s,J,-}
 \ge V(-\xi_0)-v_- -C_\delta\bigl(|s-\sigma|+|J-\bar J|\bigr)
 \ge \frac{V(-\xi_0)-v_-}{2}>0 .
\]
The flow from the intersection depends smoothly on $(s,J)$.
By \eqref{eq:profilebounds}, $V'(0)>0$ at the reference profile.
The implicit function theorem gives a unique nearby crossing of
$V=V(0)$ that depends smoothly on $(s,J)$. Place this crossing
at $\xi=0$. Then $\widehat V_{s,J}(0)=V(0)$ and
$\widehat V_{s,J}'(0)\ge\frac{V'(0)}{2}>0$ for nearby parameters.
The orbit extends backward along the unstable branch to
$\widehat q_{s,J,-}$.

At the right end, the matrix
\[
 \begin{pmatrix}
  \frac{\sigma v_+p_+}{\kappa}&0\\[1mm]
  0&\frac{R}{\mu\sigma}
 \end{pmatrix}
 K(v_+,\theta_+;\sigma)
 =
 -\frac{v_+}{\mu\kappa}
 \begin{pmatrix}
 v_+p_+\left(\sigma^2-\frac{p_+}{v_+}\right)&Rp_+\\[1mm]
 Rp_+&Rc_v
 \end{pmatrix}
\]
is symmetric with positive determinant and negative first diagonal
entry by \eqref{eq:st-spatial-matrix} and \eqref{eq:lax}.
For nearby $(s,J)$, the same weights give a uniformly negative
quadratic form for the Jacobian of \eqref{eq:st-reduced} in a small
ball about the right equilibrium. Apply the mean value formula on
the segment from the equilibrium to the profile to obtain
\[
 \begin{aligned}
 &\frac{\mathrm d}{\mathrm d\xi}\left[
     \frac{\sigma v_+p_+}{\kappa}
                (\widehat V_{s,J}-\widehat V_{s,J,+})^2
    +\frac{R}{\mu\sigma}
                (\widehat\Theta_{s,J}-\widehat\Theta_{s,J,+})^2\right]\\
 &\quad=\frac{2\sigma v_+p_+}{\kappa}
             (\widehat V_{s,J}-\widehat V_{s,J,+})\widehat V_{s,J}'
    +\frac{2R}{\mu\sigma}
             (\widehat\Theta_{s,J}-\widehat\Theta_{s,J,+})\widehat\Theta_{s,J}'\\
 &\quad\le -c_\delta\left[
                (\widehat V_{s,J}-\widehat V_{s,J,+})^2
                +(\widehat\Theta_{s,J}-\widehat\Theta_{s,J,+})^2\right]\\
 &\quad\le -c_\delta\left[
     \frac{\sigma v_+p_+}{\kappa}
                (\widehat V_{s,J}-\widehat V_{s,J,+})^2
    +\frac{R}{\mu\sigma}
                (\widehat\Theta_{s,J}-\widehat\Theta_{s,J,+})^2\right].
 \end{aligned}
\]
A sufficiently small sublevel set lies in the ball and is forward
invariant. Smooth dependence on $[0,\xi_0]$ and
\eqref{eq:profilebounds} give
\[
 \begin{aligned}
 |\widehat q_{s,J}(\xi_0)-\widehat q_{s,J,+}|
 &\le |\widehat q_{s,J}(\xi_0)-\bar q(\xi_0)|
    +|\bar q(\xi_0)-q_+|+|\widehat q_{s,J,+}-q_+|\\
 &\le C_{\delta,\xi_0}\bigl(|s-\sigma|+|J-\bar J|\bigr)
            +C_\delta e^{-c\delta\xi_0}.
 \end{aligned}
\]
Take $\xi_0$ large and then restrict $(s,J)$ so that the orbit enters
this sublevel set. The orbit extends to all $\xi\ge\xi_0$.
Integrating the differential inequality and using
$\widehat U_{s,J}-\widehat U_{s,J,+}
=-s(\widehat V_{s,J}-\widehat V_{s,J,+})$ gives
\[
 |\widehat q_{s,J}(\xi)-\widehat q_{s,J,+}|
 \le C_\delta e^{-c_\delta(\xi-\xi_0)}
       |\widehat q_{s,J}(\xi_0)-\widehat q_{s,J,+}|,
 \qquad \xi\ge\xi_0.
\]
Both tails lie in fixed neighborhoods of positive equilibria.
On $[-\xi_0,\xi_0]$, continuous dependence keeps the orbit close to
the reference profile. The entire orbit lies in a compact set of
positive states.

On the two tails, the equations for the differences from the
equilibria take the form
\[
 \begin{aligned}
 r&=\pm\xi-\xi_0\ge0,\qquad
 z=\begin{cases}
 (\widehat V_{s,J}-\widehat V_{s,J,+},
       \widehat\Theta_{s,J}-\widehat\Theta_{s,J,+}),
       &\xi\ge\xi_0,\\
 \widehat V_{s,J}-\widehat V_{s,J,-},&\xi\le-\xi_0,
 \end{cases}\\
 z_r&=G(z;s,J),\qquad G(0;s,J)=0 .
 \end{aligned}
\]
On the left, the derivative of the reversed scalar equation at
$z=0$ is the negative of the unstable eigenvalue in
\eqref{eq:st-spatial-matrix}. This derivative is uniformly negative
near the left equilibrium. On the right, use the quadratic form above.
For sufficiently small $\eta_0$, the equation
$h_r=D_zG(z;s,J)h$ gives
\[
 \begin{aligned}
 \frac{\mathrm d}{\mathrm dr}|h|^2
 &=2\partial_zG(z;s,J)|h|^2\le-8\eta_0|h|^2,
                                  &&\xi\le-\xi_0,\\
 \frac{\mathrm d}{\mathrm dr}\left(
     \frac{\sigma v_+p_+}{\kappa}h_1^2+\frac{R}{\mu\sigma}h_2^2\right)
 &\le -8\eta_0\left(
     \frac{\sigma v_+p_+}{\kappa}h_1^2+\frac{R}{\mu\sigma}h_2^2\right),
                                  &&\xi\ge\xi_0 .
 \end{aligned}
\]
The mean value formula for $G(z;s,J)-G(0;s,J)$ gives the same
inequalities for $z$. Integration gives
\[
 |z(r)|\le C_\delta e^{-4\eta_0r}|z(0)|,\qquad
 |h(r)|\le C_\delta e^{-4\eta_0(r-\rho)}|h(\rho)|,\qquad r\ge\rho\ge0 .
\]
By smooth dependence on finite intervals,
$|\partial_{(s,J)}^\nu z(0)|\le C_{\delta,\nu}$. For $|\nu|\ge1$,
assume the $e^{-3\eta_0r}$ bound at lower orders. Since
$G(0;s,J)=0$ for every $(s,J)$, the chain rule gives
\[
 \begin{aligned}
 &\left|(\partial_{(s,J)}^\nu z)_r
            -D_zG(z;s,J)\partial_{(s,J)}^\nu z\right|\\
 &\quad\le\int_0^1
     \left|(D_z\partial_{(s,J)}^\nu G)(\rho z;s,J)z\right|\dd\rho
      +C_{\delta,\nu}\sum_{1\le|\beta|<|\nu|}|\partial_{(s,J)}^\beta z|\\
 &\quad\le C_{\delta,\nu}\left(
      |z|+\sum_{1\le|\beta|<|\nu|}|\partial_{(s,J)}^\beta z|
                           \right)\le C_{\delta,\nu}e^{-3\eta_0r}.
 \end{aligned}
\]
Variation of constants then gives
\[
 \begin{aligned}
 |\partial_{(s,J)}^\nu z(r)|
 &\le C_\delta e^{-4\eta_0r}|\partial_{(s,J)}^\nu z(0)|
    +C_{\delta,\nu}\int_0^r
          e^{-4\eta_0(r-\rho)}e^{-3\eta_0\rho}\dd\rho\\
 &\le C_{\delta,\nu}\left(e^{-4\eta_0r}
             +\frac{e^{-3\eta_0r}-e^{-4\eta_0r}}{\eta_0}\right)
 \le C_{\delta,\nu}e^{-3\eta_0r}.
 \end{aligned}
\]
By the left graph, $U=-sV-J_0$ and $z_r=G(z;s,J)$, each
spatial derivative of the end-state difference is a smooth function
of $(z;s,J)$ vanishing at $z=0$. Thus
\[
 \begin{aligned}
 &e^{2\eta_0|\xi|}
   |\partial_\xi^j\partial_{(s,J)}^\nu
        (\widehat q_{s,J}(\xi)-\widehat q_{s,J,\pm})|\\
 &\quad\le C_{\delta,j,\nu} e^{2\eta_0|\xi|}
          \sum_{|\beta|\le|\nu|}
            |\partial_{(s,J)}^\beta z(\pm\xi-\xi_0)|
 \le C_{\delta,j,\nu}e^{-\eta_0|\xi|},\qquad \pm\xi\ge\xi_0 .
 \end{aligned}
\]
The derivatives are continuous on compact intervals and at the end
states. The uniform tail bound gives their continuity in
\eqref{eq:st-profile-space}. On a smaller parameter ball, Taylor's
formula at each $\xi$ and at the end states gives
\[
 \begin{aligned}
 &\left\|\partial_{(s,J)}^\nu\widehat q_{s+\dot s,J+\dot J}
       -\partial_{(s,J)}^\nu\widehat q_{s,J}
       -D_{(s,J)}\partial_{(s,J)}^\nu\widehat q_{s,J}(\dot s,\dot J)
       \right\|_{E^k_{2\eta_0}}\\
 &\quad\le\int_0^1(1-\rho)\left\|
     D_{(s,J)}^2\partial_{(s,J)}^\nu
          \widehat q_{s+\rho\dot s,J+\rho\dot J}
          [(\dot s,\dot J),(\dot s,\dot J)]
          \right\|_{E^k_{2\eta_0}}\dd\rho\\
 &\quad\le C_{\delta,k,\nu}\bigl(|\dot s|+|\dot J|\bigr)^2 .
 \end{aligned}
\]
Hence the family is smooth into $E^k_{2\eta_0}$.

For a connecting profile uniformly close to $\bar q$, the mass
relation at both ends and the fluxes \eqref{eq:fluxes} at either end
give nearby parameters. Its negative tail lies on the local unstable
graph. At $-\xi_0$ its volume differs from $V(-\xi_0)$ by less than
$\frac{V(-\xi_0)-v_-}{2}$, so it lies on the branch
$V>\widehat V_{s,J,-}$. By uniqueness for \eqref{eq:st-reduced}, it is
a translate of the corresponding profile.
\end{proof}

The profile speed is prescribed by \eqref{eq:moving-frame}.
Lemma~\ref{lem:st-raw-family} bounds the spatial and parameter
derivatives of $\widehat q_{s,J}$ on a smaller neighborhood. Since the
kernel of $\ell$ is integrable by \eqref{eq:feedbackbounds},
differentiation under the integral gives smooth dependence on
$(s,J,\tau)$. Using $\widehat q_{\sigma,\bar J}=\bar q$ and
$\ell(\bar q')<0$ from \eqref{eq:feedbackbounds}, we obtain
\begin{equation}\label{eq:st-phase-nondeg}
 \begin{aligned}
 &\partial_\tau\ell\bigl(\widehat q_{s,J}(\cdot+\tau)-\bar q\bigr)
   =\ell\bigl(\widehat q_{s,J}'(\cdot+\tau)\bigr)\\
 &\quad\le \ell(\bar q')
     +C_\delta\left(\|\widehat q_{s,J}'-\bar q'\|_\infty
                  +|\tau|\|\widehat q_{s,J}''\|_\infty\right)\\
 &\quad\le \ell(\bar q')+C_\delta\bigl(|s-\sigma|+|J-\bar J|+|\tau|\bigr)
 \le -\frac{|\ell(\bar q')|}{2}<0
 \end{aligned}
\end{equation}
for nearby $(s,J)$ and small $\tau$. The speed condition is
\begin{equation}\label{eq:st-phase}
 \ell\bigl(\widehat q_{s,J}(\cdot+\tau)-\bar q\bigr)=s-\sigma .
\end{equation}
It holds at $(s,J,\tau)=(\sigma,\bar J,0)$. The implicit function
theorem gives a unique small $\tau(s,J)$, smooth in $(s,J)$, with
$\tau(\sigma,\bar J)=0$. By \eqref{eq:st-phase-nondeg} and
\eqref{eq:st-phase},
\[
 \begin{aligned}
 \frac{|\ell(\bar q')|}{2}|\tau(s,J)|
 &\le \left|\int_0^{\tau(s,J)}
             \ell\bigl(\widehat q_{s,J}'(\cdot+\rho)\bigr)\dd\rho\right|\\
 &=\left|s-\sigma-\ell(\widehat q_{s,J}-\bar q)\right|\\
 &\le |s-\sigma|+C_\delta\|\widehat q_{s,J}-\bar q\|_{L^\infty}
 \le C_\delta\bigl(|s-\sigma|+|J-\bar J|\bigr).
 \end{aligned}
\]

Let $\chi_0\in C_c^\infty(\R)$ satisfy $\int_{\R}\chi_0\dd\xi=1$.
For $B=(f,g,h)$ with $\int_\R\omega|B|\dd\xi<\infty$, the speed
variation and the three flux variations averaged against $\chi_0$
define
\begin{equation}\label{eq:st-moments}
 \begin{aligned}
 (\mathfrak mB)_0&=\ell(B),\\
 (\mathfrak mB)_1&=\int_{\R}\chi_0(-\ell(B)V-\sigma f-g)\dd\xi,\\
 (\mathfrak mB)_2&=\int_{\R}\chi_0
     \left(-\ell(B)U-\sigma g+\frac{Rh}{V}+\frac{\mu U'-R\Theta}{V^2}f\right)
           \dd\xi
         +\mu\int_{\R}\left(\frac{\chi_0}{V}\right)'g\dd\xi,\\
 (\mathfrak mB)_3&=\int_{\R}\chi_0\Bigg[
      -\ell(B)\left(c_v\Theta+\frac{U^2}{2}\right)
      +\left(\frac{RU}{V}-\sigma c_v\right)h+(P-\sigma U)g\\
 &\quad\qquad+\frac{\kappa\Theta'+\mu UU'-R\Theta U}{V^2}f\Bigg]\dd\xi
      +\int_{\R}\left(\frac{\chi_0}{V}\right)'(\kappa h+\mu Ug)\dd\xi .
 \end{aligned}
\end{equation}
By \eqref{eq:profilebounds}, the choice of $\eta_0$ and compact support,
$V'+|\chi_0|+|(\frac{\chi_0}{V})'|
\le C_\delta e^{-\eta_0|\xi|}\le C_\delta\omega(\xi)$.
Thus \eqref{eq:st-moments}, \eqref{eq:feedbackbounds} and H\"older's
inequality give
\begin{equation}\label{eq:st-moment-bounds}
 \begin{aligned}
 |\mathfrak mB|
 &\le C_\delta\left(|\ell(B)|+\int_{\R}
             \left(|\chi_0|+\left|\left(\frac{\chi_0}{V}\right)'\right|\right)
             |B|\dd\xi\right)\\
 &\le C_\delta\int_{\R}
       \left(V'+|\chi_0|+\left|\left(\frac{\chi_0}{V}\right)'\right|\right)
                       |B|\dd\xi
 \le C_\delta\int_{\R}\omega |B|\dd\xi\\
 &\le C_\delta\min\left\{\|\omega\|_1\|B\|_\infty,
           \|\sqrt{\omega}\|_2\,m_\omega(B),\,\|\omega\|_2\|B\|_2\right\}\\
 &\le C_\delta\min\left\{\|B\|_\infty,\,m_\omega(B),\,\|B\|_2\right\}.
 \end{aligned}
\end{equation}
The weight norms are finite since $\omega(\xi)\le e^{-\eta_0|\xi|}$.

\begin{proposition}[Stationary profiles with prescribed moments]
\label{prop:stationary}
For $\alpha$ near $0$ in $\R^4$, there is a smooth family of stationary
profiles $q_\alpha$, with $\mathcal R(\alpha)=q_\alpha-\bar q$ satisfying
\begin{equation}\label{eq:st-family}
 \mathcal R(0)=0,\qquad
 \mathfrak m\mathcal R(\alpha)=\alpha,\qquad
 s_\alpha=\sigma+\ell(\mathcal R(\alpha))=\sigma+\alpha_0.
\end{equation}
If a connecting profile $q$ is sufficiently close to $\bar q$ in
$L^\infty$, has speed $\sigma+\alpha_0$ and satisfies
$\mathfrak m(q-\bar q)=\alpha$, then $q=q_\alpha$.
For every integer $k\ge0$, the map $\alpha\mapsto\mathcal R(\alpha)$
is smooth into $E^k_{2\eta_0}$ and satisfies
\begin{equation}\label{eq:st-family-bounds}
 \|\mathcal R(\alpha)\|_{E^k_{2\eta_0}}\le C_{\delta,k}|\alpha|,
 \qquad
 \|\mathcal R(\alpha)-D_\alpha\mathcal R(0)\alpha\|_{E^k_{2\eta_0}}
       \le C_{\delta,k}|\alpha|^2 .
\end{equation}
The parameter derivatives are uniformly bounded in $E^k_{2\eta_0}$
on a smaller ball. For $0\le j\le k$, $\pm\xi\ge0$ and $|\nu|\le2$,
\begin{equation}\label{eq:st-parameter-tails}
 \left|\partial_\xi^j\partial_\alpha^\nu
    \bigl(\mathcal R(\alpha)(\xi)-\mathcal R(\alpha)_\pm\bigr)\right|
 \le C_{\delta,k}e^{-2\eta_0|\xi|}
 \begin{cases}
 |\alpha|,&|\nu|=0,\\
 1,&|\nu|=1,2 .
 \end{cases}
\end{equation}
\end{proposition}

\begin{proof}
Lemma~\ref{lem:st-raw-family} and smoothness of $\tau(s,J)$ bound the
parameter derivatives of the translated profiles in $L^\infty$.
By \eqref{eq:st-moment-bounds}, we can differentiate the moments under
the integrals.
For a direction $(\dot s,\dot J)$ at $(\sigma,\bar J)$, let
$B=(f,g,h)$ be the variation of $\widehat q_{s,J}(\cdot+\tau(s,J))$.
Differentiating \eqref{eq:st-phase} gives $\dot s=\ell(B)$.
Differentiating \eqref{eq:fluxes} gives
\begin{equation}\label{eq:st-tangent-fluxes}
 \begin{aligned}
 \dot J_0&=-\ell(B)V-\sigma f-g,\\
 \dot J_1&=-\ell(B)U-\sigma g+\frac{Rh}{V}
                     +\frac{\mu U'-R\Theta}{V^2}f-\frac{\mu g'}{V},\\
 \dot J_2&=-\ell(B)\left(c_v\Theta+\frac{U^2}{2}\right)
       +\left(\frac{RU}{V}-\sigma c_v\right)h+(P-\sigma U)g\\
 &\quad+\frac{\kappa\Theta'+\mu UU'-R\Theta U}{V^2}f
                         -\frac{\kappa h'+\mu(Ug)'}{V} .
 \end{aligned}
\end{equation}
Each $\dot J_i$ is constant. The compact support of $\chi_0$ and
$\int_{\R}\chi_0\dd\xi=1$ give, by \eqref{eq:st-moments}
and \eqref{eq:st-tangent-fluxes},
\[
 \begin{aligned}
 (\mathfrak mB)_0&=\ell(B)=\dot s,\qquad
 (\mathfrak mB)_1=\int_{\R}\chi_0\dot J_0\dd\xi=\dot J_0,\\
 (\mathfrak mB)_2&=\int_{\R}\left[
        \chi_0\dot J_1+\mu\left(\frac{\chi_0 g}{V}\right)'\right]\dd\xi
        =\dot J_1,\\
 (\mathfrak mB)_3&=\int_{\R}\left[
        \chi_0\dot J_2+\left(\frac{\chi_0(\kappa h+\mu Ug)}{V}\right)'
                                           \right]\dd\xi
        =\dot J_2 .
 \end{aligned}
\]
Since $\widehat q_{\sigma,\bar J}=\bar q$ and $\tau(\sigma,\bar J)=0$,
the moment map takes the value $0$ at $(\sigma,\bar J)$ and satisfies
\[
 \left.D_{(s,J)}
       \left[\mathfrak m\bigl(\widehat q_{s,J}(\cdot+\tau(s,J))-\bar q\bigr)
       \right]\right|_{(s,J)=(\sigma,\bar J)}(\dot s,\dot J)
       =\mathfrak mB=(\dot s,\dot J).
\]
The inverse function theorem gives a smooth local inverse
$(s_\alpha,J_\alpha)$ for $\alpha$ near $0$. The associated profiles satisfy
\[
 \begin{aligned}
 q_\alpha&=\widehat q_{s_\alpha,J_\alpha}
                         (\cdot+\tau(s_\alpha,J_\alpha)),\\
 \mathfrak m\mathcal R(\alpha)&=\mathfrak m(q_\alpha-\bar q)=\alpha,
                                  \qquad \mathcal R(0)=0,\\
 s_\alpha-\sigma&=\ell(\mathcal R(\alpha))
                        =(\mathfrak m\mathcal R(\alpha))_0=\alpha_0 .
 \end{aligned}
\]
If a connecting profile $q$ is sufficiently close to $\bar q$ in
$L^\infty$, has speed $s=\sigma+\alpha_0$ and satisfies
$\mathfrak m(q-\bar q)=\alpha$, Lemma~\ref{lem:st-raw-family} gives
$q=\widehat q_{s,J}(\cdot+\tau)$ with nearby $(s,J)$. Its volume at
$-\tau$ equals $V(0)$, so
\[
 \left(\min_{|\xi|\le1}V'(\xi)\right)\min\{|\tau|,1\}
 \le |V(-\tau)-V(0)|\le\|q-\bar q\|_\infty .
\]
Thus $\tau$ is small. Since $s-\sigma=\alpha_0=\ell(q-\bar q)$,
\eqref{eq:st-phase} gives $\tau=\tau(s,J)$, and the local inverse
gives $(s,J)=(s_\alpha,J_\alpha)$. Hence $q=q_\alpha$.

For $|\tau|\le1$ and $\pm\xi\ge\xi_0+1$, we have
$\pm(\xi+\tau)\ge\pm\xi-|\tau|\ge\xi_0$. Hence
Lemma~\ref{lem:st-raw-family} gives
\[
 \begin{aligned}
 &e^{2\eta_0|\xi|}
 \left|\partial_\xi^j\partial_{(s,J)}^\nu\partial_\tau^m
   \bigl(\widehat q_{s,J}(\xi+\tau)-\widehat q_{s,J,\pm}\bigr)\right|\\
 &\quad=e^{2\eta_0|\xi|}
 \left|\partial_\xi^{j+m}\partial_{(s,J)}^\nu
   \bigl(\widehat q_{s,J}(\xi+\tau)-\widehat q_{s,J,\pm}\bigr)\right|\\
 &\quad\le C_{\delta,j+m,\nu}e^{2\eta_0|\xi|-3\eta_0|\xi+\tau|}
 \le C_{\delta,j+m,\nu}e^{-\eta_0|\xi|}.
 \end{aligned}
\]
For $|\xi|\le\xi_0+1$, the argument $\xi+\tau$ lies in
$[-\xi_0-2,\xi_0+2]$. The mixed derivatives depend smoothly on the
parameters there and at the end states. The tail bound gives
continuity in \eqref{eq:st-profile-space}. Taylor's formula, as in
Lemma~\ref{lem:st-raw-family}, gives joint smoothness of the translated
family into $E^k_{2\eta_0}$. The smooth phase and local inverse
give uniform bounds on the $\alpha$-derivatives on a smaller ball.
Since $\mathcal R(0)=0$, Taylor's formula in $E^k_{2\eta_0}$ gives
\[
 \begin{aligned}
 \|\mathcal R(\alpha)\|_{E^k_{2\eta_0}}
 &\le \int_0^1
       \|D_\alpha\mathcal R(\rho\alpha)\alpha\|_{E^k_{2\eta_0}}\dd\rho
 \le C_{\delta,k}|\alpha|,\\
 \|\mathcal R(\alpha)-D_\alpha\mathcal R(0)\alpha\|_{E^k_{2\eta_0}}
 &\le \int_0^1(1-\rho)
       \|D_\alpha^2\mathcal R(\rho\alpha)[\alpha,\alpha]\|_{E^k_{2\eta_0}}\dd\rho
 \le C_{\delta,k}|\alpha|^2 .
 \end{aligned}
\]
Evaluation at the end states is bounded in
\eqref{eq:st-profile-space} and commutes with parameter
differentiation. The derivative bounds and \eqref{eq:st-family-bounds}
give
\[
 \begin{aligned}
 &\left|\partial_\xi^j\partial_\alpha^\nu
          \bigl(\mathcal R(\alpha)(\xi)-\mathcal R(\alpha)_\pm\bigr)\right|\\
 &\quad\le e^{-2\eta_0|\xi|}
                \|\partial_\alpha^\nu\mathcal R(\alpha)\|_{E^k_{2\eta_0}}
 \le C_{\delta,k}e^{-2\eta_0|\xi|}
 \begin{cases}
 |\alpha|,&|\nu|=0,\\
 1,&|\nu|=1,2 ,
 \end{cases}
 \qquad \pm\xi\ge0,\quad 0\le j\le k .
 \end{aligned}
\]
This is \eqref{eq:st-parameter-tails}.
\end{proof}

\begin{lemma}[End states and incoming traces]\label{lemma:incoming}
The end-state map of the family in Proposition~\ref{prop:stationary}
satisfies
\begin{equation}\label{eq:st-end-injectivity}
 |\alpha|\le C_\delta
    \bigl(|\mathcal R(\alpha)_-|+|\mathcal R(\alpha)_+|\bigr).
\end{equation}
Moreover, the three right traces and the left $3$-trace in
\eqref{eq:normalform-definition} give local coordinates near $\alpha=0$.
\end{lemma}

\begin{proof}
For small $\alpha$, \eqref{eq:st-family-bounds} and the mass
Rankine--Hugoniot relation \eqref{eq:RH} for each profile give
\[
 \begin{aligned}
 v_{\alpha,+}-v_{\alpha,-}
 &\ge \delta-|\mathcal R(\alpha)_+|-|\mathcal R(\alpha)_-|
 \ge \delta-C_{\delta,0}|\alpha|\ge\frac{\delta}{2},\\
 |s_\alpha-\sigma|
 &=\frac{\big|(u_{\alpha,+}-u_+)-(u_{\alpha,-}-u_-)
    +\sigma\big[(v_{\alpha,+}-v_+)-(v_{\alpha,-}-v_-)\big]\big|}
         {v_{\alpha,+}-v_{\alpha,-}}\\
 &\le \frac{2(1+\sigma)}{\delta}
           \bigl(|\mathcal R(\alpha)_-|+|\mathcal R(\alpha)_+|\bigr).
 \end{aligned}
\]
At the right end, \eqref{eq:profilebounds} and
\eqref{eq:st-parameter-tails} make the derivative terms in
\eqref{eq:fluxes} vanish. By \eqref{eq:st-family-bounds},
$v_{\alpha,+}\ge\frac{v_+}{2}$. Evaluating \eqref{eq:fluxes} there
and using \eqref{eq:st-family}, \eqref{eq:st-moment-bounds} and
smooth dependence of the translated profiles on $(s,J)$ gives
\[
 \begin{aligned}
 |\alpha|=|\mathfrak m\mathcal R(\alpha)|
 &\le C_\delta\|q_\alpha-\bar q\|_\infty
 \le C_\delta\bigl(|s_\alpha-\sigma|+|J_\alpha-\bar J|\bigr)\\
 &\le C_\delta\left(
     |s_\alpha-\sigma|+|\mathcal R(\alpha)_+|
      +\left|\frac{R\theta_{\alpha,+}}{v_{\alpha,+}}
                  -\frac{R\theta_+}{v_+}\right|\right)\\
 &\le C_\delta\left(
     |s_\alpha-\sigma|+|\mathcal R(\alpha)_+|
      +\frac{R|\theta_{\alpha,+}-\theta_+|}{v_{\alpha,+}}
      +\frac{R\theta_+|v_{\alpha,+}-v_+|}{v_{\alpha,+}v_+}\right)\\
 &\le C_\delta\bigl(|s_\alpha-\sigma|+|\mathcal R(\alpha)_+|\bigr)
 \le C_\delta\bigl(|\mathcal R(\alpha)_-|+|\mathcal R(\alpha)_+|\bigr).
 \end{aligned}
\]

In the conservative variables $\mathcal U=(v,u,E)$,
$p=\frac{\gamma-1}{v}\left(E-\frac{u^2}{2}\right)$. At the left
state, the Jacobian $\mathcal A$ of the Euler flux $(-u,p,pu)$
and its third left eigenvector are
\[
 \begin{aligned}
 \mathcal A_-
 &=\begin{pmatrix}
   0&-1&0\\[1mm]
   -\frac{p_-}{v_-}&-\frac{(\gamma-1)u_-}{v_-}
                         &\frac{\gamma-1}{v_-}\\[2mm]
   -\frac{u_-p_-}{v_-}&p_--\frac{(\gamma-1)u_-^2}{v_-}
                         &\frac{(\gamma-1)u_-}{v_-}
  \end{pmatrix},\\
 l_3^-&=\left(\frac{1}{2\gamma},
 -\frac{1}{2c_-}+\frac{(\gamma-1)u_-}{2\gamma p_-},
 -\frac{\gamma-1}{2\gamma p_-}\right).
 \end{aligned}
\]
Using $c_-^2=\frac{\gamma p_-}{v_-}$, multiplication gives
\[
 \begin{aligned}
 l_3^-\mathcal A_-
 &=\left(\frac{p_-}{2v_-c_-},
       -\frac{1}{2}+\frac{(\gamma-1)u_-}{2v_-c_-},
       -\frac{\gamma-1}{2v_-c_-}\right)
 =c_-l_3^-,\\
 \mathcal A_-(1,-c_-,-p_--u_-c_-)^{\mathsf T}
 &=\left(c_-,-\frac{\gamma p_-}{v_-},
           -c_-p_--\frac{\gamma p_-u_-}{v_-}\right)^{\mathsf T}\\
 &=c_-(1,-c_-,-p_--u_-c_-)^{\mathsf T},\\
 l_3^-(1,-c_-,-p_--u_-c_-)^{\mathsf T}
 &=\frac{1}{2\gamma}+\frac{1}{2}+\frac{\gamma-1}{2\gamma}=1 .
 \end{aligned}
\]
By \eqref{eq:RH},
\[
 \begin{aligned}
 [v]&=\delta,\qquad [u]=-\sigma\delta,\qquad [p]=\sigma[u],\\
 [E]-u_-[u]
 &=\frac{[pu]-u_-[p]}{\sigma}
   =\frac{p_+[u]}{\sigma}=-\delta p_+ .
 \end{aligned}
\]
Consequently,
\begin{equation}\label{eq:st-incoming-jump}
 \begin{aligned}
 l_3^-[\mathcal U]
 &=\frac{[v]}{2\gamma}-\frac{[u]}{2c_-}
                 -\frac{\gamma-1}{2\gamma p_-}\bigl([E]-u_-[u]\bigr)\\
 &=\delta\left(\frac{1}{2\gamma}+\frac{\sigma}{2c_-}
                  +\frac{(\gamma-1)p_+}{2\gamma p_-}\right)
 \ge\frac{\delta}{2\gamma}>0 .
 \end{aligned}
\end{equation}

By \eqref{eq:lax}, the moving-frame speeds on the right are all
negative, and the incoming speed on the left is $c_--\sigma>0$.
Consider the trace map
\[
 \alpha\longmapsto
 \bigl(\mathcal U_{\alpha,+}-\mathcal U_+,\,
       l_3^-(\mathcal U_{\alpha,-}-\mathcal U_-)\bigr).
\]
Let $\dot\alpha$ lie in the kernel of the derivative at $\alpha=0$
and let $B$ be its profile variation. Then $\dot{\mathcal U}_+=0$
and $l_3^-\dot{\mathcal U}_-=0$. Differentiate \eqref{eq:RH},
apply $l_3^-$ and use \eqref{eq:st-incoming-jump} to obtain
\[
 \begin{aligned}
 &(\mathcal A_+-\sigma I)\dot{\mathcal U}_+
  -(\mathcal A_--\sigma I)\dot{\mathcal U}_-
      =[\mathcal U]\,\dot s,\\
 &(l_3^-[\mathcal U])\dot s
    =l_3^-(\mathcal A_+-\sigma I)\dot{\mathcal U}_+
                   -(c_--\sigma)l_3^-\dot{\mathcal U}_-=0,
       \qquad \dot s=0,\\
 &\dot J=(\mathcal A_+-\sigma I)\dot{\mathcal U}_+
                    -\mathcal U_+\dot s=0,\qquad
       \dot\alpha=\mathfrak mB=(\dot s,\dot J)=0 .
 \end{aligned}
\]
The last line uses \eqref{eq:fluxes} at the right equilibrium and
the moment identity in Proposition~\ref{prop:stationary}.

For the fixed-end coordinates in \eqref{eq:normalform-definition},
differentiation at their reference state gives
\begin{align*}
 D_q\eta(q_*)&=\left(\frac{R}{v_*},0,\frac{c_v}{\theta_*}\right),\\
 D_qw_e(q_*)&=(0,1,e\mathfrak h'(\theta_*))
                 -2e k_0\mathfrak h(\theta_*)D_q\eta(q_*)\\*
 &=\left(-\frac{ec_*}{\gamma},1,\frac{eR}{v_*c_*}\right),\\
 \det D_qN_{q_+}(q_+)
 &=2\left(\frac{c_+c_v}{\gamma\theta_+}
                   +\frac{R^2}{v_+^2c_+}\right)
       =\frac{2R(c_v+R)}{v_+^2c_+}>0,\\
 D_qw_+(q_-)&=\left(-\frac{c_-}{\gamma},1,\frac{R}{v_-c_-}\right)
                 =-2c_-\,l_3^-D_q\mathcal U(q_-).
\end{align*}
At $\alpha=0$, the derivative of
$\alpha\mapsto \bigl(N_{q_+}(q_{\alpha,+}),[N_{q_-}(q_{\alpha,-})]_+\bigr)$
is the conservative trace derivative followed by multiplication
of its right component by $D_qN_{q_+}(q_+)(D_q\mathcal U(q_+))^{-1}$
and its left component by $-2c_-$. Both factors are invertible since
$\det D_q\mathcal U(q_+)=c_v>0$ and $\det D_qN_{q_+}(q_+)>0$.
The inverse function theorem gives a local inverse. The mean value
theorem for this inverse gives, for small $\alpha$,
\[
 |\alpha|\le C_\delta\left(
     |N_{q_+}(q_{\alpha,+})|+|[N_{q_-}(q_{\alpha,-})]_+|\right).
\]
\end{proof}

\subsection{Two levels of spatial recovery}

By \eqref{eq:global-regularity}, $Z(t)\in H^2(\R)$, so $Z_\pm(t)=0$.
The moments and profile in \eqref{eq:st-moment-bounds}--\eqref{eq:st-family}
give
\begin{equation}\label{eq:st-residual}
 \begin{aligned}
 \alpha(t)&=\mathfrak mZ(t),\qquad |\alpha(t)|\le C_\delta\|Z(t)\|_\infty,\\
 r(t)&=Z(t)-\mathcal R(\alpha(t))=q(t)-q_{\alpha(t)},
                   \qquad r_\pm(t)=-\mathcal R(\alpha(t))_\pm,\\
 \mathfrak mr(t)&=\mathfrak mZ(t)-\alpha(t)=0,\qquad
                   \ell(r(t))=(\mathfrak mr(t))_0=0 .
 \end{aligned}
\end{equation}
In the coordinate \eqref{eq:moving-frame}, the spatial operator
for \eqref{eq:NSF} at speed $b$ is
\begin{equation}\label{eq:st-spatial-operator}
 \mathscr N_b(v,u,\theta)=
 \begin{pmatrix}
  bv'+u'\\
  bu'-p'+\mu\left(\frac{u'}{v}\right)'\\
  b\theta'-\frac{p}{c_v}u'+\frac{\kappa}{c_v}
          \left(\frac{\theta'}{v}\right)'+\frac{\mu(u')^2}{c_vv}
 \end{pmatrix}.
\end{equation}
At a fixed time, $F=Z_t=q_t=\mathscr N_b(q)$, and
\eqref{eq:st-family} and \eqref{eq:st-residual} give
\[
 \begin{aligned}
 s_\alpha&=\sigma+\ell(\mathcal R(\alpha))
          =\sigma+\ell(Z)=b,\qquad \mathscr N_b(q_\alpha)=0,\\
 \sigma+\ell(q_\alpha+\rho r-\bar q)
 &=s_\alpha+\rho\ell(r)=b,\qquad 0\le\rho\le1 .
 \end{aligned}
\]
The segment lies in the positive state region for small $\|Z\|_{H^2}$.
Taylor's formula at fixed $b$ gives
\begin{equation}\label{eq:st-residual-equation}
 \begin{aligned}
 F&=\mathscr N_b(q)-\mathscr N_b(q_\alpha)
      =\mathscr L_\alpha r+\mathscr E_\alpha(r),\\
 \mathscr L_\alpha&=D_q\mathscr N_b(q_\alpha),\qquad
 \mathscr E_\alpha(r)=\int_0^1(1-\rho)
     D_q^2\mathscr N_b(q_\alpha+\rho r)[r,r]\dd\rho .
 \end{aligned}
\end{equation}

Define $r_{\rm par}=(r_2,r_3)$ and
\begin{equation}\label{eq:st-X-norms}
\begin{aligned}
 X_1(r)&=m_\omega(r)+\|r'\|_2+\|r_{\rm par}''\|_2,\\
 X_2(r)&=m_\omega(r)+m_\omega(r')+\|r''\|_2
                                      +\|r_{\rm par}'''\|_2 .
\end{aligned}
\end{equation}

By \eqref{eq:global-regularity} and \eqref{eq:st-parameter-tails},
$r$ is bounded and $r'\in H^1$ at each time. Also
$r_{\rm par}'''\in L^2$ almost everywhere in time. Since
$m_\omega(r)\le\|\sqrt{\omega}\|_2\|r\|_\infty$,
$X_1(r)<\infty$ at each time and $X_2(r)<\infty$ almost everywhere.
At a time with $F'\in L^2$, differentiate the last two rows of
\eqref{eq:st-spatial-operator} in distributions to get
$\left(\frac{u'}{v}\right)'',\left(\frac{\theta'}{v}\right)''\in L^2$.
The product rule for $u'=v\frac{u'}{v}$ and
$\theta'=v\frac{\theta'}{v}$ gives
\[
 \begin{aligned}
 \|(u''',\theta''')\|_2
 &\le C_\delta\left(
   \left\|\left(\left(\frac{u'}{v}\right)'',
               \left(\frac{\theta'}{v}\right)''\right)\right\|_2
                       +\|q'\|_2+\|q''\|_2\right)\\
 &\le C_\delta\left(\|F'\|_2+\|q'\|_2+\|q''\|_2
   +\|p''\|_2+\|(pu')'\|_2+\left\|\left(\frac{(u')^2}{v}\right)'\right\|_2\right)\\
 &\le C_\delta\bigl(\|F'\|_2+\|q'\|_2+\|q''\|_2\bigr).
 \end{aligned}
\]
The last bound follows from $q'\in L^\infty\cap L^2$, $q''\in L^2$
and the positive lower bound of $v$. Thus $X_2(r)<\infty$ at each
time with $F'\in L^2$.

\begin{proposition}[Spatial recovery]\label{prop:spatial}
For $\|Z\|_{H^2}$ sufficiently small, the residual in
\eqref{eq:st-residual} satisfies
\begin{equation}\label{eq:st-spatial-estimates}
 X_1(r)\le C_\delta\|F\|_2,\qquad
 X_2(r)\le C_\delta\bigl(\|F'\|_2+m_\omega(F)\bigr).
\end{equation}
The second inequality holds whenever its right-hand side is finite.
This holds at almost every time for the strong solution.
\end{proposition}

\begin{proof}
Write $q_\alpha=(V_\alpha,U_\alpha,\Theta_\alpha)$.
Differentiating the parabolic rows of \eqref{eq:st-spatial-operator}
at fixed $b$ gives
\[
 \begin{aligned}
 (\mathscr L_\alpha r)_2
 &=b r_2'
   -\left(\frac{Rr_3}{V_\alpha}
      -\frac{R\Theta_\alpha}{V_\alpha^2}r_1\right)'
   +\mu\left(\frac{r_2'}{V_\alpha}
      -\frac{U_\alpha'}{V_\alpha^2}r_1\right)',\\
 c_v(\mathscr L_\alpha r)_3
 &=c_v b r_3'-\frac{R\Theta_\alpha}{V_\alpha}r_2'
      -U_\alpha'\left(\frac{Rr_3}{V_\alpha}
                     -\frac{R\Theta_\alpha}{V_\alpha^2}r_1\right)\\
 &\quad+\kappa\left(\frac{r_3'}{V_\alpha}
                      -\frac{\Theta_\alpha'}{V_\alpha^2}r_1\right)'
      +\mu\left(\frac{2U_\alpha'}{V_\alpha}r_2'
                      -\frac{(U_\alpha')^2}{V_\alpha^2}r_1\right).
 \end{aligned}
\]
Thus
\[
 \mathscr L_\alpha r
 =C_\alpha(\xi)r'+D_\alpha(\xi)r''+E_\alpha(\xi)r,\qquad
 D_\alpha=\operatorname{diag}
      \left(0,\frac{\mu}{V_\alpha},\frac{\kappa}{c_vV_\alpha}\right).
\]
The second derivatives in $\mathscr L_\alpha r$ and the third
derivatives in $(\mathscr L_\alpha r)'$ act on $r_{\rm par}$.
By \eqref{eq:profilebounds} and \eqref{eq:st-parameter-tails},
$C_\alpha,D_\alpha$ approach their end values exponentially.
The same bounds give exponential decay of $E_\alpha$ and of the
spatial derivatives of the coefficients. For small $\alpha$,
\eqref{eq:st-family} gives $b=s_\alpha\ge\frac{\sigma}{2}$.
The mass row gives
\begin{equation}\label{eq:st-linear-mass}
 \begin{aligned}
 (\mathscr L_\alpha r)_1&=br_1'+r_2',\\
 \|\partial_\xi^j r_1\|_2
 &=\frac{1}{b}\|\partial_\xi^{j-1}(\mathscr L_\alpha r)_1
                       -\partial_\xi^j r_2\|_2\\
 &\le\frac{2}{\sigma}\left(\|\partial_\xi^{j-1}(\mathscr L_\alpha r)_1\|_2
           +\|\partial_\xi^j r_2\|_2\right),\qquad j=1,2 .
 \end{aligned}
\end{equation}
For the residual in \eqref{eq:st-residual-equation},
$(\mathscr E_\alpha(r))_1=0$ and $W=(\mathscr L_\alpha r)_1$.

At either end state $(v_*,u_*,\theta_*)$, the limiting matrices are
\[
 C_*=\begin{pmatrix}
  b&1&0\\[1mm]
  \frac{p_*}{v_*}&b&-\frac{R}{v_*}\\[2mm]
  0&-\frac{p_*}{c_v}&b
 \end{pmatrix},\qquad
 D_*=\operatorname{diag}\left(0,\frac{\mu}{v_*},\frac{\kappa}{c_vv_*}\right),
 \qquad p_*=\frac{R\theta_*}{v_*}.
\]
Expanding in the first row and using $c_*^2=\frac{\gamma p_*}{v_*}$
gives
\begin{equation}\label{eq:st-tail-determinant}
 \begin{aligned}
 \det(C_*+ikD_*)
 &=b\left[\left(b+\frac{ik\mu}{v_*}\right)
          \left(b+\frac{ik\kappa}{c_vv_*}\right)-\frac{Rp_*}{c_vv_*}\right]
           -\frac{p_*}{v_*}\left(b+\frac{ik\kappa}{c_vv_*}\right)\\
 &=b(b^2-c_*^2)-\frac{b\mu\kappa}{c_vv_*^2}k^2
   +ik\left[\frac{b^2\mu}{v_*}
      +\frac{\kappa}{c_vv_*}\left(b^2-\frac{p_*}{v_*}\right)\right].
 \end{aligned}
\end{equation}
The strict inequalities in \eqref{eq:lax} hold for small $\alpha$.
At the left end, $b^2<c_*^2$ makes the real part negative for all
real $k$. At the right end, $b^2>c_*^2>\frac{p_*}{v_*}$ makes the
imaginary part nonzero for $k\ne0$ and gives $\det C_*>0$ at $k=0$.

For $z\in\mathbb C^3$ and $z_{\rm par}=(z_2,z_3)^{\mathsf T}$,
the first row and the positive entries of $D_*$ give
\[
 \begin{aligned}
 |z|&\le\frac{1}{b}|(C_*+ikD_*)z|
                      +\left(1+\frac{1}{b}\right)|z_{\rm par}|,\\
 |k|\,|z_{\rm par}|
 &\le C_\delta|ikD_*z|
 \le C_\delta\bigl(|(C_*+ikD_*)z|+|z|\bigr)
 \le C_\delta\bigl(|(C_*+ikD_*)z|+|z_{\rm par}|\bigr).
 \end{aligned}
\]
For large $|k|$, the left side absorbs the last term on the right.
On a bounded frequency interval, \eqref{eq:st-tail-determinant}
gives a uniformly bounded inverse for $\alpha$ in a small closed
ball. Thus for every real $k$,
\begin{equation}\label{eq:st-tail-multiplier}
 |z|+|k|\,|z_{\rm par}|
       \le C_\delta|(C_*+ikD_*)z|.
\end{equation}
By Plancherel's identity, for $g\in L^2$ with $g_{\rm par}'\in L^2$,
\[
 \|g\|_2+\|g_{\rm par}'\|_2
 \le C_\delta\|C_*g+D_*g'\|_2 .
\]

For $L>1$, take translates of fixed smooth cutoffs
$0\le\chi_\pm\le1$ with $\chi_\pm=0$ on $\pm\xi\le L$ and
$\chi_\pm=1$ on $\pm\xi\ge L+1$. Put $K=[-L-1,L+1]$ and denote the end
matrices by $C_\pm,D_\pm$. By \eqref{eq:profilebounds} and
\eqref{eq:st-parameter-tails},
\[
 |C_\alpha-C_\pm|+|D_\alpha-D_\pm|+\frac{|E_\alpha|}{\sqrt{\omega}}
 \le C_\delta e^{-\eta_0 L},\qquad \pm\xi\ge L .
\]
The product rule and \eqref{eq:st-tail-multiplier} give
\[
 \begin{aligned}
 &C_\pm(\chi_\pm r')+D_\pm(\chi_\pm r')'\\
 &\quad=\chi_\pm\bigl[\mathscr L_\alpha r+(C_\pm-C_\alpha)r'
              +(D_\pm-D_\alpha)r''-E_\alpha r\bigr]
              +D_\pm\chi_\pm' r',\\
 &\|\chi_\pm r'\|_2+\|(\chi_\pm r_{\rm par}')'\|_2\\
 &\quad\le C_\delta\left(\|\mathscr L_\alpha r\|_2
                +\|r_{\rm par}'\|_{L^2(K)}
                +e^{-\eta_0 L}X_1(r)\right).
 \end{aligned}
\]
Here $D_\pm-D_\alpha$ and $D_\pm\chi_\pm'$ act only on the parabolic
components, and $\chi_\pm'$ is supported in $K$.

On $K$, eliminate $r_1'$ from the parabolic rows using
\eqref{eq:st-linear-mass}. The positive diagonal entries of
$D_\alpha$ in those rows give
\[
 \begin{aligned}
 \|r'\|_{L^2(K)}+\|r_{\rm par}''\|_{L^2(K)}
 &\le C_\delta\left(\|\mathscr L_\alpha r\|_{L^2(K)}
             +\|r_{\rm par}'\|_{L^2(K)}+\|r\|_{L^2(K)}\right),\\
 m_\omega(r)^2
 &\le\int_{\R}\omega(\xi)\left[
      \|r\|_{L^2(-1,1)}^2+2(1+|\xi|)\|r'\|_2^2\right]\dd\xi
 \le C_\delta\bigl(\|r\|_{L^2(K)}^2+\|r'\|_2^2\bigr).
 \end{aligned}
\]
The second line follows by averaging as in \eqref{eq:anchored-weight}.
The corresponding cutoff equals one on $\R\setminus K$.
Add the estimates on and outside $K$ and interpolate
$\|r_{\rm par}'\|_{L^2(K)}$ between
$\|r_{\rm par}\|_{L^2(K)}$ and $\|r_{\rm par}''\|_{L^2(K)}$ to obtain
\[
 \begin{aligned}
 X_1(r)
 &\le C_\delta\left(\|\mathscr L_\alpha r\|_2
          +\|r\|_{L^2(K)}+\|r_{\rm par}'\|_{L^2(K)}\right)
          +C_\delta e^{-\eta_0 L}X_1(r)\\
 &\le C_{\delta,\epsilon_0,K}\left(\|\mathscr L_\alpha r\|_2
                        +\|r\|_{L^2(K)}\right)
     +C_\delta\left(\epsilon_0\|r_{\rm par}''\|_{L^2(K)}
                              +e^{-\eta_0 L}X_1(r)\right)\\
 &\le C_{\delta,\epsilon_0,K}\left(\|\mathscr L_\alpha r\|_2
                  +\|r\|_{L^2(K)}\right)
        +C_\delta(\epsilon_0+e^{-\eta_0 L})X_1(r).
 \end{aligned}
\]
The constant $C_\delta$ multiplying $\epsilon_0+e^{-\eta_0 L}$
is independent of $L$ and $\epsilon_0$. Take $L$ large and
$\epsilon_0$ small to bound the last term by $\frac{1}{2}X_1(r)$.
Then
\begin{equation}\label{eq:st-garding-one}
 X_1(r)\le C_\delta
       \bigl(\|\mathscr L_\alpha r\|_2+\|r\|_{L^2(K)}\bigr).
\end{equation}

Take $L$ sufficiently large and enlarge $K$ to contain $[-L,L]$.
The estimate \eqref{eq:st-garding-one} holds on the larger $K$
with the same constant. The product rule and the coefficient bounds give
\[
 \begin{aligned}
 \mathscr L_\alpha(r')
 &= (\mathscr L_\alpha r)'
      -D_\alpha' r''-C_\alpha' r'-E_\alpha' r,\\
 |D_\alpha'|+\frac{|C_\alpha'|+|E_\alpha'|}{\sqrt{\omega}}
 &\le C_\delta e^{-\eta_0 L},\qquad |\xi|\ge L,\\
 \|\mathscr L_\alpha(r')\|_2
 &\le \|(\mathscr L_\alpha r)'\|_2
       +C_\delta\bigl(\|r\|_{H^1(K)}+\|r_{\rm par}''\|_{L^2(K)}\bigr)
       +C_\delta e^{-\eta_0 L}X_2(r).
 \end{aligned}
\]
In $D_\alpha'r''$, only the components $r_{\rm par}''$ occur.
Use the fundamental theorem of calculus for $r'$ and the integral
Taylor formula for $r$. Averaging over $y\in[-1,1]$ gives
\[
 \begin{aligned}
 |r'(\xi)|+\frac{|r(\xi)|}{1+|\xi|}
 &\le C\int_{-1}^1\bigl(|r(y)|+|r'(y)|\bigr)\dd y
       +C\sqrt{1+|\xi|}\,\|r''\|_2\\
 &\le C\bigl(\|r\|_{H^1(-1,1)}+\sqrt{1+|\xi|}\,\|r''\|_2\bigr)
 \le C_\delta\sqrt{1+|\xi|}\,X_2(r),\\
 m_\omega(r)+m_\omega(r')
 &\le C_\delta\bigl(\|r\|_{H^1(K)}+\|r''\|_2\bigr).
 \end{aligned}
\]
We used the positive lower bound of $\omega$ on $[-1,1]$ and
$\int_\R(1+|\xi|)^3\omega(\xi)\dd\xi<\infty$.
Since $X_1(r')\le X_2(r)$, we can apply \eqref{eq:st-garding-one} to
$r'$. Interpolation on $K$ then gives
\[
 \begin{aligned}
 X_2(r)=X_1(r')+m_\omega(r)
 &\le C_\delta\bigl(X_1(r')+\|r\|_{H^1(K)}\bigr)\\
 &\le C_\delta\left(\|(\mathscr L_\alpha r)'\|_2+\|r\|_{H^1(K)}
                   +\|r_{\rm par}''\|_{L^2(K)}\right)
                   +C_\delta e^{-\eta_0 L}X_2(r)\\
 &\le C_{\delta,\epsilon_0,K}\left(
       \|(\mathscr L_\alpha r)'\|_2+\|r\|_{H^1(K)}\right)
        +C_\delta\left(\epsilon_0\|r_{\rm par}'''\|_{L^2(K)}
                              +e^{-\eta_0 L}X_2(r)\right)\\
 &\le C_{\delta,\epsilon_0,K}\left(
       \|(\mathscr L_\alpha r)'\|_2+\|r\|_{H^1(K)}\right)
        +C_\delta(\epsilon_0+e^{-\eta_0 L})X_2(r).
 \end{aligned}
\]
Here $C_\delta$ in the last term is independent of $L$ and $\epsilon_0$.
Take $L$ large and $\epsilon_0$ small so that
$C_\delta(\epsilon_0+e^{-\eta_0 L})\le\frac{1}{2}$. We obtain
\begin{equation}\label{eq:st-garding-two}
 X_2(r)\le C_\delta
       \bigl(\|(\mathscr L_\alpha r)'\|_2+\|r\|_{H^1(K)}\bigr).
\end{equation}

Consider a distributional solution $B\in H^1_{\mathrm{loc}}$ of
$\mathscr L_\alpha B=0$ with $|B(\xi)|\le C(1+|\xi|)^{\frac{3}{2}}$.
The variations $\dot J_i$ of \eqref{eq:fluxes} at fixed $s=b$ are
locally in $L^2$. Linearizing the conservation laws at the stationary
profile $q_\alpha$ gives, in distributions,
\[
 \begin{pmatrix}\dot J_0'\\ \dot J_1'\\ \dot J_2'\end{pmatrix}
 =-\begin{pmatrix}
     (\mathscr L_\alpha B)_1\\
     (\mathscr L_\alpha B)_2\\
     c_v(\mathscr L_\alpha B)_3+U_\alpha(\mathscr L_\alpha B)_2
    \end{pmatrix}=0 .
\]
Each $\dot J_i$ is therefore constant. Subtract the profile variation
with these fluxes at fixed speed
\[
 (f,g,h)=B-\left.D_J
       \bigl[\widehat q_{s,J}(\cdot+\tau(s,J))\bigr]
       \right|_{(s,J)=(s_\alpha,J_\alpha)}\dot J .
\]
The difference has zero flux variation. Linearizing
\eqref{eq:st-reduced} gives $g=-bf$ and the same first-order system
for $(f,h)^{\mathsf T}$ as for $(V_\alpha',\Theta_\alpha')^{\mathsf T}$.
The equations hold in distributions. Smooth coefficients give a
classical solution.

On the left unstable graph, $V_\alpha'\ne0$ and
$\frac{\Theta_\alpha'}{V_\alpha'}$ tends to the slope of the unstable
eigenvector. The derivative of the determinant of the two solutions
is the trace of the coefficient matrix times that determinant.
Using the first row of the variational equation for $V_\alpha'$
and \eqref{eq:st-reduced}, we obtain
\[
 \begin{aligned}
 \left(h-\frac{\Theta_\alpha'}{V_\alpha'}f\right)'
 &=\frac{(V_\alpha'h-\Theta_\alpha'f)'}{V_\alpha'}
   -\frac{V_\alpha''}{V_\alpha'}
       \left(h-\frac{\Theta_\alpha'}{V_\alpha'}f\right)\\
 &=\left(-\frac{bc_vV_\alpha}{\kappa}
        +\frac{R\Theta_\alpha'}{\mu bV_\alpha'}\right)
        \left(h-\frac{\Theta_\alpha'}{V_\alpha'}f\right).
 \end{aligned}
\]
At the left equilibrium the coefficient tends to the trace of
\eqref{eq:st-spatial-matrix} minus its positive eigenvalue. This is
its negative eigenvalue. For sufficiently large $\xi_0$, the
coefficient is at most $-c_\delta$ on $\xi\le-\xi_0$.
Lemma~\ref{lem:st-raw-family} bounds the subtracted profile variation.
For $\rho\le\xi\le-\xi_0$,
\[
 \begin{aligned}
 \left|\left(h-\frac{\Theta_\alpha'}{V_\alpha'}f\right)(\xi)\right|
 &\le e^{-c_\delta(\xi-\rho)}
     \left|\left(h-\frac{\Theta_\alpha'}{V_\alpha'}f\right)(\rho)\right|\\
 &\le C e^{-c_\delta(\xi-\rho)}(1+|\rho|)^{\frac{3}{2}} .
 \end{aligned}
\]
Letting $\rho$ tend to $-\infty$ gives
$h=\frac{\Theta_\alpha'}{V_\alpha'}f$ there. By uniqueness for the
linear system and $g=-bf$, $(f,g,h)$ is a constant multiple
of $q_\alpha'$ on $\R$.

Differentiating \eqref{eq:st-phase} at fixed $s$ gives
$\ell(f,g,h)=\ell(B)$. By \eqref{eq:st-phase-nondeg},
$\ell(q_\alpha')\le-\frac{|\ell(\bar q')|}{2}<0$.
If $\mathfrak mB=0$, then $\ell(B)=0$ and hence $(f,g,h)=0$.
It follows that
\begin{equation}\label{eq:st-fixed-speed-moments}
 \begin{aligned}
 0=\mathfrak mB
 &=\left.D_{(s,J)}
       \left[\mathfrak m\bigl(\widehat q_{s,J}(\cdot+\tau(s,J))-\bar q\bigr)
       \right]\right|_{(s,J)=(s_\alpha,J_\alpha)}(0,\dot J).
 \end{aligned}
\end{equation}
The derivative is invertible on a small closed ball of $\alpha$
by Proposition~\ref{prop:stationary}. Thus $\dot J=0$ and $B=0$.

If either estimate in \eqref{eq:st-linear-recovery} fails, there are
$\alpha_n$ in the closed parameter ball above and $\mathfrak mr_n=0$
such that one of the following holds,
\[
 \begin{aligned}
 X_1(r_n)&=1,\qquad
     \lim_{n\to\infty}\|\mathscr L_{\alpha_n}r_n\|_2=0,\\
 X_2(r_n)&=1,\qquad
     \lim_{n\to\infty}\bigl(\|(\mathscr L_{\alpha_n}r_n)'\|_2
                              +m_\omega(\mathscr L_{\alpha_n}r_n)\bigr)=0.
 \end{aligned}
\]
The averaging estimates and the positive minimum of $\omega$ on
each compact interval $I$ give
\[
 \begin{aligned}
 \|r\|_{H^1(I)}+\sup_{\xi\in\R}
       \frac{|r(\xi)|}{(1+|\xi|)^{\frac{1}{2}}}
 &\le C_{\delta,I}\bigl(m_\omega(r)+\|r'\|_2\bigr)
 \le C_{\delta,I}X_1(r),\\
 \|r\|_{H^2(I)}+\sup_{\xi\in\R}
       \frac{|r(\xi)|}{(1+|\xi|)^{\frac{3}{2}}}
 &\le C_{\delta,I}\bigl(m_\omega(r)+m_\omega(r')+\|r''\|_2\bigr)
 \le C_{\delta,I}X_2(r).
 \end{aligned}
\]
Pass to a subsequence with $\alpha_n$ converging to $\alpha_*$ in
the closed ball. By compactness on bounded intervals, $r_n$ converges
along a diagonal subsequence to $r_*\in H^1_{\mathrm{loc}}$.
The convergence is strong in $L^2_{\mathrm{loc}}$ in the first
case and in $H^1_{\mathrm{loc}}$ in the second. Pass to a further
subsequence for convergence almost everywhere. Since
$m_\omega(f)\le\|f\|_2$, $m_\omega(\mathscr L_{\alpha_n}r_n)$ tends
to zero in both cases. For $\varphi\in C_c^\infty(I;\R^3)$, local
convergence and smooth dependence of the coefficients give
\[
 \begin{aligned}
 \langle\mathscr L_{\alpha_*}r_*,\varphi\rangle
 &=\lim_{n\to\infty}\int_\R r_n\cdot
      \bigl[(D_{\alpha_n}^{\mathsf T}\varphi)''
       -(C_{\alpha_n}^{\mathsf T}\varphi)'
       +E_{\alpha_n}^{\mathsf T}\varphi\bigr]\dd\xi\\
 &=\lim_{n\to\infty}\langle\mathscr L_{\alpha_n}r_n,\varphi\rangle=0,
 \end{aligned}
\]
where $\|\mathscr L_{\alpha_n}r_n\|_{L^2(I)}
\le C_{\delta,I}m_\omega(\mathscr L_{\alpha_n}r_n)$.
The growth bounds give
$\omega|r_*-r_n|\le C_\delta\omega(1+|\xi|)^{\frac{3}{2}}$.
By \eqref{eq:st-moment-bounds} and dominated convergence,
\[
 |\mathfrak mr_*|
 \le C_\delta\lim_{n\to\infty}\int_{\R}\omega|r_*-r_n|\dd\xi=0.
\]
The kernel calculation gives $r_*=0$. Apply
\eqref{eq:st-garding-one} in the first case and
\eqref{eq:st-garding-two} in the second to get
\[
 \begin{aligned}
 1&\le C_\delta\lim_{n\to\infty}
         \bigl(\|\mathscr L_{\alpha_n}r_n\|_2+\|r_n\|_{L^2(K)}\bigr)=0,\\
 1&\le C_\delta\lim_{n\to\infty}
         \bigl(\|(\mathscr L_{\alpha_n}r_n)'\|_2+\|r_n\|_{H^1(K)}\bigr)=0 .
 \end{aligned}
\]
It follows that
\begin{equation}\label{eq:st-linear-recovery}
 \begin{aligned}
 X_1(r)&\le C_\delta\|\mathscr L_\alpha r\|_2,\\
 X_2(r)&\le C_\delta\bigl(\|(\mathscr L_\alpha r)'\|_2
                               +m_\omega(\mathscr L_\alpha r)\bigr),
 \qquad \mathfrak mr=0 .
 \end{aligned}
\end{equation}

Sobolev embedding and \eqref{eq:global-H2}, \eqref{eq:st-residual}
and \eqref{eq:st-family-bounds} give
$\|r\|_\infty+\|r'\|_\infty\le C_\delta(\|Z\|_{H^2}+|\alpha|)
\le C_\delta\eps$. By \eqref{eq:profilebounds} and
\eqref{eq:st-parameter-tails},
$\sum_{j=1}^3|\partial_\xi^j q_\alpha|\le C_\delta\omega$.

Since $v=V_\alpha+r_1$, we have
$\frac{r_1}{vV_\alpha}=\frac{1}{V_\alpha}-\frac{1}{v}$.
Subtracting the linear terms in the momentum and temperature rows
of \eqref{eq:st-residual-equation} gives
\[
 \begin{aligned}
 (\mathscr E_\alpha(r))_2
 &=\left[\frac{r_1}{vV_\alpha}
      \left(R\left(r_3-\frac{\Theta_\alpha}{V_\alpha}r_1\right)
        -\mu\left(r_2'-\frac{U_\alpha'}{V_\alpha}r_1\right)\right)\right]',\\
 c_v(\mathscr E_\alpha(r))_3
 &=-\kappa\left[\frac{r_1}{vV_\alpha}
          \left(r_3'-\frac{\Theta_\alpha'}{V_\alpha}r_1\right)\right]'\\
 &\quad-\frac{R}{v}\left(r_3-\frac{\Theta_\alpha}{V_\alpha}r_1\right)
                   \left(r_2'-\frac{U_\alpha'}{V_\alpha}r_1\right)
        +\frac{\mu}{v}\left(r_2'-\frac{U_\alpha'}{V_\alpha}r_1\right)^2 .
 \end{aligned}
\]
The positive lower bounds for $v,V_\alpha$, the profile bounds and
$\omega|r||r'|\le\frac{1}{2}(\omega|r|^2+|r'|^2)$ give by the product rule
\[
 \begin{aligned}
 \left|\left(\frac{r_1}{vV_\alpha}\right)'\right|
 &\le C_\delta\bigl(|r'|+\omega|r|\bigr),\\
 \left|\left(\frac{r_1}{vV_\alpha}\right)''\right|
 &\le C_\delta\bigl(|r''|+|r'|^2+\omega(|r|+|r'|)\bigr),\\
 |\mathscr E_\alpha(r)|
 &\le C_\delta\bigl(
      |r||r_{\rm par}''|+|r'|^2+|r||r'|+\omega|r|^2\bigr),\\
 |(\mathscr E_\alpha(r))'|
 &\le C_\delta\bigl(
      |r||r_{\rm par}'''|+(|r|+|r'|)|r''|+|r'|^2+\omega|r|^2\bigr).
 \end{aligned}
\]
Taking norms and using $0<\omega\le1$ and \eqref{eq:st-GN-bounded}, we obtain
\begin{equation}\label{eq:st-nonlinear-error}
 \begin{aligned}
 &\|\mathscr E_\alpha(r)\|_2\\
 &\quad\le C_\delta\left[\|r\|_\infty
   \bigl(\|r_{\rm par}''\|_2+\|r'\|_2+m_\omega(r)\bigr)
                        +\|r'\|_\infty\|r'\|_2\right]
 \le C_\delta\eps X_1(r),\\
 &\|(\mathscr E_\alpha(r))'\|_2+m_\omega(\mathscr E_\alpha(r))\\
 &\quad\le C_\delta\left[\|r\|_\infty
   \bigl(\|r_{\rm par}'''\|_2+m_\omega(r)\bigr)
   +(\|r\|_\infty+\|r'\|_\infty)
           \bigl(\|r''\|_2+m_\omega(r')\bigr)+\|r'\|_4^2\right]\\
 &\quad\le C_\delta(\|r\|_\infty+\|r'\|_\infty)X_2(r)
 \le C_\delta\eps X_2(r).
 \end{aligned}
\end{equation}
Since $r'\in H^1(\R)$ and $r\in L^\infty$, the quantity
$(r\cdot r')|r'|^2$ vanishes at both infinities. Integration by parts
and the Cauchy--Schwarz inequality give
\begin{equation}\label{eq:st-GN-bounded}
 \begin{aligned}
 \|r'\|_4^4
 &=-\int_{\R}
       \bigl((r\cdot r'')|r'|^2+2(r\cdot r')(r'\cdot r'')\bigr)\dd\xi\\
 &\le 3\|r\|_\infty\|r'\|_4^2\|r''\|_2,\\
 \|r'\|_4^2&\le 3\|r\|_\infty\|r''\|_2 .
 \end{aligned}
\end{equation}

By \eqref{eq:st-linear-recovery}, \eqref{eq:st-residual-equation}
and \eqref{eq:st-nonlinear-error},
\[
 \begin{aligned}
 X_1(r)
 &\le C_\delta\bigl(\|F\|_2+\|\mathscr E_\alpha(r)\|_2\bigr)
 \le C_\delta\|F\|_2+C_\delta\eps X_1(r),\\
 X_2(r)
 &\le C_\delta\bigl(\|F'\|_2+m_\omega(F)
            +\|(\mathscr E_\alpha(r))'\|_2+m_\omega(\mathscr E_\alpha(r))\bigr)\\
 &\le C_\delta\bigl(\|F'\|_2+m_\omega(F)\bigr)
        +C_\delta\eps X_2(r).
 \end{aligned}
\]
Taking $C_\delta\eps\le\frac{1}{2}$ gives \eqref{eq:st-spatial-estimates}.
\end{proof}

\subsection{The window norm and residual gradient potentials}

The mass row \eqref{eq:st-linear-mass} gives $W=br_1'+r_2'$ and
$br_1''=W'-r_2''$. Since $b\ge\frac{\sigma}{2}$,
\eqref{eq:st-spatial-estimates} and \eqref{eq:time-energy-equivalence}
give, at every time,
\begin{equation}\label{eq:st-endpoint-recovery}
 \begin{aligned}
 m_\omega(r)+\|r'\|_{H^1}+\|r_{\rm par}''\|_2
 &\le X_1(r)+\|r''\|_2\\
 &\le X_1(r)+\frac{\|W'\|_2+\|r_2''\|_2}{b}
                         +\|r_{\rm par}''\|_2\\
 &\le C_\delta\bigl(\|F\|_2+\|W'\|_2\bigr)
               \le C_\delta\sqrt{\mathcal T}.
 \end{aligned}
\end{equation}
Averaging the fundamental theorem of calculus over $[-1,1]$ and
using the lower bound of $\omega$ there give
\[
 \begin{aligned}
 |r(\xi)|
 &\le\frac{1}{2}\int_{-1}^1
             \bigl(|r(y)|+\sqrt{|\xi-y|}\,\|r'\|_2\bigr)\dd y
 \le C_\delta m_\omega(r)+\sqrt{1+|\xi|}\,\|r'\|_2\\
 &\le C_\delta\sqrt{1+|\xi|}\,X_1(r)
 \le C_\delta\sqrt{1+|\xi|}\,\|F\|_2
 \le C_\delta\sqrt{1+|\xi|}\sqrt{\mathcal T}.
 \end{aligned}
\]
For $f\in H^1_{\mathrm{loc}}(\R)$ with $f'\in L^2$, the same averaging
argument as in \eqref{eq:anchored-weight} gives
\begin{equation}\label{eq:st-anchor}
 \begin{aligned}
 m_\omega(f)^2
 &\le C_\delta\bigl(\|f\|_{L^2(-1,1)}^2+\|f'\|_2^2\bigr)\\
 &\le C_\delta\left(\|f'\|_2^2+\int_{\R}V'|f|^2\dd\xi\right),
 \end{aligned}
\end{equation}
where the last inequality uses the positive minimum of $V'$ on
$[-1,1]$. By \eqref{eq:profilebounds} and \eqref{eq:feedbackbounds},
$V'\le C_\delta\omega$ and
$|\ell(F)|^2\le C_\delta\int_{\R}V'|F|^2\dd\xi$.
Together with \eqref{eq:st-anchor}, this gives
\[
 \begin{aligned}
 c_\delta\Gamma(t)^2
 &\le \|F'\|_2^2+\int_{\R}V'|F|^2\dd\xi+|\ell(F)|^2\\
 &\le \|F'\|_2^2+C_\delta\int_{\R}V'|F|^2\dd\xi
       \le C_\delta\Gamma(t)^2.
 \end{aligned}
\]

For $I=[s,t]$ with $0\le s<t$, define
\begin{equation}\label{eq:st-h}
 \begin{aligned}
 h(t)={}&m_\omega(r)+m_\omega(r')+\|r'\|_\infty
       +\|r_{\rm par}''\|_\infty+\|F\|_\infty+|\alpha'(t)|,\\
 H_I^2={}&\int_I h(\rho)^2\dd\rho .
 \end{aligned}
\end{equation}

\begin{lemma}\label{lem:st-window-h}
At almost every time,
\begin{equation}\label{eq:st-h-product}
 h(t)^2\le C_\delta\|F(t)\|_2\Gamma(t).
\end{equation}
Consequently, for $\beta_I$ and $Q_I$ as in \eqref{eq:time-window-notation},
\begin{equation}\label{eq:st-H-window}
 H_I^2\le C_\delta\beta_IQ_I.
\end{equation}
\end{lemma}

\begin{proof}
Since $m_\omega(r')\le\|r'\|_2$, we have
$m_\omega(r)+m_\omega(r')\le\min\{X_1(r),X_2(r)\}$.
The identity $Z_t=F$ in $L^2$ and \eqref{eq:st-moment-bounds} give
$\alpha'=\mathfrak mF$ and
$|\alpha'|^2\le C_\delta\|F\|_2m_\omega(F)$.
The one-dimensional Sobolev inequality for $r'$, $r_{\rm par}''$
and $F$, and \eqref{eq:st-spatial-estimates}, then give
\[
 \begin{aligned}
 h(t)^2
 &\le C_\delta\left[
       \bigl(m_\omega(r)+m_\omega(r')\bigr)^2
       +\|r'\|_2\|r''\|_2+\|r_{\rm par}''\|_2\|r_{\rm par}'''\|_2
       +\|F\|_2\bigl(\|F'\|_2+m_\omega(F)\bigr)\right]\\
 &\le C_\delta\bigl(X_1(r)X_2(r)+\|F\|_2\Gamma(t)\bigr)
 \le C_\delta\|F\|_2\Gamma(t).
 \end{aligned}
\]
Integrate \eqref{eq:st-h-product} and use $\|F\|_2^2\le C_\delta B$
from \eqref{eq:time-energy-equivalence}. The Cauchy--Schwarz inequality
gives
\[
 \begin{aligned}
 H_I^2
 &=\int_I h(\rho)^2\dd\rho
 \le C_\delta\int_I\|F(\rho)\|_2\Gamma(\rho)\dd\rho\\
 &\le C_\delta\left(\int_I B(\rho)\dd\rho\right)^{\frac{1}{2}}
              \left(\int_I\Gamma(\rho)^2\dd\rho\right)^{\frac{1}{2}}
 =C_\delta\beta_IQ_I .
 \end{aligned}
\]
\end{proof}

Let $K(\tau,\rho;x,y)$ be a nonnegative measurable kernel with
$\int_{\R}K(\tau,\rho;x,y)\dd y\le1$ for $\rho<\tau$ in $I$ and
$x\in\R$. By \eqref{eq:st-h} and \eqref{eq:st-H-window},
\begin{equation}\label{eq:st-residual-potential}
 \begin{aligned}
 &\sup_{\tau\in I,\ x\in\R}
   \int_{I\cap(-\infty,\tau)}
   \int_{\R}K(\tau,\rho;x,y)|r'(\rho,y)|^2\dd y\dd\rho\\
 &\quad\le\sup_{\tau\in I,\ x\in\R}\int_{I\cap(-\infty,\tau)}
          \|r'(\rho)\|_\infty^2
          \left(\int_{\R}K(\tau,\rho;x,y)\dd y\right)\dd\rho\\
 &\quad\le\int_I\|r'(\rho)\|_\infty^2\dd\rho
       \le H_I^2\le C_\delta\beta_IQ_I .
 \end{aligned}
\end{equation}
For the Gaussian with variance proportional to $\tau-\rho$, the
spatial integral equals one for each $\rho<\tau$. The term
$\mathcal R'(\alpha)$ in $Z'=r'+\mathcal R'(\alpha)$ belongs to the
stationary comparison profile.

\subsection{A nonlinear Nash inequality}

\begin{proposition}\label{prop:nash}
For the same small perturbations,
\begin{equation}\label{eq:st-nash}
 \|Z\|_\infty^2\le C_\delta\|Z\|_2\|F\|_2.
\end{equation}
\end{proposition}

\begin{proof}
If $\|Z\|_2=0$, then $Z=0$. If $\|F\|_2=0$,
\eqref{eq:st-spatial-estimates} gives $r=0$, so
$\mathcal R(\alpha)_\pm=Z_\pm=0$ and $\alpha=0$ by
\eqref{eq:st-end-injectivity}. Hence $Z=0$ also in this case.
Suppose $\|Z\|_2$ and $\|F\|_2$ are positive. The pointwise bound
after \eqref{eq:st-endpoint-recovery} gives
$|r(\xi)|\le C_\delta\sqrt{1+|\xi|}\,\|F\|_2$.
For $L\ge1$, average $Z=\mathcal R(\alpha)+r$ on $[L,2L]$
and $[-2L,-L]$ and use \eqref{eq:st-parameter-tails} to get
\[
 \begin{aligned}
 |\mathcal R(\alpha)_\pm|
 &\le \frac{1}{L}\int_L^{2L}
    \bigl(|Z(\pm y)|+|r(\pm y)|
       +|\mathcal R(\alpha)(\pm y)-\mathcal R(\alpha)_\pm|\bigr)\dd y\\
 &\le \frac{\|Z\|_2}{\sqrt L}+C_\delta\sqrt L\|F\|_2
                    +C_\delta|\alpha|e^{-2\eta_0 L}.
 \end{aligned}
\]
Choose $L_0\ge1$ with $C_\delta e^{-2\eta_0 L_0}\le\frac{1}{2}$.
Then \eqref{eq:st-end-injectivity} and absorption of the last term
give for all $L\ge L_0$
\begin{equation}\label{eq:st-alpha-average}
 |\alpha|\le C_\delta\left(\frac{\|Z\|_2}{\sqrt L}
                              +\sqrt L\|F\|_2\right),\qquad L\ge L_0.
\end{equation}
Take $L=\frac{\|Z\|_2}{\|F\|_2}$. When $L\ge L_0$, the two terms
in \eqref{eq:st-alpha-average} are equal. For $0<L<L_0$, use
$\alpha=\mathfrak m Z$ and $\|Z\|_2^2=L\|Z\|_2\|F\|_2$ in
\eqref{eq:st-moment-bounds}. Thus
\begin{equation}\label{eq:st-alpha-nash}
 |\alpha|^2\le C_\delta
 \begin{cases}
  \frac{\|Z\|_2^2}{L}+L\|F\|_2^2,&L\ge L_0,\\[2mm]
  \|Z\|_2^2,&0<L<L_0
 \end{cases}
 \le C_\delta\|Z\|_2\|F\|_2 .
\end{equation}
With the same $L$, averaging $Z=\mathcal R(\alpha)+r$ over
$[\xi,\xi+L]$ and using \eqref{eq:st-family-bounds},
\eqref{eq:st-spatial-estimates} and \eqref{eq:st-alpha-nash} give
\[
 \begin{aligned}
 |Z(\xi)|
 &\le\frac{1}{L}\int_\xi^{\xi+L}
       \bigl(|Z(y)|+|r(\xi)-r(y)|+|\mathcal R(\alpha)(\xi)
                                  -\mathcal R(\alpha)(y)|\bigr)\dd y\\
 &\le\frac{\|Z\|_2}{\sqrt L}+\sqrt L\|r'\|_2
                                  +2\|\mathcal R(\alpha)\|_\infty\\
 &\le\frac{\|Z\|_2}{\sqrt L}
              +C_\delta\sqrt L\|F\|_2+C_\delta|\alpha|
 \le C_\delta\sqrt{\|Z\|_2\|F\|_2}.
 \end{aligned}
\]
Taking the supremum in $\xi$ gives \eqref{eq:st-nash}.
\end{proof}

\section{Characteristic variables and positive comparison systems}
\label{sec:normalform}

\subsection{An exact form of the Euler coupling}
Let $c_p=c_v+R$, $k_0=\frac{1}{4c_p}$ and
$\mathfrak h(\theta)=\frac{2\sqrt{\gamma R\theta}}{\gamma-1}$.
For a reference state $q_*=(v_*,u_*,\theta_*)$ with $v_*,\theta_*>0$,
define
\begin{equation}\label{eq:normalform-definition}
 \begin{gathered}
 \eta=c_v\log\left(\frac{\theta}{\theta_*}\right)
         +R\log\left(\frac{v}{v_*}\right),\qquad N_{q_*}(q)=(w_+,\eta,w_-),\\
 w_e=e^{-k_0\eta}\{u-u_*+e \mathfrak h(\theta)\}
       -e \mathfrak h(\theta_*)e^{k_0\eta},\qquad e\in\{+1,-1\}.
 \end{gathered}
\end{equation}
The sum and difference of $w_+$ and $w_-$ and the entropy relation give
\[
 \begin{aligned}
 u-u_*&=\frac{e^{k_0\eta}}{2}(w_++w_-),\\
 \mathfrak h(\theta)&=\frac{e^{k_0\eta}}{2}(w_+-w_-)
                         +\mathfrak h(\theta_*)e^{2k_0\eta},\\
 v&=v_*e^{\frac{\eta}{R}}\left(\frac{\theta_*}{\theta}\right)^{\frac{c_v}{R}} .
 \end{aligned}
\]
The formula for $\mathfrak h(\theta)$ has the positive value
$\mathfrak h(\theta_*)$ at $N_{q_*}(q)=0$. It stays positive near zero,
and $\mathfrak h'(\theta)>0$ gives a smooth local inverse.

For the Euler part of \eqref{eq:NSF}, write
$D_e=\partial_t+ec\partial_x$, with $c=c(q)$ as in \eqref{eq:lax}.
Then
\[
 \begin{aligned}
 p\mathfrak h'(\theta)
 &=\frac{R}{\gamma-1}c=c_vc,\qquad
 c\mathfrak h'(\theta)=\frac{\gamma R}{(\gamma-1)v}=\frac{c_p}{v},\\
 2k_0c\mathfrak h(\theta)
 &=\frac{\gamma R}{(\gamma-1)c_p}\frac{\theta}{v}=\frac{\theta}{v}.
 \end{aligned}
\]
The Euler equations then give
\[
 \begin{aligned}
 \eta_t
 &=\frac{c_v\theta_t}{\theta}+\frac{Rv_t}{v}
   =\left(-\frac{p}{\theta}+\frac{R}{v}\right)u_x=0,\qquad
 D_e\eta=ec\eta_x,\\
 D_e(u-u_*+e\mathfrak h(\theta))
 &=-p_x+e\left(c-\frac{p\mathfrak h'(\theta)}{c_v}\right)u_x
                         +c\mathfrak h'(\theta)\theta_x\\
 &=-p_x+\frac{c_p}{v}\theta_x
   =\frac{R\theta}{v^2}v_x+\frac{c_v}{v}\theta_x
   =\frac{\theta}{v}\eta_x .
 \end{aligned}
\]
Differentiating \eqref{eq:normalform-definition} and using
$2k_0c\mathfrak h(\theta)=\frac{\theta}{v}$ give
\begin{align*}
 D_e w_e
 &=e^{-k_0\eta}\frac{\theta}{v}\eta_x
    -e k_0c\bigl[w_e+2e\mathfrak h(\theta_*)e^{k_0\eta}\bigr]\eta_x\\
 &=k_0c\bigl[2e^{-k_0\eta}\mathfrak h(\theta)
             -2\mathfrak h(\theta_*)e^{k_0\eta}-e w_e\bigr]\eta_x\\
 &=k_0c(w_+-w_--e w_e)\eta_x
   =-e k_0c\,w_{-e}\eta_x .
\end{align*}
In the moving coordinate \eqref{eq:moving-frame}, $D_e$ becomes
$\partial_t+(ec-b)\partial_\xi$, and the entropy operator becomes
$\partial_t-b\partial_\xi$. The full NSF equations for
$\mathsf w=(w_+,\eta,w_-)$ become
\begin{equation}\label{eq:triangular-Euler}
 \mathsf A(\mathsf w)=
 \begin{pmatrix}
 c(\mathsf w)&k_0c(\mathsf w)w_-&0\\
 0&0&0\\
 0&-k_0c(\mathsf w)w_+&-c(\mathsf w)
 \end{pmatrix},
 \qquad \mathsf w_t+(\mathsf A(\mathsf w)-bI)\mathsf w_\xi
     =\mathsf B(\mathsf w)\mathsf w_{\xi\xi}+\mathsf C(\mathsf w)[\mathsf w_\xi,\mathsf w_\xi].
\end{equation}
Since $q_*$ is fixed, the chain rule gives
\[
 \begin{aligned}
 \mathsf w_t&=DN_{q_*}(q)q_t,\qquad
 \mathsf w_\xi=DN_{q_*}(q)q_\xi,\\
 \mathsf w_{\xi\xi}
 &=DN_{q_*}(q)q_{\xi\xi}+D^2N_{q_*}(q)[q_\xi,q_\xi].
 \end{aligned}
\]
The diffusion matrix comes from the second derivatives in
\eqref{eq:st-spatial-operator}
\[
 \mathsf B(\mathsf w)
 =DN_{q_*}(q)\operatorname{diag}\left(0,\frac{\mu}{v},
                       \frac{\kappa}{c_vv}\right)(DN_{q_*}(q))^{-1}.
\]
The product rule for the fluxes gives
\[
 \begin{aligned}
 &DN_{q_*}(q)\begin{pmatrix}
     0\\[1mm]
     \mu\left(\frac{u_\xi}{v}\right)_\xi\\[2mm]
     \frac{\kappa}{c_v}\left(\frac{\theta_\xi}{v}\right)_\xi
                       +\frac{\mu u_\xi^2}{c_vv}
   \end{pmatrix}
       -\mathsf B(\mathsf w)\mathsf w_{\xi\xi}\\
 &\quad=DN_{q_*}(q)\begin{pmatrix}
      0\\[1mm]
      -\frac{\mu v_\xi u_\xi}{v^2}\\[2mm]
      -\frac{\kappa v_\xi\theta_\xi}{c_vv^2}
                       +\frac{\mu u_\xi^2}{c_vv}
    \end{pmatrix}
    -\mathsf B(\mathsf w)D^2N_{q_*}(q)[q_\xi,q_\xi]\\
 &\quad=\mathsf C(\mathsf w)[\mathsf w_\xi,\mathsf w_\xi].
 \end{aligned}
\]
Substituting $q_\xi=(DN_{q_*}(q))^{-1}\mathsf w_\xi$ gives a
quadratic expression in $\mathsf w_\xi$. Polarization gives the smooth
symmetric bilinear map $\mathsf C$.

The zero middle row of $\mathsf A$ and direct multiplication give
\begin{equation}\label{eq:affine-diagonalizer}
 \begin{aligned}
 T_0(\mathsf w)
 &=\begin{pmatrix}1&k_0w_-&0\\0&1&0\\0&k_0w_+&1\end{pmatrix},
 \qquad \det T_0=1,\qquad D_{\mathsf w}^2T_0=0,\\
 T_0\mathsf A&=\mathsf A=\operatorname{diag}(c,0,-c)T_0,\\
 T_0(\mathsf A-bI)
 &=\mathsf A-bT_0=\operatorname{diag}(c-b,-b,-c-b)T_0 .
 \end{aligned}
\end{equation}

By \eqref{eq:normalform-definition} and \eqref{eq:affine-diagonalizer},
the acoustic and contact rows of $T_0DN_{q_*}(q)$ are
\[
 \begin{aligned}
 L_e&=D_qw_e+k_0w_{-e}D_q\eta\\
 &=e^{-k_0\eta}(0,1,e\mathfrak h'(\theta))
       -k_0\bigl(w_e-w_{-e}+2e\mathfrak h(\theta_*)e^{k_0\eta}\bigr)D_q\eta\\
 &=e^{-k_0\eta}\bigl[(0,1,e\mathfrak h'(\theta))
                     -2e k_0\mathfrak h(\theta)D_q\eta\bigr]
  =e^{-k_0\eta}\left(-\frac{ec}{\gamma},1,\frac{eR}{vc}\right),\\
 L_0&=D_q\eta=\left(\frac{R}{v},0,\frac{c_v}{\theta}\right).
 \end{aligned}
\]
Add and subtract the acoustic rows and use the contact row to obtain
the inverse in the order $(+,0,-)$
\[
 (T_0DN_{q_*}(q))^{-1}
 =\begin{pmatrix}
 -{\displaystyle\frac{e^{k_0\eta}}{2c}}&{\displaystyle\frac{v}{c_p}}&{\displaystyle\frac{e^{k_0\eta}}{2c}}\\[2mm]
 {\displaystyle\frac{e^{k_0\eta}}{2}}&0&{\displaystyle\frac{e^{k_0\eta}}{2}}\\[2mm]
 {\displaystyle\frac{p e^{k_0\eta}}{2c_vc}}&{\displaystyle\frac{\theta}{c_p}}
                             &-{\displaystyle\frac{p e^{k_0\eta}}{2c_vc}}
 \end{pmatrix}.
\]
Substituting the inverse into the formula for $\mathsf B$ gives
\begin{align*}
 T_0\mathsf B T_0^{-1}
 &=(T_0DN_{q_*}(q))\operatorname{diag}\left(0,\frac{\mu}{v},
                      \frac{\kappa}{c_vv}\right)(T_0DN_{q_*}(q))^{-1}\\
 &=\begin{pmatrix}
 \beta_a&{\displaystyle\frac{p\kappa e^{-k_0\eta}}{c c_p c_vv}}
   &{\displaystyle\frac{\mu}{2v}}-{\displaystyle\frac{(\gamma-1)\kappa}{2c_pv}}\\[2mm]
 {\displaystyle\frac{R\kappa e^{k_0\eta}}{2c_vv^2c}}&\beta_0
   &-{\displaystyle\frac{R\kappa e^{k_0\eta}}{2c_vv^2c}}\\[2mm]
 {\displaystyle\frac{\mu}{2v}}-{\displaystyle\frac{(\gamma-1)\kappa}{2c_pv}}
   &-{\displaystyle\frac{p\kappa e^{-k_0\eta}}{c c_p c_vv}}&\beta_a
 \end{pmatrix},
\end{align*}
where the diagonal coefficients are
\begin{equation}\label{eq:effective-diffusion}
 \begin{aligned}
 \beta_a
 &=\frac{\mu}{2v}+\frac{R}{vc}\frac{\kappa}{c_vv}\frac{p}{2c_vc}
   =\frac{\mu}{2v}+\frac{(\gamma-1)\kappa}{2c_pv}>0,\\
 \beta_0
 &=\frac{c_v}{\theta}\frac{\kappa}{c_vv}\frac{\theta}{c_p}
   =\frac{\kappa}{c_pv}>0 .
 \end{aligned}
\end{equation}
Using the rows $L_+,L_0,L_-$, we have
\[
 \begin{aligned}
 &(T_0DN_{q_*}(q))^{\mathsf T}
    \operatorname{diag}\left(\frac{e^{2k_0\eta}}{2},
               \frac{\theta}{c_p},\frac{e^{2k_0\eta}}{2}\right)(T_0DN_{q_*}(q))\\
 &\quad=\frac{e^{2k_0\eta}}{2}(L_+^{\mathsf T}L_++L_-^{\mathsf T}L_-)
                +\frac{\theta}{c_p}L_0^{\mathsf T}L_0
       =\operatorname{diag}\left(\frac{R\theta}{v^2},1,
                                                \frac{c_v}{\theta}\right).
 \end{aligned}
\]
Multiplying by the diffusion matrix gives
\[
 \begin{aligned}
 &\operatorname{diag}\left(\frac{e^{2k_0\eta}}{2},
             \frac{\theta}{c_p},\frac{e^{2k_0\eta}}{2}\right)
     T_0\mathsf B T_0^{-1}\\
 &\quad=\bigl((T_0DN_{q_*}(q))^{-1}\bigr)^{\mathsf T}
      \operatorname{diag}\left(0,\frac{\mu}{v},
                        \frac{\kappa}{v\theta}\right)(T_0DN_{q_*}(q))^{-1},
 \end{aligned}
\]
which is symmetric positive semidefinite and has rank two.

\subsection{Corrected fields on the two half-lines}
Let $q_{\alpha,j}$ be the end state of $q_\alpha$ on side
$j\in\{-,+\}$, and put $\mathcal R_j(\alpha)=q_{\alpha,j}-q_j$.
On $\R_j$, with $\R_\pm=\{\xi:\pm\xi>0\}$, define
\begin{equation}\label{eq:corrected-end-field}
 Q^j=q_{\alpha,j}+r
     =q_j+Z-[\mathcal R(\alpha)-\mathcal R_j(\alpha)],
 \qquad U^j=N_{q_j}(Q^j).
\end{equation}
By \eqref{eq:st-residual} and \eqref{eq:st-parameter-tails},
\[
 |Q^j-q_j|\le |Z|+C_\delta|\alpha|e^{-2\eta_0|\xi|},
 \qquad \lim_{\xi\to j\infty}U^j(t,\xi)=0 .
\]
Differentiating \eqref{eq:corrected-end-field} gives
\[
 \begin{aligned}
 Q^j_\xi&=r_\xi,\qquad
 Q^j_t=F-D_\alpha[\mathcal R(\alpha)-\mathcal R_j(\alpha)]\alpha',\\
 |Q^j_t-F|&\le C_\delta e^{-2\eta_0|\xi|}|\alpha'|.
 \end{aligned}
\]
At speed $b=s_\alpha$, we have
$\mathscr N_b(q_\alpha)=\mathscr N_b(q_{\alpha,j})=0$.
Since $F=\mathscr N_b(q_\alpha+r)$, the source is
\[
 \begin{aligned}
 G^j&=\mathscr N_b(q_\alpha+r)-\mathscr N_b(q_\alpha)
      -\bigl[\mathscr N_b(q_{\alpha,j}+r)-\mathscr N_b(q_{\alpha,j})\bigr]\\
    &\quad-D_\alpha[\mathcal R(\alpha)-\mathcal R_j(\alpha)]\alpha'.
 \end{aligned}
\]
Applying \eqref{eq:triangular-Euler} to $Q^j$ gives
\begin{equation}\label{eq:corrected-end-PDE}
 \begin{aligned}
 Q^j_t&=\mathscr N_b(Q^j)+G^j,\\
 U^j_t+(\mathsf A(U^j)-bI)U^j_\xi
 &=\mathsf B(U^j)U^j_{\xi\xi}
  +\mathsf C(U^j)[U^j_\xi,U^j_\xi]+DN_{q_j}(Q^j)G^j.
 \end{aligned}
\end{equation}
The viscosity and heat terms in the last two rows of $G^j$ come from
differentiating the flux difference
\[
 \begin{aligned}
 &\frac{(\mu(U_{\alpha,\xi}+r_{2,\xi}),\frac{\kappa}{c_v}(\Theta_{\alpha,\xi}+r_{3,\xi}))}
            {V_\alpha+r_1}
     -\frac{(\mu U_{\alpha,\xi},\frac{\kappa}{c_v}\Theta_{\alpha,\xi})}{V_\alpha}
     -\frac{(\mu r_{2,\xi},\frac{\kappa}{c_v}r_{3,\xi})}{v_{\alpha,j}+r_1}\\
 &\quad=-\frac{r_1(\mu U_{\alpha,\xi},\frac{\kappa}{c_v}\Theta_{\alpha,\xi})}
                  {V_\alpha(V_\alpha+r_1)}
     +\frac{(v_{\alpha,j}-V_\alpha)(\mu r_{2,\xi},\frac{\kappa}{c_v}r_{3,\xi})}
                  {(V_\alpha+r_1)(v_{\alpha,j}+r_1)} .
 \end{aligned}
\]
Estimating the pressure and heating terms by the same coefficient
differences and using \eqref{eq:profilebounds} and
\eqref{eq:st-parameter-tails}, we get
\[
 \begin{aligned}
 |G^j|
 &\le C_\delta |q_\alpha-q_{\alpha,j}|
       \bigl(|r_\xi|+|\partial_\xi^2 r_{\rm par}|+|r_\xi|^2\bigr)
    +C_\delta\bigl(|q_{\alpha,\xi}|+|\partial_\xi^2 q_\alpha|\bigr)|r|\\
 &\quad+C_\delta|q_{\alpha,\xi}|\,|r_\xi|
       +\bigl|D_\alpha[\mathcal R(\alpha)-\mathcal R_j(\alpha)]\bigr|
                |\alpha'|\\
 &\le C_\delta e^{-2\eta_0|\xi|}
     \bigl(|r|+|r_\xi|+|\partial_\xi^2 r_{\rm par}|
                         +|r_\xi|^2+|\alpha'|\bigr).
 \end{aligned}
\]
Take norms on $\R_j$ for $G^j,U^j$ and on $\R$ for $r,F$. The bound
after \eqref{eq:corrected-end-field} keeps $Q^j$ in a compact set
where $DN_{q_j}$ and $D^2N_{q_j}$ are bounded. From
\eqref{eq:st-moment-bounds} and $\alpha'=\mathfrak mF$, we have
$|\alpha'|\le C_\delta\min\{\|F\|_2,m_\omega(F)\}$.
Use $e^{-2\eta_0|\xi|}\le C_\delta\omega(\xi)$ and
$\|r_\xi\|_\infty\le C_\delta\eps$ in the bound for $G^j$.
Then \eqref{eq:st-spatial-estimates} gives
\begin{equation}\label{eq:end-source-bounds}
 \begin{aligned}
 &\left\|\frac{G^j}{\sqrt{\omega}}\right\|_2
     +\left\|\frac{DN_{q_j}(Q^j)G^j}{\sqrt{\omega}}\right\|_2\\
 &\quad\le C_\delta\bigl(m_\omega(r)+m_\omega(r_\xi)
                   +m_\omega(\partial_\xi^2r_{\rm par})+|\alpha'|\bigr)\\
 &\quad\le C_\delta\bigl(\min\{X_1(r),X_2(r)\}+|\alpha'|\bigr)
       \le C_\delta\min\{\|F\|_2,\Gamma\}.
 \end{aligned}
\end{equation}

Differentiating $U^j=N_{q_j}(Q^j)$ gives
\[
 \begin{aligned}
 U^j_\xi&=DN_{q_j}(Q^j)r_\xi,\\
 U^j_{\xi\xi}
 &=DN_{q_j}(Q^j)r_{\xi\xi}+D^2N_{q_j}(Q^j)[r_\xi,r_\xi],\\
 U^j_t
 &=DN_{q_j}(Q^j)\bigl(F-D_\alpha[\mathcal R(\alpha)
                                  -\mathcal R_j(\alpha)]\alpha'\bigr).
 \end{aligned}
\]
Apply \eqref{eq:st-GN-bounded} to the quadratic term and use
\eqref{eq:st-spatial-estimates}, \eqref{eq:st-endpoint-recovery}
and \eqref{eq:st-h} to get
\begin{equation}\label{eq:end-field-derivatives}
 \begin{aligned}
 \|U^j_\xi\|_2
 &\le C_\delta\|r_\xi\|_2\le C_\delta\|F\|_2,\\
 \|U^j_{\xi\xi}\|_2
 &\le C_\delta\bigl(\|r_{\xi\xi}\|_2+\|r_\xi\|_4^2\bigr)\\
 &\le C_\delta(1+\|r\|_\infty)\|r_{\xi\xi}\|_2
       \le C_\delta\min\{\Gamma,\sqrt{\mathcal T}\},\\
 \|U^j_\xi\|_\infty+\|U^j_t\|_\infty
 &\le C_\delta\bigl(\|r_\xi\|_\infty+\|F\|_\infty+|\alpha'|\bigr)
       \le C_\delta h.
 \end{aligned}
\end{equation}
The bounds involving $\|F\|_2$ and $\sqrt{\mathcal T}$ hold at each
time. Those involving $\Gamma$ and $h$ hold almost everywhere.

\subsection{A positive parabolic system for the frozen modes}\label{sec:lift}

On $\R_j$, use $x=j\xi$ and freeze the coefficients in
\eqref{eq:corrected-end-PDE} at $\vartheta$ in a compact neighborhood
of the origin in the $U^j$ coordinates. The bounds are uniform on
these neighborhoods. The neighborhood sizes, switching rate and
frequency cutoff may depend on the fixed shock and on $\mu,\kappa$.

The first-order part is $j\{b(\tau)I-\mathsf A(\vartheta)\}\partial_x$.
For $z=ik$, the full symbol is $jb(\tau)zI+\mathsf P(\vartheta,z)$, where
\begin{equation}\label{eq:lift-physical-symbol}
 \begin{gathered}
 \mathsf P(\vartheta,z)=z\mathsf P_1(\vartheta)
                         +z^2\mathsf P_2(\vartheta),\qquad
 \mathsf P_1(\vartheta)=-j\mathsf A(\vartheta),\qquad
 \mathsf P_2(\vartheta)=\mathsf B(\vartheta).
 \end{gathered}
\end{equation}
By \eqref{eq:affine-diagonalizer}, the eigenvalues of $\mathsf P_1$ are
$-jc(\vartheta),0,jc(\vartheta)$ in the order $(+,0,-)$. They are
uniformly separated on the compact neighborhoods.

Let $a$ be the eigenvalue of $\mathsf P_1$ for branch $e\in\{+,0,-\}$.
The normalized eigenvalue equations are
\[
 t(\vartheta,z)(\mathsf P_1+z\mathsf P_2)
       =\frac{\lambda(\vartheta,z)}{z}\,t(\vartheta,z),\qquad
 (t(\vartheta,z)T_0^{-1})_e=1,
\]
At $z=0$, the eigenvalue variable $\frac{\lambda(\vartheta,z)}{z}$
equals $a$ and $t(\vartheta,0)$ is the corresponding row of $T_0$.
After fixing $(tT_0^{-1})_e=1$, the linearization in the two other
components and the eigenvalue variable is diagonal. Its entries
are the two nonzero eigenvalue gaps and $-1$. The implicit function
theorem and compactness give a common disk on which $t$ and
$\frac{\lambda}{z}$ are analytic in $z$ and smooth in $\vartheta$.
For real $z$ the equations are real, so uniqueness gives real values.

Multiplying the normalized equation by $T_0^{-1}$ and taking
its $e$-th component gives
\[
 \lambda(\vartheta,z)=a z
          +z^2\bigl[t(\vartheta,z)\mathsf B(\vartheta)T_0^{-1}\bigr]_e .
\]
Taylor expansion and \eqref{eq:effective-diffusion} yield
\begin{equation}\label{eq:lift-root-expansion}
 \begin{aligned}
 \lambda(\vartheta,z)&=a z+2Dz^2+c_3z^3+c_4z^4+O(z^5),\\
 2D&=\bigl(T_0\mathsf B T_0^{-1}\bigr)_{ee}
   =\begin{cases}
       \beta_a,&e=+,-,\\
       \beta_0,&e=0
     \end{cases}
   \ge c_\delta>0 .
 \end{aligned}
\end{equation}
The lower bound and the remainder, including its $\vartheta$
derivatives through order two, are uniform on the compact
neighborhoods. The coefficients are real, so on a smaller common
frequency neighborhood
\[
 \operatorname{Re}\lambda(\vartheta,ik)
   =-2Dk^2+O(k^4)\le -c_\delta k^2 .
\]

\begin{lemma}[Positive realization of four coefficients]\label{lemma:lift}
For each branch, use $D$ from \eqref{eq:lift-root-expansion}.
There are smooth functions $x_-<0<x_+$, smooth positive probabilities
$\pi_-,\pi_0,\pi_+$ and a constant $\rho>0$ for which the matrix
\begin{equation}\label{eq:lift-auxiliary-symbol}
 \begin{aligned}
 \mathsf M(\vartheta,z)
 &=Dz^2I_3+z\operatorname{diag}(a+x_-,a,a+x_+)+\mathsf Q,\\
 \mathsf Q&=\rho(\mathbf1\pi^{\mathsf T}-I_3)
 \end{aligned}
\end{equation}
has a simple slow root $\mu(\vartheta,z)$ with
\begin{equation}\label{eq:lift-order-five}
 \lambda-\mu=O(z^5).
\end{equation}
The bound in \eqref{eq:lift-order-five} holds also after up to two
parameter derivatives. The corresponding system in physical space is
cooperative and uniformly parabolic. For a sufficiently small choice of
$\rho$, each auxiliary drift has the sign of its physical branch after
adding the common moving-frame velocity.
\end{lemma}

\begin{proof}
The positive lower bound of $D$ and the bounds on $c_3,c_4$ let us
choose $0<\rho\le1$ for both ends and all three branches with
\begin{equation}\label{eq:lift-rate-choice}
 \rho\sup_{j,e,\vartheta}|Dc_4-c_3^2|
       \le\frac{1}{2}\inf_{j,e,\vartheta} D^3 .
\end{equation}
Define the moments
\[
 m_2=D\rho,\qquad m_3=c_3\rho^2,\qquad
 m_4=c_4\rho^3+2D^2\rho^2
\]
and the nodes
\begin{equation}\label{eq:lift-three-nodes}
 \alpha_*=\frac{m_3}{m_2}=\frac{c_3}{D}\rho,\qquad
 \beta_*=\frac{m_4}{m_2}-\alpha_*^2,\qquad
 x_\pm=\frac{\alpha_*\pm\sqrt{\alpha_*^2+4\beta_*}}{2},
 \qquad x_0=0.
\end{equation}
By \eqref{eq:lift-rate-choice},
\begin{equation}\label{eq:lift-node-positivity}
 \begin{aligned}
 \beta_*&=2m_2+\frac{Dc_4-c_3^2}{D^2}\rho^2,\qquad
           \frac{3}{2}m_2\le\beta_*\le\frac{5}{2}m_2,\\
 |x_\pm|&\le|\alpha_*|+\sqrt{\beta_*}\le C_\delta\sqrt{\rho}.
 \end{aligned}
\end{equation}
In particular, $x_-<0<x_+$ and $x_+x_-=-\beta_*$. Take
\begin{equation}\label{eq:lift-probabilities}
 \pi_+=\frac{m_2}{x_+(x_+-x_-)},\qquad
 \pi_-=-\frac{m_2}{x_-(x_+-x_-)},\qquad
 \pi_0=1-\frac{m_2}{\beta_*}.
\end{equation}
These weights satisfy
\[
 \frac{1}{3}\le\pi_0\le\frac{3}{5},\qquad
 \pi_\pm=\frac{m_2}{|x_\pm|(x_+-x_-)}
       \ge\frac{m_2}{C_\delta\rho}\ge c_\delta>0.
\]
Using $x_\pm^2=\alpha_*x_\pm+\beta_*$ and summing over
$\ell\in\{-,0,+\}$ gives
\begin{equation}\label{eq:lift-moments}
 \begin{aligned}
 \sum_\ell\pi_\ell
 &=\pi_0+\frac{m_2}{x_+-x_-}
     \left(\frac{1}{x_+}-\frac{1}{x_-}\right)
   =\pi_0+\frac{m_2}{\beta_*}=1,\\
 \sum_\ell\pi_\ell x_\ell
 &=\frac{m_2}{x_+-x_-}-\frac{m_2}{x_+-x_-}=0,\\
 \sum_\ell\pi_\ell x_\ell^2
 &=\frac{m_2(x_+-x_-)}{x_+-x_-}=m_2,\\
 \sum_\ell\pi_\ell x_\ell^3
 &=\sum_{\ell=\pm}\pi_\ell x_\ell(\alpha_*x_\ell+\beta_*)
   =\alpha_*m_2=m_3,\\
 \sum_\ell\pi_\ell x_\ell^4
 &=\sum_{\ell=\pm}\pi_\ell x_\ell^2(\alpha_*x_\ell+\beta_*)
   =\alpha_*m_3+\beta_*m_2=m_4 .
 \end{aligned}
\end{equation}
The off-diagonal entries of $\mathsf Q$ are $\rho\pi_\ell>0$ and
its row sums are zero. This gives a cooperative system. The lower
bound on $D$ gives uniform parabolicity.

By \eqref{eq:lift-moments},
\[
 \mathsf Q\mathbf1=0,\qquad
 \mathsf Q|_{\ker\pi^{\mathsf T}}=-\rho I,\qquad
 \mathsf Q^2=-\rho\mathsf Q .
\]
Thus $0$ is simple and $-\rho$ is semisimple with multiplicity two.
Writing $\mu=az+Dz^2+y$, the matrix determinant lemma gives
\[
 \begin{aligned}
 \mu I_3-\mathsf M
 &=\operatorname{diag}_{\ell}(\rho+y-x_\ell z)
                      -\rho\mathbf1\pi^{\mathsf T},\\
 \det(\mu I_3-\mathsf M)
 &=\prod_\ell(\rho+y-x_\ell z)
      \left(1-\rho\sum_\ell\frac{\pi_\ell}{\rho+y-x_\ell z}\right).
 \end{aligned}
\]
When $|y|$ and $|z|$ are small, the diagonal factors have modulus at least
$\frac{\rho}{2}$ on the parameter neighborhoods. The slow root satisfies
\begin{equation}\label{eq:lift-reset-equation}
 1-\rho\sum_\ell\frac{\pi_\ell}{\rho+y-x_\ell z}=0.
\end{equation}
At $(y,z)=(0,0)$ the left-hand side is zero. Its $y$ derivative is
$\frac{1}{\rho}$ and its $z$ derivative is
$-\frac{1}{\rho}\sum_\ell\pi_\ell x_\ell=0$.
The implicit function theorem gives an analytic root $y$ on a common
small disk, unique near zero and smooth in $\vartheta$, with
$y(\vartheta,0)=\partial_z y(\vartheta,0)=0$.

Since $y=O(z^2)$, expansion of \eqref{eq:lift-reset-equation} and
\eqref{eq:lift-moments} give
\[
 \begin{aligned}
 0&=\sum_\ell\pi_\ell
  \left(\frac{y-x_\ell z}{\rho}-\frac{(y-x_\ell z)^2}{\rho^2}
            +\frac{(y-x_\ell z)^3}{\rho^3}
            -\frac{(y-x_\ell z)^4}{\rho^4}\right)+O(z^5)\\
  &=\frac{y}{\rho}-\frac{y^2+m_2z^2}{\rho^2}
       +\frac{3m_2yz^2-m_3z^3}{\rho^3}
       -\frac{m_4z^4}{\rho^4}+O(z^5).
 \end{aligned}
\]
In particular, $y=\frac{m_2}{\rho}z^2+O(z^3)$. Substituting this
expression in the terms $y^2$ and $yz^2$ yields
\[
 \begin{aligned}
 y&=\frac{m_2}{\rho}z^2+\frac{m_3}{\rho^2}z^3
            +\frac{m_4-2m_2^2}{\rho^3}z^4+O(z^5)
   =Dz^2+c_3z^3+c_4z^4+O(z^5),\\
 \mu&=az+Dz^2+y=az+2Dz^2+c_3z^3+c_4z^4+O(z^5).
 \end{aligned}
\]
On a smaller common disk the Taylor bounds are uniform for $y$
and its $\vartheta$ derivatives through order two. Comparing with
\eqref{eq:lift-root-expansion} proves \eqref{eq:lift-order-five},
including these derivatives.

In the outward coordinate $x=j\xi$, the physical coefficient
of $\partial_x$ is $jb(\tau)+a=j(b(\tau)-e c(\vartheta))$, with
$e=+1,0,-1$. The three reference coefficients at $\vartheta=0$ and
$b(\tau)=\sigma$ are
\[
 \bigl(j(\sigma-c_j),j\sigma,j(\sigma+c_j)\bigr)
 =
 \begin{cases}
  (\sigma-c_+,\sigma,\sigma+c_+),&j=+,\\
  (c_--\sigma,-\sigma,-\sigma-c_-),&j=- .
 \end{cases}
\]
By \eqref{eq:lax}, the moduli are at least
$\min\{\sigma-c_+,c_--\sigma\}>0$.
Equations \eqref{eq:moving-frame}, \eqref{eq:feedbackbounds} and
\eqref{eq:global-H2} give $|b(\tau)-\sigma|\le C_\delta\eps$ for all
times. By \eqref{eq:lift-node-positivity}, small $\rho$, small
parameter neighborhoods and small initial norm give
\[
 \begin{aligned}
 |jb(\tau)+a+x_\ell-j(\sigma-e c_j)|
 &\le |b(\tau)-\sigma|+|c(\vartheta)-c_j|+|x_\ell|\\
 &\le C_\delta(\eps+\sqrt{\rho})+|c(\vartheta)-c_j|\\
 &\le\frac{1}{2}\min\{\sigma-c_+,c_--\sigma\}.
 \end{aligned}
\]
The physical coefficient and the three auxiliary coefficients thus
have the sign of $j(\sigma-e c_j)$ and stay bounded away from zero.
A positive coefficient of $\partial_x$ transports toward $x=0$.
The three right branches and the left $+$ branch are incoming.
\end{proof}

For each branch define the column
\begin{equation}\label{eq:lift-exact-column}
 j_\ell(\vartheta,z)
 =\frac{\rho}{\lambda(\vartheta,z)-Dz^2+\rho-(a+x_\ell)z},
 \qquad \ell\in\{-,0,+\}.
\end{equation}
The denominators with $\lambda$ or $\mu$ equal $\rho$ at $z=0$
and have modulus at least $\frac{\rho}{2}$ on a smaller common disk.
Thus $j$ is analytic in $z$ and smooth in $\vartheta$, with
$j(\vartheta,0)=\mathbf1$. Substituting $y=\mu-az-Dz^2$ in
\eqref{eq:lift-reset-equation} and using \eqref{eq:lift-auxiliary-symbol}
give
\begin{equation}\label{eq:lift-column-residual}
 \begin{aligned}
 &(\lambda I_3-\mathsf M)j
    =\rho(1-\pi^{\mathsf T}j)\mathbf1\\
 &\quad=\rho^2(\lambda-\mu)\sum_\ell
 \frac{\pi_\ell\mathbf1}{
       [\lambda-Dz^2+\rho-(a+x_\ell)z]
       [\mu-Dz^2+\rho-(a+x_\ell)z]}
    =O(z^5).
 \end{aligned}
\end{equation}
The factor multiplying $\lambda-\mu$ has bounded $\vartheta$
derivatives through order two on this disk. By \eqref{eq:lift-order-five},
the remainder in \eqref{eq:lift-column-residual} has the same uniform
order after these derivatives.

Choose a fixed real even cutoff $\chi\in C_c^\infty((-k_*,k_*))$
equal to one near zero and independent of $\vartheta$. Take $k_*$
small enough for \eqref{eq:lift-exact-column}. At $z=ik$, the three
branches give the $9\times3$ multiplier
\begin{equation}\label{eq:lift-multiplier}
 \mathsf L(\vartheta,k)
 =\chi(k)\begin{pmatrix}j_+t_+\\j_0t_0\\j_-t_-\end{pmatrix},
 \qquad
 \mathsf M=\operatorname{diag}(\mathsf M_+,\mathsf M_0,\mathsf M_-).
\end{equation}
Here the branch columns and rows are evaluated at $(\vartheta,ik)$.
Using $t_e\mathsf P=\lambda_e t_e$ and
\eqref{eq:lift-column-residual}, we obtain
\begin{equation}\label{eq:lift-intertwining}
 \mathsf L\mathsf P-\mathsf M\mathsf L
 =\chi(k)\begin{pmatrix}
   \bigl[(\lambda_+I_3-\mathsf M_+)j_+\bigr]t_+\\
   \bigl[(\lambda_0I_3-\mathsf M_0)j_0\bigr]t_0\\
   \bigl[(\lambda_-I_3-\mathsf M_-)j_-\bigr]t_-
 \end{pmatrix}
 =O(k^5).
\end{equation}
By \eqref{eq:lift-column-residual}, the quotient by $(ik)^5$ extends
smoothly across $k=0$, including its $\vartheta$ derivatives through
order two. Since $\chi$ is real and even and $j_e,t_e$ are real for
real $z$, we have $\mathsf L(\vartheta,-k)=\overline{\mathsf L(\vartheta,k)}$.

At $k=0$, $j_e=\mathbf1$ and $t_e$ is the corresponding row of
$T_0$. By \eqref{eq:affine-diagonalizer}, with the input coordinates
ordered as $(+,0,-)$,
\[
 \mathsf L(\vartheta,0)
 =\begin{pmatrix}
    \mathbf1&k_0\vartheta_-\mathbf1&\mathbf0\\
    \mathbf0&\mathbf1&\mathbf0\\
    \mathbf0&k_0\vartheta_+\mathbf1&\mathbf1
  \end{pmatrix},
 \qquad
 \partial_{\vartheta_0}\mathsf L(\vartheta,0)=0,\quad
 D_\vartheta^2\mathsf L(\vartheta,0)=0 .
\]

\begin{lemma}[Difference kernel]\label{lemma:lift-kernel}
Define
\begin{equation}\label{eq:lift-difference-definition}
 \widehat{\mathsf K}(\vartheta,k,\tau)
 =\mathsf L(\vartheta,k)e^{\tau\mathsf P(\vartheta,ik)}
       -e^{\tau\mathsf M(\vartheta,ik)}\mathsf L(\vartheta,k).
\end{equation}
For every integer $n\ge-1$, every parameter multi-index $\beta$ with
$|\beta|\le2$, and $1\le p\le\infty$,
\begin{equation}\label{eq:lift-kernel-bound}
 \|\partial_x^n\partial_\vartheta^\beta
             \mathsf K(\vartheta,\cdot,\tau)\|_{L^p}
 \le C_{n,\beta,p}
 \begin{cases}
  \tau,&0\le\tau\le1,\\
  \tau^{-\frac{3+n-|\beta|}{2}-\frac{1}{2}(1-\frac{1}{p})},&\tau\ge1.
 \end{cases}
\end{equation}
For $n=-1$, divide the Fourier symbol by $ik$, using its vanishing at
$k=0$. The same bounds hold after the common moving-frame translation.
Keep that translation fixed under parameter differentiation.
\end{lemma}

\begin{proof}
At $z=0$, $\mathsf M=\mathsf Q$. By \eqref{eq:lift-moments},
$\mathbf1\pi^{\mathsf T}$ is the projection onto the zero eigenspace
of $\mathsf Q$. On the circle $|\zeta+\rho|=\frac{\rho}{2}$,
\[
 (\zeta I_3-\mathsf Q)^{-1}
 =\frac{\mathbf1\pi^{\mathsf T}}{\zeta}
   +\frac{I_3-\mathbf1\pi^{\mathsf T}}{\zeta+\rho},
 \qquad \|(\zeta I_3-\mathsf Q)^{-1}\|\le\frac{C}{\rho}.
\]
For small $|z|$, we have
$\|(\zeta I_3-\mathsf Q)^{-1}(\mathsf M-\mathsf Q)\|\le\frac{1}{2}$
uniformly in $\vartheta$ on this circle. The Neumann series gives
\[
 \|(\zeta I_3-\mathsf M)^{-1}\|
 \le 2\|(\zeta I_3-\mathsf Q)^{-1}\|\le\frac{C}{\rho}.
\]
Choose the cutoff in this disk with $|\lambda|,|\mu|\le\frac{\rho}{4}$.
As at $z=0$, the circle encloses two eigenvalues counted with
multiplicity. The slow root $\mu$ lies outside. Let $\Pi_{\rm f}$
be the Riesz projection onto the enclosed spectrum. Then
\[
 e^{\tau\mathsf M(\vartheta,z)}\Pi_{\rm f}(\vartheta,z)
 =\frac{1}{2\pi i}\int_{|\zeta+\rho|=\frac{\rho}{2}}
    e^{\tau\zeta}(\zeta I_3-\mathsf M(\vartheta,z))^{-1}\dd\zeta .
\]
At $\tau=0$, the formula gives $\Pi_{\rm f}$, analytic in $z$
and smooth in $\vartheta$. Differentiating the resolvent gives
products of resolvents and derivatives of $\mathsf M$, all bounded
on the circle. For $m\ge0$ and $|\beta|\le2$, this gives
\[
 \begin{aligned}
 &\bigl\|\partial_z^m\partial_\vartheta^\beta
           (e^{\tau\mathsf M}\Pi_{\rm f})\bigr\|\\
 &\quad\le\frac{1}{2\pi}\int_{|\zeta+\rho|=\frac{\rho}{2}}
   e^{\tau\operatorname{Re}\zeta}
   \bigl\|\partial_z^m\partial_\vartheta^\beta
                (\zeta I_3-\mathsf M)^{-1}\bigr\|\,|\dd\zeta|
   \le C_{m,\beta}e^{-\frac{\rho\tau}{2}}.
 \end{aligned}
\]
On the circle, $|\lambda-\zeta|\ge\frac{\rho}{4}$. Since $\lambda$
lies outside, the left side of
\[
 \frac{(\zeta I_3-\mathsf M)^{-1}(\lambda I_3-\mathsf M)}
      {\lambda-\zeta}
 =(\zeta I_3-\mathsf M)^{-1}+\frac{I_3}{\lambda-\zeta}
\]
integrates to $2\pi i\Pi_{\rm f}$. Apply the identity to $j$ and use
\eqref{eq:lift-column-residual} to obtain
\begin{equation}\label{eq:lift-fast-column}
 \Pi_{\rm f}j
 =\frac{1}{2\pi i}\int_{|\zeta+\rho|=\frac{\rho}{2}}
       \frac{(\zeta I_3-\mathsf M)^{-1}(\lambda I_3-\mathsf M)j}
            {\lambda-\zeta}\dd\zeta
 =O(z^5).
\end{equation}
The quotient $\frac{\Pi_{\rm f}j}{z^5}$ is analytic at zero and
bounded uniformly after up to two $\vartheta$ derivatives.

On each branch, $t e^{\tau\mathsf P}=e^{\tau\lambda}t$, and the
slow projection satisfies
$e^{\tau\mathsf M}(I_3-\Pi_{\rm f})=e^{\tau\mu}(I_3-\Pi_{\rm f})$.
The column multiplying $\chi(k)t$ in \eqref{eq:lift-difference-definition}
is therefore
\begin{equation}\label{eq:lift-kernel-splitting}
 \begin{aligned}
 (e^{\tau\lambda}I_3-e^{\tau\mathsf M})j
 &=e^{\tau\lambda}j-e^{\tau\mu}(I_3-\Pi_{\rm f})j
                           -e^{\tau\mathsf M}\Pi_{\rm f}j\\
 &=(e^{\tau\lambda}-e^{\tau\mu})j
       +(e^{\tau\mu}I_3-e^{\tau\mathsf M})\Pi_{\rm f}j.
 \end{aligned}
\end{equation}
Both slow roots have positive quadratic coefficients. On a smaller
common support this gives
\[
 \operatorname{Re}\lambda(\vartheta,ik),\quad
 \operatorname{Re}\mu(\vartheta,ik)\le-\theta k^2.
\]
For the auxiliary matrix, we can also use
\[
 \operatorname{Re}\langle f,\mathsf Qf\rangle_\pi
 =\rho|\pi^{\mathsf T}f|^2-\rho\|f\|_\pi^2\le0.
\]
The real diagonal drifts are skew-adjoint when multiplied by $ik$.

For $0\le\tau\le1$, compute \eqref{eq:lift-difference-definition}
with the full matrices by differentiating
$e^{(\tau-r)\mathsf M}\mathsf L e^{r\mathsf P}$ in $r$ and integrating
from $0$ to $\tau$
\[
 \begin{aligned}
 \widehat{\mathsf K}
 &=\int_0^\tau e^{(\tau-r)\mathsf M}
               (\mathsf L\mathsf P-\mathsf M\mathsf L)
                     e^{r\mathsf P}\dd r,\\
 (ik)^n\partial_\vartheta^\beta\widehat{\mathsf K}
 &=(ik)^{n+5}\int_0^\tau \partial_\vartheta^\beta
   \left[e^{(\tau-r)\mathsf M}
      \frac{\mathsf L\mathsf P-\mathsf M\mathsf L}{(ik)^5}
       e^{r\mathsf P}\right]\dd r .
 \end{aligned}
\]
Equation \eqref{eq:lift-intertwining} bounds the integrand and its
first $k$ derivative uniformly on the cutoff support for
$0\le r\le\tau\le1$ and $|\beta|\le2$. When $n=-1$, the factor is
$(ik)^4$ and extends smoothly at zero. Fourier inversion gives the
$L^\infty$ bound. For the $L^1$ bound, use Cauchy--Schwarz with
weight $1+x^2$ and then Plancherel
\[
 \begin{aligned}
 &\|\partial_x^n\partial_\vartheta^\beta\mathsf K\|_{L^\infty_x}
       +\|\partial_x^n\partial_\vartheta^\beta\mathsf K\|_{L^1_x}\\
 &\quad\le C\Bigl(
    \|(ik)^n\partial_\vartheta^\beta\widehat{\mathsf K}\|_{L^1_k}
      +\|(ik)^n\partial_\vartheta^\beta\widehat{\mathsf K}\|_{H^1_k}
             \Bigr)
       \le C_{n,\beta}\tau .
 \end{aligned}
\]
Interpolation gives \eqref{eq:lift-kernel-bound} on $0\le\tau\le1$.

For $\tau\ge1$, set $k=\frac{\zeta}{\sqrt{\tau}}$ in each block
$\widehat{\mathsf K}_e$. For the first term in
\eqref{eq:lift-kernel-splitting}, use
\[
 e^{\tau\lambda}-e^{\tau\mu}
 =\tau(\lambda-\mu)\int_0^1
                 e^{\tau[\ell\lambda+(1-\ell)\mu]}\dd\ell .
\]
By \eqref{eq:lift-root-expansion} and \eqref{eq:lift-order-five},
for $m=0,1$, $|\beta'|\le2$ and $0\le\ell\le1$,
\[
 \begin{aligned}
 \bigl|\partial_\zeta^m\partial_\vartheta^{\beta'}
      \{\tau[\ell\lambda+(1-\ell)\mu-ia k]\}\bigr|
 &\le C_{\beta'}|\zeta|^{2-m},\\
 \operatorname{Re}\{\tau[\ell\lambda+(1-\ell)\mu-ia k]\}
 &\le-\theta\zeta^2 .
 \end{aligned}
\]
The centered roots are $O(k^2)$ and their first $k$ derivatives
are $O(k)$, also after the stated parameter derivatives. The chain
rule bounds the centered exponential and its derivatives above by
$C(1+|\zeta|)^5 e^{-\theta\zeta^2}$. For the phase we have
\[
 \sum_{m=0}^1
 \bigl|\partial_\zeta^m\bigl(
 e^{-ia\sqrt{\tau}\zeta}\partial_\vartheta^{\beta'}
                       e^{ia\sqrt{\tau}\zeta}\bigr)\bigr|
 \le C_{\beta'}\tau^{\frac{|\beta'|}{2}}(1+|\zeta|)^{|\beta'|}.
\]
In each block, take the parameter derivatives before removing
the phase $e^{ia\sqrt{\tau}\zeta}$.

The factors $\chi$, $j$, $t$, $\frac{\lambda-\mu}{(ik)^5}$ and
$\frac{\Pi_{\rm f}j}{(ik)^5}$ have uniform bounds on the support
through one $k$ derivative and two parameter derivatives.
The powers of $k$ in the two slow terms give
\[
 \tau(ik)^{n+5}=\tau^{-\frac{3+n}{2}}(i\zeta)^{n+5},
 \qquad
 (ik)^{n+5}=\tau^{-\frac{5+n}{2}}(i\zeta)^{n+5}.
\]
The second factor is for $e^{\tau\mu}\Pi_{\rm f}j$, which uses
$\ell=0$ in the centered bound. Differentiating a cutoff or a quotient
in $\zeta$ gives a factor $\frac{1}{\sqrt{\tau}}$. Apply the contour
bound to the fast term $(e^{\tau\mathsf M}\Pi_{\rm f})(\Pi_{\rm f}j)$.
On $|\zeta|\le k_*\sqrt{\tau}$ we have
\[
 e^{-\frac{\rho\tau}{2}}
 \le e^{-\frac{\rho\tau}{4}}e^{-\frac{\rho}{4k_*^2}\zeta^2}.
\]
The time exponential absorbs the $\sqrt{\tau}$ from differentiating
$e^{-ia\sqrt{\tau}\zeta}$ once in $\zeta$. Inserting these bounds
in \eqref{eq:lift-kernel-splitting} gives, for $m=0,1$,
\[
 \begin{aligned}
 &\left|\partial_\zeta^m\left[
     \tau^{\frac{3+n-|\beta|}{2}}e^{-ia\sqrt{\tau}\zeta}
     (ik)^n\partial_\vartheta^\beta
       \widehat{\mathsf K}_e\left(\vartheta,
                          \frac{\zeta}{\sqrt{\tau}},\tau\right)
                    \right]\right|\\
 &\quad\le C_{n,\beta}(1+|\zeta|)^{n+12}e^{-c_\delta\zeta^2}.
 \end{aligned}
\]
For $n=-1$, the polynomial is $(i\zeta)^4$. The bracketed expression
is bounded in $L^1_\zeta$ and $H^1_\zeta$, so its inverse transform
has bounded $L^p$ norm by the same endpoint estimates as above.
Inverting in $k$ introduces the factor $\frac{1}{\sqrt{\tau}}$
and the rescaled argument $\frac{x+a\tau}{\sqrt{\tau}}$. Hence
\[
 \begin{aligned}
 \|\partial_x^n\partial_\vartheta^\beta\mathsf K_e\|_{L^p_x}
 &\le C_{n,\beta,p}\tau^{-\frac{3+n-|\beta|}{2}}
                    \tau^{-\frac{1}{2}+\frac{1}{2p}}\\
 &= C_{n,\beta,p}
     \tau^{-\frac{3+n-|\beta|}{2}-\frac{1}{2}(1-\frac{1}{p})}.
 \end{aligned}
\]
The three block bounds give \eqref{eq:lift-kernel-bound}.
The common drift on $[s,t]$ gives the Fourier factor
$\exp(ikj\int_s^t b(r)\dd r)$, held fixed under $\vartheta$ derivatives.
It is a spatial translation and preserves the $L^p$ norms.
\end{proof}

The table lists the decay exponents in \eqref{eq:lift-kernel-bound}
for $\tau\ge1$.
\begin{equation}\label{eq:lift-kernel-table}
\begin{array}{c|rrrrrr}
 &\mathsf K&\partial_x^{-1}\mathsf K&\mathsf K_\vartheta
 &\partial_x^{-1}\mathsf K_\vartheta&\mathsf K_{\vartheta\vartheta}
 &\partial_x^{-1}\mathsf K_{\vartheta\vartheta}\\ \hline
 L^1_x&\frac{3}{2}&1&1&\frac{1}{2}&\frac{1}{2}&0\\
 L^2_x&\frac{7}{4}&\frac{5}{4}&\frac{5}{4}&\frac{3}{4}&\frac{3}{4}&\frac{1}{4}\\
 L^\infty_x&2&\frac{3}{2}&\frac{3}{2}&1&1&\frac{1}{2}
\end{array}
\end{equation}
All these kernels are $O(\tau)$ for $0\le\tau\le1$.

\section{Comparison on finite time intervals}\label{sec:window-comparison}

\subsection{Transverse estimates for the actual comparison kernels}
For $I=[s,t]$, retain the notation in \eqref{eq:time-window-notation}
and \eqref{eq:st-h}. Set
\begin{equation}\label{eq:window-notation}
 \begin{aligned}
 M_I&=\operatorname*{ess\,sup}_{\tau\in I}\sqrt{\mathcal T(\tau)},
 &\Theta_I&=M_I+Q_I,\\
 L_I&=1+\log(2+t-s),
 &K_I&=\mathfrak b_s+Q_I+S_I\beta_I .
 \end{aligned}
\end{equation}
By \eqref{eq:time-energy-equivalence}, \eqref{eq:def-B} and
\eqref{eq:feedbackbounds},
$\mathfrak b_s\le C_\delta\sqrt{B(s)}\le C_\delta\|Z(s)\|_{H^2}$.
The bounds \eqref{eq:time-window-estimate} and \eqref{eq:st-H-window},
together with \eqref{eq:time-B-integral} and \eqref{eq:global-H2}, give
\begin{equation}\label{eq:window-inputs}
 \begin{aligned}
 \Theta_I&\le C_\delta(\mathfrak b_s+S_I\beta_I),\qquad
 H_I^2\le C_\delta\beta_IQ_I,\\
 \beta_I+Q_I+H_I+M_I
 &\le C_\delta(\beta_I+\Theta_I)
 \le C_\delta\bigl(\mathfrak b_s+(1+S_I)\beta_I\bigr)\\
 &\le C_\delta(1+S_I)\|Z(s)\|_{H^2}\le C_\delta\eps.
 \end{aligned}
\end{equation}
The constants are independent of $s$ and $t$.

In the outward coordinate $x=j\xi$, one block of
\eqref{eq:comparison-auxiliary-generator} has the form
\[
 z_t=D z_{xx}-\operatorname{diag}(a_1,a_2,a_3)z_x+\mathsf Q z,
 \qquad a_a=-v_{e,a},\qquad \mathsf Q\mathbf1=0 ,
\]
where $D$ is bounded and uniformly positive and $\mathsf Q$ has
positive off-diagonal entries. At the parameter in
\eqref{eq:comparison-natural-parameter}, the coefficients are continuous
in time and smooth in space, with bounded spatial derivatives. Each
scalar parabolic propagator with diagonal term $\mathsf Q_{aa}\le0$
has a nonnegative kernel of mass at most one. Iterating Duhamel's
formula with the bounded nonnegative off-diagonal entries of
$\mathsf Q$ gives the system kernel. On each finite interval the
series converges in the bounded-function norm.

Fix the output component, time and position. Write $k_a(\tau,y)\ge0$
for the three input kernels. Since the system preserves the
constant vector $\mathbf1$, the sum $k=\sum_a k_a$ has mass one.
The input-variable adjoint equations give
\begin{equation}\label{eq:actual-total-kernel}
 \begin{aligned}
 -k_\tau
 &=\sum_a\left[(Dk_a)_{yy}+(a_a k_a)_y
                       +\sum_\ell \mathsf Q_{\ell a}k_\ell\right]
   =(Dk)_{yy}+(\bar a k)_y,\\
 \bar a&=\frac{\sum_a a_a k_a}{k}.
 \end{aligned}
\end{equation}
Where $k=0$, take $\bar a$ in the convex hull of the three speeds.

Let $m$ be a fixed nonnegative Schwartz kernel of mass one.
Convolving \eqref{eq:actual-total-kernel} in $y$, we find that
$\widetilde k=k*m$ satisfies
\begin{equation}\label{eq:filtered-FP}
 \begin{aligned}
 -\widetilde k_\tau
 &=\bigl((Dk)*m\bigr)_{yy}+\bigl((\bar a k)*m\bigr)_y
   =(\widetilde D\widetilde k)_{yy}+(\widetilde a\widetilde k)_y,\\
 \widetilde D&=\frac{(Dk)*m}{k*m},\qquad
 \widetilde a=\frac{(\bar a k)*m}{k*m},\qquad
 \|\widetilde k\|_\infty\le\|m\|_\infty .
 \end{aligned}
\end{equation}
The convolution preserves mass and $\widetilde D$ has the
ellipticity bounds of $D$. Where $\widetilde k=0$, take $\widetilde D=D$
and $\widetilde a=\bar a$. At the output time $k$ is a point mass
at the output position. Convolution defines the filtered kernel at
that time as well.

For the physical branch $i\ne e$, use $\lambda_i(q_j)=i c_j$ and
$a_e(\vartheta)=-j e c(\vartheta)$ in
\eqref{eq:comparison-auxiliary-generator}. By
\eqref{eq:lift-node-positivity},
\[
 \begin{aligned}
 \bigl|j(\lambda_i(q_j)-b)-a_a-j(i-e)c_j\bigr|
 &=\bigl|je(c_j-c(\vartheta))+x_{e,a}(\vartheta)\bigr|\\
 &\le |c_j-c(\vartheta)|+C_\delta\sqrt{\rho}
 \le\frac{c_j}{2},
 \end{aligned}
\]
Take the parameter neighborhood and the switching rate $\rho$ in
Lemma~\ref{lemma:lift} sufficiently small for this bound. The
coefficients $\bar a$ and $\widetilde a$ in
\eqref{eq:actual-total-kernel} and \eqref{eq:filtered-FP} are averages
with nonnegative weights. The speed $j(\lambda_i(q_j)-b)$ is constant
in $y$. Accordingly, $j(\lambda_i(q_j)-b)-\bar a$ and
$j(\lambda_i(q_j)-b)-\widetilde a$ both have the sign of $j(i-e)$ and
absolute value at least $\frac{c_j}{2}$.

\begin{lemma}[A cumulative-kernel estimate]\label{lemma:cumulative}
Let $k$ be nonnegative and uniformly bounded with unit mass. Assume
\eqref{eq:actual-total-kernel} holds with bounded $D$ and $\bar a$.
Let $f\in L^\infty(I;L^2)\cap L^2(I;H^1)$ and
$H,g\in L^2(I\times\R)$ satisfy
\[
 f_\tau+a_i(\tau)f_y=H_y+g
\]
in distributions, with $a_i$ bounded. If $a_i-\bar a$ has a fixed
sign and absolute value at least $g_0>0$, then
\begin{equation}\label{eq:cumulative-estimate}
 \int_I\!\!\int k f^2
 \le C\left[
 \sup_I\|f\|_2^2+\int_I(\|f_y\|_2^2+\|H\|_2^2)
                  +\int_I\!\!\int|fg|\right].
\end{equation}
The constant depends only on $g_0$, the indicated coefficient bounds
and $\|k\|_\infty$.
\end{lemma}
\begin{proof}
Set $\Phi(\tau,y)=\int_{-\infty}^y k(\tau,z)\dd z$. Then
$0\le\Phi\le1$, $\Phi_y=k$, and \eqref{eq:actual-total-kernel} gives
$\Phi_\tau=-(Dk)_y-\bar a k$ in distributions. To eliminate the spatial
integration constant, translate a test function toward $y$ tending
to $-\infty$. There $\Phi$ tends to zero, and the tails of
$Dk,\bar a k\in L^\infty(I;L^1)$ vanish in the test.

Mollify $\Phi$ in $y$ with a smooth nonnegative kernel of compact
support and unit mass. At fixed mollification scale the weight and
its first derivatives are bounded. Since the equation gives
$f_\tau\in L^2(I;H^{-1})$, use the mollified weight times $f$ in the
$L^2$ energy identity. Integrate the convolved $(Dk)_y$ term by parts
in $y$. The convolutions of $k$, $Dk$ and $\bar a k$ converge almost
everywhere with uniform bounds. The convolution of $\Phi$
converges uniformly in $y$ because $\Phi_y$ is bounded. The products
$f^2,ff_y,fH,f_yH,fg$ are integrable. Passing to the limit at
Lebesgue times $s,t$ gives
\begin{equation}\label{eq:cumulative-identity}
 \begin{aligned}
 \frac{1}{2}\int_I\!\!\int(a_i-\bar a)kf^2
 &=\left[\frac{1}{2}\int\Phi f^2\right]_s^t\\
 &\quad+\int_I\!\!\int\bigl[kf(H-Df_y)+\Phi f_yH-\Phi fg\bigr].
 \end{aligned}
\end{equation}
Apply Cauchy--Schwarz and Young's inequalities with the fixed sign
of $a_i-\bar a$ and the uniform bound on $k$. This gives
\[
 \begin{aligned}
 \frac{g_0}{2}\int_I\!\!\int kf^2
 &\le \frac{1}{2}\left|\int_I\!\!\int(a_i-\bar a)kf^2\right|\\
 &\le \sup_I\|f\|_2^2
   +C\left(\int_I\!\!\int kf^2\right)^{\frac{1}{2}}
       \left(\int_I(\|f_y\|_2^2+\|H\|_2^2)\right)^{\frac{1}{2}}\\
 &\quad+\frac{1}{2}\int_I(\|f_y\|_2^2+\|H\|_2^2)+\int_I\!\!\int|fg|\\
 &\le\frac{g_0}{4}\int_I\!\!\int kf^2+\sup_I\|f\|_2^2
     +C\int_I(\|f_y\|_2^2+\|H\|_2^2)+\int_I\!\!\int|fg|.
 \end{aligned}
\]
Absorbing the first term on the right gives
\eqref{eq:cumulative-estimate}.

Let the endpoint Lebesgue times tend to the ends of $I$.
\end{proof}

Let $T$ be a fixed spatial filter with Schwartz kernel $m_T$ and
choose an even nonnegative Schwartz kernel $m$ of mass one with
$|m_T|\le C_Tm$. The Cauchy--Schwarz inequality and Fubini's theorem
give
\[
 \begin{aligned}
 \int_{\R}k|Tf|^2
 &\le\|m_T\|_1\int_{\R}k(|m_T|*|f|^2)
 \le C_T\|m_T\|_1\int_{\R}k(m*|f|^2)\\
 &=C_T\|m_T\|_1\int_{\R}(k*m)|f|^2.
 \end{aligned}
\]
The filtered kernel satisfies \eqref{eq:filtered-FP}, so
Lemma~\ref{lemma:cumulative} applies with the signed speed gap above.

\begin{lemma}[Entropy and time-derivative potentials]\label{lemma:transverse}
Fix an end and an acoustic output family $e\in\{+1,-1\}$ of the
positive system. Evaluate the coefficients at a fixed spatial
mollification of the corrected field $U^j$, and extend them boundedly
to the whole line as below. For each fixed spatial filter $T$,
the following bounds are uniform in output time, position and
microstate.
\begin{align}
 \int_I\!\!\int k\,|T(\mathbf1_{\R_j}U^j_0)|^2
     &\le C_\delta\eps^2,\label{eq:entropy-potential}\\
 \int_I\!\!\int k\,|T(\mathbf1_{\R_j}U^j_{-e,t})|^2
 +\int_I\!\!\int k\,|T(\mathbf1_{\R_j}U^j_{-e,\xi})|^2
     &\le C_\delta K_I^2,\label{eq:time-gradient-potential}\\
 \int_I\!\!\int k\,|T(\mathbf1_{\R_j}U^j_{\xi\xi})|^2
     &\le C_\delta Q_I^2.\label{eq:second-gradient-potential}
\end{align}
On the left half-line one may reflect to the outward coordinate.
Time integrals stop at the chosen output time.
\end{lemma}
\begin{proof}
The signed speed gap following \eqref{eq:filtered-FP} applies to
$i=0$ and $i=-e$ in the outward coordinate $x=j\xi$.

\emph{The entropy amplitude.}
Write $Q^j=(v^j,u^j,\theta^j)$. On $\R_j$, \eqref{eq:normalform-definition}
gives $\eta=U^j_0=c_v\log\left(\frac{\theta^j}{\theta_j}\right)
+R\log\left(\frac{v^j}{v_j}\right)$. By \eqref{eq:corrected-end-PDE}
and $p(Q^j)=\frac{R\theta^j}{v^j}$,
\begin{equation}\label{eq:corrected-entropy}
 \begin{aligned}
 \eta_t-b\eta_\xi
 &=\left(\frac{R}{v^j}-\frac{p(Q^j)}{\theta^j}\right)u^j_\xi
    +\frac{1}{\theta^j}\left(\frac{\kappa\theta^j_\xi}{v^j}\right)_\xi
    +\frac{\mu|u^j_\xi|^2}{v^j\theta^j}+\frac{RG^j_1}{v^j}+\frac{c_vG^j_3}{\theta^j}\\
 &=(H_\eta)_\xi+P_\eta+\frac{RG^j_1}{v^j}+\frac{c_vG^j_3}{\theta^j},\\
 H_\eta&=\frac{\kappa\theta^j_\xi}{v^j\theta^j},\qquad
 P_\eta=\frac{\kappa|\theta^j_\xi|^2}{v^j(\theta^j)^2}
                +\frac{\mu|u^j_\xi|^2}{v^j\theta^j}.
 \end{aligned}
\end{equation}
By \eqref{eq:corrected-end-field} and \eqref{eq:st-parameter-tails},
$|\eta|\le C_\delta(|Z|+|\alpha|e^{-2\eta_0|\xi|})$.
Since $Q^j_\xi=r_\xi$, we also have
$|\eta_\xi|+|H_\eta|+\sqrt{P_\eta}\le C_\delta|r_\xi|$.
On $\R_j$, \eqref{eq:st-moment-bounds},
\eqref{eq:st-spatial-estimates}, \eqref{eq:time-energy-equivalence}
and \eqref{eq:end-source-bounds} give
\[
 \begin{aligned}
 \|\eta\|_2&\le C_\delta\|Z\|_2\le C_\delta\eps,
 &\|\eta\|_\infty&\le C_\delta S_I,\\
 \|\sqrt\omega\eta\|_2+\|r_\xi\|_2
 &\le C_\delta\bigl(m_\omega(Z)+\|F\|_2\bigr)\le C_\delta\sqrt B,
 &\left\|\frac{G^j}{\sqrt\omega}\right\|_2&\le C_\delta\Gamma .
 \end{aligned}
\]

Take a smooth $0\le\chi\le1$ with $\chi=0$ for $j\xi\le1$ and
$\chi=1$ for $j\xi\ge2$. Multiplying \eqref{eq:corrected-entropy}
by $\chi$ gives
\[
 \begin{aligned}
 (\chi\eta)_t-b(\chi\eta)_\xi
 &=(\chi H_\eta)_\xi
       +\chi\bigl(P_\eta+\frac{RG^j_1}{v^j}+\frac{c_vG^j_3}{\theta^j}\bigr)-\chi'(b\eta+H_\eta).
 \end{aligned}
\]
Extend the cutoff field and flux by zero. The outward coordinate
$y=j\xi$ gives speed $-jb$ and flux $j\chi H_\eta$.
The source terms belong to $L^2(I\times\R)$ by
\eqref{eq:end-source-bounds} and
\[
 \int_I\|P_\eta\|_2^2
 \le C_\delta\sup_I\|r_\xi\|_2^2\int_I\|r_\xi\|_\infty^2
 \le C_\delta M_I^2 H_I^2.
\]
Use $|\chi'|+|\mathbf1_{\R_j}-\chi|^2\le C_\delta\omega$ and the
filter bound. Lemma~\ref{lemma:cumulative} with $\widetilde k=k*m$
then gives
\begin{align*}
 &\int_I\!\!\int k\,|T(\mathbf1_{\R_j}\eta)|^2\\*
 &\quad\le C_\delta\left[\int_I\!\!\int\widetilde k|\chi\eta|^2
                    +\int_I\|\sqrt\omega\eta\|_2^2\right]\\
 &\quad\le C_\delta\Bigl[\sup_I\|\chi\eta\|_2^2+
     \int_I\bigl(\|(\chi\eta)_\xi\|_2^2+\|\chi H_\eta\|_2^2
                                  +\|\sqrt\omega\eta\|_2^2\bigr)\Bigr]\\*
 &\qquad+C_\delta\int_I\!\!\int_{\R_j}|\chi\eta|
      \bigl|\chi(P_\eta+\frac{RG^j_1}{v^j}+\frac{c_vG^j_3}{\theta^j})-\chi'(b\eta+H_\eta)\bigr|\\
 &\quad\le C_\delta\Bigl[\sup_I\|\eta\|_2^2+
       \int_I\Bigl(\|\sqrt\omega\eta\|_2^2
                  +(1+\|\eta\|_\infty)\|r_\xi\|_2^2
                  +\|\sqrt\omega\eta\|_2
                       \left\|\frac{G^j}{\sqrt\omega}\right\|_2\Bigr)\Bigr]\\
 &\quad\le C_\delta\bigl(\eps^2+(1+S_I)\beta_I^2+\beta_IQ_I\bigr)
       \le C_\delta\eps^2 .
\end{align*}
The last line uses \eqref{eq:window-inputs} and
$S_I\le C_\delta\eps$ from \eqref{eq:global-H2}, proving
\eqref{eq:entropy-potential}.

\emph{Time derivatives.}
In the conserved variables $\mathcal U=(v,u,E)$, where
$E=c_v\theta+\frac{u^2}{2}$, set $P=\mathcal U_t=D_q\mathcal U(q)F$.
Multiply the momentum equation by $u$ and add it to the temperature
equation in \eqref{eq:NSF}. The right-hand side is the derivative of
$\frac{\kappa\theta_\xi+\mu u u_\xi}{v}$. In the moving coordinate,
\[
 \begin{gathered}
 \mathcal U_t-b\mathcal U_\xi+\partial_\xi(-u,p,pu)^{\mathsf T}
  =\partial_\xi\bigl(\mathsf B_c(\mathcal U)\mathcal U_\xi\bigr),\\
 \mathsf B_c(\mathcal U)=
 \begin{pmatrix}
 0&0&0\\[1mm]
 0&\frac{\mu}{v}&0\\[1mm]
 0&\frac{u\bigl(\mu-\frac{\kappa}{c_v}\bigr)}{v}
       &\frac{\kappa}{c_vv}
 \end{pmatrix}.
 \end{gathered}
\]
Write $\mathsf A_c$ for the Jacobian of $(-u,p,pu)$ in $\mathcal U$ and
$\mathsf A_{c,j}$ for its value at $q_j$. Differentiate in time and
use $b_t=\ell(F)$ from \eqref{eq:time-derivatives} to get
\begin{equation}\label{eq:conserved-time}
 \begin{aligned}
 P_t+(\mathsf A_{c,j}-bI)P_\xi
 &=\partial_\xi\bigl[\mathsf B_cP_\xi
       +D\mathsf B_c[P]\mathcal U_\xi
       -(\mathsf A_c-\mathsf A_{c,j})P\bigr]+b_t\mathcal U_\xi .
 \end{aligned}
\end{equation}
The coefficient bounds and \eqref{eq:time-energy-equivalence} give
$\sup_I\|P\|_2\le C_\delta M_I$. In
$P_\xi=D_q\mathcal U(q)F_\xi+D_q^2\mathcal U(q)[q_\xi,F]$,
write $q_\xi=\bar q'+Z_\xi$. On $\R_j$,
\eqref{eq:profilebounds} gives
$|\bar q'|+|\bar q-q_j|\le C_\delta\omega$, hence
$|\mathsf A_c-\mathsf A_{c,j}|\le C_\delta(\omega+|Z|)$.
With $\|Z_\xi\|_\infty^2\le C B$, this yields
\begin{equation}\label{eq:time-potential-input}
 \begin{aligned}
 &\|P_\xi\|_{L^2(I\times\R)}
 \le C_\delta\left[\int_I\bigl(\Gamma^2
                 +\|Z_\xi\|_\infty^2\|F\|_2^2\bigr)\right]^{\frac{1}{2}}
 \le C_\delta(Q_I+\beta_IM_I),\\
 &\|D\mathsf B_c[P]\mathcal U_\xi
       -(\mathsf A_c-\mathsf A_{c,j})P\|_{L^2(I\times\R_j)}\\
 &\quad\le C_\delta\left[\int_I\Bigl(m_\omega(F)^2+
                (\|Z_\xi\|_\infty^2+S_I^2)\|F\|_2^2\Bigr)\right]^{\frac{1}{2}}
 \le C_\delta(Q_I+\beta_IM_I+S_I\beta_I).
 \end{aligned}
\end{equation}

At the end state, \eqref{eq:affine-diagonalizer} gives the rows of
$L_c^j=DN_{q_j}(q_j)\bigl[D_q\mathcal U(q_j)\bigr]^{-1}$
as left eigenvectors of $\mathsf A_{c,j}$, in the order $(+,0,-)$.
For $i\ne e$, set $f=\chi(L_c^jP)_i$ with the same cutoff $\chi$.
Then $f_t+(\lambda_i(q_j)-b)f_\xi=H_\xi+g$, with
\[
 \begin{aligned}
 H&=\chi\bigl[L_c^j\bigl(\mathsf B_cP_\xi
      +D\mathsf B_c[P]\mathcal U_\xi-(\mathsf A_c-\mathsf A_{c,j})P\bigr)\bigr]_i,\\
 g&=\chi b_t(L_c^j\mathcal U_\xi)_i
      +(\lambda_i(q_j)-b)\chi'(L_c^jP)_i\\
 &\quad-\chi'\bigl[L_c^j\bigl(\mathsf B_cP_\xi
      +D\mathsf B_c[P]\mathcal U_\xi-(\mathsf A_c-\mathsf A_{c,j})P\bigr)\bigr]_i.
 \end{aligned}
\]
Extend the fields by zero and use $y=j\xi$, with transport speed
$j(\lambda_i(q_j)-b)$ and flux $jH$. The source $g$ belongs to
$L^2(I\times\R)$ by \eqref{eq:time-potential-input},
$b_t\in L^2(I)$ and $\sup_I\|\mathcal U_\xi\|_2<\infty$. With
$|\chi'|\le C_\delta\omega$, $|\bar q'|\le C_\delta\omega$ and the bound
$|b_t|\le C_\delta m_\omega(F)$ from \eqref{eq:feedbackbounds},
we have
\[
 \begin{aligned}
 \int_I\!\!\int |fg|
 &\le C_\delta\int_I\Bigl[
 |b_t|\bigl(m_\omega(F)+\|F\|_2\|Z_\xi\|_2\bigr)+m_\omega(F)^2\\
 &\qquad+m_\omega(F)\bigl(\|P_\xi\|_2+
 \|D\mathsf B_c[P]\mathcal U_\xi
           -(\mathsf A_c-\mathsf A_{c,j})P\|_{L^2(\R_j)}\bigr)\Bigr]\\
 &\le C_\delta\bigl[Q_I^2+M_I\beta_IQ_I
                  +Q_I(Q_I+\beta_IM_I+S_I\beta_I)\bigr]
 \le C_\delta K_I^2 .
 \end{aligned}
\]
We used \eqref{eq:time-potential-input} and \eqref{eq:window-inputs}.
Applying Lemma~\ref{lemma:cumulative} and
$|\mathbf1_{\R_j}-\chi|^2\le C_\delta\omega$ gives
\[
 \begin{aligned}
 &\int_I\!\!\int \widetilde k\,|\mathbf1_{\R_j}(L_c^jP)_i|^2\\
 &\quad\le C_\delta\left[\int_I\!\!\int\widetilde k |f|^2
                                +\int_I m_\omega(F)^2\right]\\
 &\quad\le C_\delta\left[M_I^2+(Q_I+\beta_IM_I+S_I\beta_I)^2
                               +\int_I\!\!\int |fg|\right]
 \le C_\delta K_I^2 .
 \end{aligned}
\]

On $\R_j$, \eqref{eq:corrected-end-field} and
\eqref{eq:st-parameter-tails} give
\[
 \begin{aligned}
 &|DN_{q_j}(q_j)F-L_c^jP|+|U^j_t-DN_{q_j}(q_j)F|\\
 &\quad\le C_\delta\bigl[(|q-q_j|+|Q^j-q_j|)|F|
                            +e^{-2\eta_0|\xi|}|\alpha'|\bigr]\\
 &\quad\le C_\delta(S_I|F|+\omega|F|+\omega|\alpha'|).
 \end{aligned}
\]
By the filter bound with $\widetilde k=k*m$ and
$|\alpha'|\le C_\delta m_\omega(F)$ from \eqref{eq:st-moment-bounds}, we have
\[
 \begin{aligned}
 \int_I\!\!\int k\,|T(\mathbf1_{\R_j}U^j_{-e,t})|^2
 &\le C_\delta\int_I\!\!\int\widetilde k\,
                                 |\mathbf1_{\R_j}U^j_{-e,t}|^2\\
 &\le C_\delta\int_I\!\!\int\widetilde k\,
                          |\mathbf1_{\R_j}(L_c^jP)_{-e}|^2\\
 &\quad+C_\delta\int_I\bigl(S_I^2\|F\|_2^2+m_\omega(F)^2+|\alpha'|^2\bigr)\\
 &\le C_\delta(K_I^2+S_I^2\beta_I^2+Q_I^2)\le C_\delta K_I^2 .
 \end{aligned}
\]
The time-derivative bound in \eqref{eq:time-gradient-potential}
also holds for the $-e$ component of $DN_{q_j}(q_j)F$.

\emph{Spatial derivatives.}
For the reference state $q_j$, $\mathsf A(0)=\operatorname{diag}(c_j,0,-c_j)$.
The $i$-th component of \eqref{eq:corrected-end-PDE} gives
\begin{equation}\label{eq:transverse-spatial-recovery}
 \begin{aligned}
 (\lambda_i(q_j)-b)U^j_{i,\xi}
 &=\bigl[\mathsf B(0)U^j_{\xi\xi}-U^j_t\bigr]_i\\
 &\quad+\bigl[(\mathsf A(0)-\mathsf A(U^j))U^j_\xi
              +(\mathsf B(U^j)-\mathsf B(0))U^j_{\xi\xi}\bigr]_i\\
 &\quad+\bigl[\mathsf C(U^j)[U^j_\xi,U^j_\xi]+DN_{q_j}(Q^j)G^j\bigr]_i.
 \end{aligned}
\end{equation}
By \eqref{eq:lax}, \eqref{eq:feedbackbounds} and \eqref{eq:global-H2},
\[
 |\lambda_i(q_j)-b|
 \ge |\lambda_i(q_j)-\sigma|-|b-\sigma|
 \ge \frac{1}{2}\min\{\sigma-c_+,c_--\sigma\}>0
\]
for sufficiently small $\eps$. The bounds \eqref{eq:st-family-bounds},
\eqref{eq:st-residual} and \eqref{eq:corrected-end-field} give
$\|U^j\|_\infty+\|r\|_\infty\le C_\delta S_I$. Hence
$|\mathsf A(0)-\mathsf A(U^j)|+|\mathsf B(U^j)-\mathsf B(0)|
\le C_\delta S_I$.

The identity $U^j_\xi=DN_{q_j}(Q^j)r_\xi$ bounds the quadratic term in
\eqref{eq:transverse-spatial-recovery} by $C_\delta|r_\xi|^2$.
Apply the filter bound to the zero extension of $U^j_{-e,\xi}$. In
\eqref{eq:transverse-spatial-recovery}, use the bound on $U^j_{-e,t}$
with $\widetilde k=k*m$ and apply \eqref{eq:end-field-derivatives},
\eqref{eq:end-source-bounds} and \eqref{eq:st-GN-bounded}.
This gives
\[
 \begin{aligned}
 &\int_I\!\!\int k\,|T(\mathbf1_{\R_j}U^j_{-e,\xi})|^2\\
 &\quad\le C_\delta\int_I\!\!\int\widetilde k\,
                         |\mathbf1_{\R_j}U^j_{-e,t}|^2\\
 &\qquad+C_\delta\int_I\Bigl[(1+S_I^2)\|U^j_{\xi\xi}\|_2^2
                +S_I^2\|U^j_\xi\|_2^2+\|r_\xi\|_4^4
                +\|DN_{q_j}(Q^j)G^j\|_2^2\Bigr]\\
 &\quad\le C_\delta\Bigl[
       K_I^2+(1+S_I^2)Q_I^2+S_I^2\beta_I^2
              +\sup_I\|r\|_\infty^2\int_I\|r_{\xi\xi}\|_2^2\Bigr]\\
 &\quad\le C_\delta\bigl[K_I^2+(1+S_I^2)Q_I^2+S_I^2\beta_I^2\bigr]
 \le C_\delta K_I^2 .
 \end{aligned}
\]
The norms of $U^j$ and $G^j$ are on $\R_j$ and those of $r$ are
on $\R$. This proves the spatial bound in
\eqref{eq:time-gradient-potential}. For the second derivative,
\eqref{eq:end-field-derivatives} gives
\[
 \int_I\!\!\int k\,|T(\mathbf1_{\R_j}U^j_{\xi\xi})|^2
 \le C_\delta\int_I\|U^j_{\xi\xi}\|_{L^2(\R_j)}^2
 \le C_\delta Q_I^2,
\]
This is \eqref{eq:second-gradient-potential}.
\end{proof}

\begin{corollary}[The affine zero-frequency chain term]
\label{cor:zero-chain}
In the outward coordinate $x=j\xi$, set $U=U^j$ on $\R_+$ and
$\vartheta^\natural=\varphi*\widetilde U$ as in
\eqref{eq:comparison-natural-parameter}. Let $P^0$ solve the whole-line
system for $\mathscr C_{\vartheta^\natural}$ with zero initial data
at $s$, zero contact source, and acoustic source
\[
 \begin{aligned}
 \bigl[(\partial_t-\mathscr C_{\vartheta^\natural})P^0\bigr]_{e,a}
 &= k_0\mathbf1_{\R_+}\bigl[\chi(D)(\mathbf1_{\R_+}U_0)\bigr]\\
 &\quad\times\bigl[\varphi*(\mathbf1_{\R_+}U_{-e,t})
    -v_{e,a}\varphi*(\mathbf1_{\R_+}U_{-e,x})
    -D_e\varphi*(\mathbf1_{\R_+}U_{-e,xx})\bigr].
 \end{aligned}
\]
Then $\|P^0\|_{L^\infty(I\times\R)}\le C_\delta\eps K_I$. For every
$N>1$ the difference between \eqref{eq:comparison-zero-chain} and
this source on $x>0$ is bounded by $C_{\delta,N}S_IM_I(1+x)^{-N}$.
With zero initial data, a source given by a product of two filtered
first spatial derivatives of $\widetilde U$ with bounded coefficients
produces a solution bounded by $C_\delta H_I^2$.
\end{corollary}
\begin{proof}
Fix the output time $\tau$, position $x$ and component in acoustic
block $e$. Write $k$ for its total kernel. Uniform boundedness of
$v_{e,a}$ and $D_e$, positivity and Cauchy--Schwarz together with
\eqref{eq:entropy-potential}--\eqref{eq:second-gradient-potential}
give
\[
 \begin{aligned}
 &|P^0_{e,a}(\tau,x)|\\
 &\quad\le\int_s^\tau\!\!\int\sum_\ell k_\ell
               \bigl|[(\partial_t-\mathscr C_{\vartheta^\natural})P^0]_{e,\ell}\bigr|\\
 &\quad\le C_\delta
 \left(\int_s^\tau\!\!\int k\,
         |\chi(D)(\mathbf1_{\R_+}U_0)|^2\right)^{\frac{1}{2}}\\
 &\qquad\times
 \left(\int_s^\tau\!\!\int k\,\bigl[
       |\varphi*(\mathbf1_{\R_+}U_{-e,t})|^2
      +|\varphi*(\mathbf1_{\R_+}U_{-e,x})|^2
      +|\varphi*(\mathbf1_{\R_+}U_{-e,xx})|^2\bigr]\right)^{\frac{1}{2}}\\
 &\quad\le C_\delta\eps(K_I+Q_I)\le C_\delta\eps K_I .
 \end{aligned}
\]
Each filtered first spatial derivative is bounded by $C_\delta h$
by \eqref{eq:comparison-extension-bounds}. Since $\int k=1$, the
response to their product is bounded by
$C_\delta\int_s^\tau h(\rho)^2\dd\rho\le C_\delta H_I^2$.

The extension matches $U$ and $U_x$ at $x=0$, so the derivatives
of $\vartheta^\natural$ in \eqref{eq:comparison-zero-chain} are the
convolutions of $\varphi$ with
$\widetilde U_t,\widetilde U_x,\widetilde U_{xx}$. The reflected
time derivatives on $(-1,0)$ and
$|\bar U_t|\le\|\varrho\|_2\|U_t\|_{L^2(0,1)}$ control the local
norm of $\widetilde U_t$. Using
\eqref{eq:comparison-extension-bounds} and differentiating
\eqref{eq:corrected-end-field} in time gives
\[
 \begin{aligned}
 \|\widetilde U_x\|_2+\|\widetilde U_{xx}\|_2
      +\sup_{z\in\R}\|\widetilde U_t\|_{L^2(z,z+1)}
 &\le C_\delta\bigl(\sqrt{\mathcal T}+\|U_t\|_{L^2(\R_+)}\bigr)\\
 &\le C_\delta\bigl(\sqrt{\mathcal T}+\|F\|_2+|\alpha'|\bigr)
 \le C_\delta M_I .
 \end{aligned}
\]
For the last inequality use \eqref{eq:time-energy-equivalence} and
\eqref{eq:st-moment-bounds}. For $f=\widetilde U_t,\widetilde U_x$ or
$\widetilde U_{xx}$, Cauchy--Schwarz on the compact support of
$\varphi$ gives, for $x>0$,
\[
 |\varphi*f(\rho,x)|\le C_\delta M_I,\qquad
 |\varphi*(\mathbf1_{\R_-}f)(\rho,x)|
     \le C_{\delta,N}M_I(1+x)^{-N}.
\]
In the second bound the convolution restricted to $x>0$ is supported
in a fixed interval. Write $m_\chi$ for the Schwartz kernel of
$\chi(D)$. The bound on $\|\widetilde U\|_\infty$ in
\eqref{eq:comparison-extension-bounds} yields
\[
 \begin{aligned}
 |\chi(D)(\mathbf1_{\R_-}\widetilde U_0)(\rho,x)|
 &\le \|\widetilde U(\rho)\|_\infty\int_x^\infty |m_\chi(z)|\dd z
 \le C_{\delta,N}S_I(1+x)^{-N}.
 \end{aligned}
\]
Since $\|\chi(D)(\mathbf1_{\R_+}U_0)\|_\infty\le C_\delta S_I$,
subtracting the two products on $x>0$ gives
\[
 \begin{aligned}
 \bigl|R^0_{e,a}
       -[(\partial_t-\mathscr C_{\vartheta^\natural})P^0]_{e,a}\bigr|
 &\le C_{\delta,N}
       \bigl[S_I(1+x)^{-N}M_I+S_IM_I(1+x)^{-N}\bigr]\\
 &\le C_{\delta,N}S_IM_I(1+x)^{-N}.
 \end{aligned}
\]
\end{proof}

\subsection{Comparison on a finite time interval}\label{sec:comparison}

For $I=[s,t]$, use \eqref{eq:time-window-notation}, \eqref{eq:st-h}
and \eqref{eq:window-notation}--\eqref{eq:window-inputs}. Choose $s$
where the endpoint spatial estimates hold. Constants are
independent of the length and position of $I$. They may depend on
the fixed shock, the cutoffs and the switching rate from
Lemma~\ref{lemma:lift}.

\begin{proposition}[Finite-interval comparison]\label{prop:comparison}
For sufficiently small initial $H^2$ norm, the solution satisfies
\begin{equation}\label{eq:comparison-final}
 S_I\le C_\delta\left\{
 S(s)+\mathfrak b_s+\Theta_I+H_I^2
       +H_I\mathfrak b_s+L_IH_IQ_I\right\}.
\end{equation}
The proof uses Proposition~\ref{prop:spatial} for spatial recovery,
Lemma~\ref{lemma:incoming} for the boundary traces and
Lemma~\ref{lemma:transverse} for the variable-coefficient potentials.
\end{proposition}

\subsubsection{Extension and frozen equations}

On one half-line use $x=j\xi>0$ and set $U=U^j$ from
\eqref{eq:corrected-end-field}. The frozen physical generator is
\begin{equation}\label{eq:comparison-frozen-generator}
 \mathscr A_\vartheta(\tau,D)
 =j\{b(\tau)I-\mathsf A(\vartheta)\}\partial_x
                  +\mathsf B(\vartheta)\partial_x^2.
\end{equation}
Hold the common drift $jb(\tau)$ fixed in parameter differentiation.
This term translates the frozen kernels. The auxiliary operator is
\begin{equation}\label{eq:comparison-auxiliary-generator}
 \mathscr C_\vartheta(\tau,D)
 =\operatorname{diag}_e\left\{
 D_e(\vartheta)\partial_x^2 I_3
 +\operatorname{diag}_{a}\bigl(v_{e,a}(\vartheta,\tau)\bigr)
                                      \partial_x+\mathsf Q_e(\vartheta)
 \right\}.
\end{equation}
Here $v_{e,a}=jb(\tau)+a_e(\vartheta)+x_{e,a}(\vartheta)$, with the
reflection included in $a_e$. Each block has a common scalar
diffusion coefficient, and $\mathsf Q_e\mathbf1=0$.

Take a nonnegative smooth function $\varrho$ of unit integral
supported in $(0,1)$ and set
$\bar U(\tau)=\int_0^1\varrho(y)U(\tau,y)\dd y$.
Take a smooth cutoff $\theta_E$ on $[-1,0]$, equal to one near
zero and zero near $-1$. Set
\[
 \widetilde U(\tau,x)=
 \begin{cases}
 U(\tau,x),&x\ge0,\\
 \bar U(\tau)+\theta_E(x)
 \bigl[3U(\tau,-x)-2U(\tau,-2x)-\bar U(\tau)\bigr],&-1<x<0,\\
 \bar U(\tau),&x\le-1.
 \end{cases}
\]
At $x=0$ the two traces from the extension interval are
\[
 \begin{pmatrix}\widetilde U(\tau,0-)\\ \widetilde U_x(\tau,0-)\end{pmatrix}
 =\begin{pmatrix}(3-2)U(\tau,0+)\\ (-3+4)U_x(\tau,0+)\end{pmatrix}
 =\begin{pmatrix}U(\tau,0+)\\ U_x(\tau,0+)\end{pmatrix}.
\]
The extension is constant near $x=-1$. The matched traces give
$\widetilde U\in H^2_{\mathrm{loc}}(\R)$. At a fixed time the weak
derivatives on $-1<x<0$ are
\[
 \begin{aligned}
 \widetilde U_x
 &=\theta_E'\bigl[3U(-x)-2U(-2x)-\bar U\bigr]
    +\theta_E\bigl[-3U_x(-x)+4U_x(-2x)\bigr],\\
 \widetilde U_{xx}
 &=\theta_E''\bigl[3U(-x)-2U(-2x)-\bar U\bigr]
    +2\theta_E'\bigl[-3U_x(-x)+4U_x(-2x)\bigr]\\
 &\quad+\theta_E\bigl[3U_{xx}(-x)-8U_{xx}(-2x)\bigr].
 \end{aligned}
\]
For $0\le y\le2$ and $p=2,\infty$, the fundamental theorem of calculus
gives
\[
 \begin{aligned}
 |U(y)-\bar U|
 &\le\int_0^1\varrho(z)|U(y)-U(z)|\dd z
 \le\int_0^2|U_x(z)|\dd z
 \le 2^{1-\frac{1}{p}}\|U_x\|_{L^p(0,2)} .
 \end{aligned}
\]
Apply the same extension to $U_t$ to obtain $\widetilde U_t$, with
$\bar U_t=\int_0^1\varrho(y)U_t(y)\dd y$.
By \eqref{eq:corrected-end-field} and \eqref{eq:st-residual},
$\|U\|_\infty\le C_\delta S(\tau)$. Since $U^j_\xi=DN_{q_j}(Q^j)r_\xi$,
the weight $\omega$ on $|\xi|\le2$ gives
$\|U_x\|_{L^2(0,2)}\le C_\delta m_\omega(r_\xi)$.
Consequently, \eqref{eq:time-energy-equivalence},
\eqref{eq:st-spatial-estimates} and \eqref{eq:end-field-derivatives}
give
\begin{equation}\label{eq:comparison-extension-bounds}
 \begin{aligned}
 \|\widetilde U\|_\infty
 &\le C_\delta\|U\|_\infty\le C_\delta S(\tau),\\
 \|\widetilde U_x\|_\infty+\|\widetilde U_t\|_\infty
 &\le C_\delta\bigl(\|U_x\|_\infty+\|U_t\|_\infty\bigr)
 \le C_\delta h(\tau),\\
 \|\widetilde U_x\|_2+\|\widetilde U_{xx}\|_2
 &\le C_\delta\bigl(\|U_x\|_2+\|U_{xx}\|_2\bigr)
 \le C_\delta\sqrt{\mathcal T(\tau)},\\
 \|\widetilde U_x\|_{L^2(-1,0)}+\|\widetilde U_{xx}\|_2
 &\le C_\delta\bigl(\|U_x\|_{L^2(0,2)}+\|U_{xx}\|_2\bigr)\\
 &\le C_\delta\bigl(m_\omega(r_\xi)+\|U_{xx}\|_2\bigr)
 \le C_\delta\Gamma(\tau).
 \end{aligned}
\end{equation}
Unrestricted norms of $U$ are on $\R_+$ and those of
$\widetilde U$ on $\R$. The bounds involving $h$ and $\Gamma$ hold
almost everywhere in time.

Fix a nonnegative $\varphi\in C_c^\infty(\R)$ with integral one and set
\begin{equation}\label{eq:comparison-natural-parameter}
 \vartheta^\natural(\tau,x)=\varphi*\widetilde U(\tau,x).
\end{equation}
By \eqref{eq:comparison-extension-bounds},
$\|\vartheta^\natural(\tau)\|_\infty\le\|\widetilde U(\tau)\|_\infty
\le C_E S(\tau)$ for a fixed $C_E$. We take parameter suprema over
$|\vartheta|\le C_E S_I$. For sufficiently small $\eps$, this ball
lies in the fixed neighborhood of Lemma~\ref{lemma:lift}.
Since
\[
 \vartheta^\natural_x=\varphi*\widetilde U_x,\qquad
 \vartheta^\natural_{xx}=\varphi'*\widetilde U_x
                        =\varphi*\widetilde U_{xx},\qquad
 \vartheta^\natural_t=\varphi*\widetilde U_t,
\]
Young's inequality and \eqref{eq:comparison-extension-bounds} give
\begin{equation}\label{eq:comparison-parameter-derivatives}
 \begin{aligned}
 &\|\vartheta^\natural_x\|_{L^2(I;L^\infty)}
 +\|\vartheta^\natural_{xx}\|_{L^2(I;L^\infty)}
 +\|\vartheta^\natural_t\|_{L^2(I;L^\infty)}\\
 &\quad\le(\|\varphi\|_1+\|\varphi'\|_1)
                  \|\widetilde U_x\|_{L^2(I;L^\infty)}
         +\|\varphi\|_1\|\widetilde U_t\|_{L^2(I;L^\infty)}
 \le C_\delta H_I .
 \end{aligned}
\end{equation}
The last bound in \eqref{eq:comparison-extension-bounds} also gives
$\|\widetilde U_{xx}\|_{L^2(I;L^2)}\le C_\delta Q_I$.

In \eqref{eq:comparison-frozen-generator}, evaluate the coefficients
at $\widetilde U(\tau,x)$ to define $\mathscr A_{\widetilde U}$.
The extended equation is
\begin{equation}\label{eq:comparison-source-split}
 \begin{aligned}
 \widetilde U_t-\mathscr A_\vartheta\widetilde U
 &=E_\vartheta=f_\vartheta+n+g+\bar U'(\tau)\mathbf1_{\{x<0\}},\\
 f_\vartheta&=(\mathscr A_{\widetilde U}-\mathscr A_\vartheta)\widetilde U\\
 &=-j\bigl[\mathsf A(\widetilde U)-\mathsf A(\vartheta)\bigr]\widetilde U_x
       +\bigl[\mathsf B(\widetilde U)-\mathsf B(\vartheta)\bigr]\widetilde U_{xx}.
 \end{aligned}
\end{equation}
On $x>0$, \eqref{eq:corrected-end-PDE} gives
$n=\mathsf C(U)[U_x,U_x]$ and $g=DN_{q_j}(Q^j)G^j$ in the coordinate
$x=j\xi$. Set $n=0$ on $x<0$. On $-1<x<0$ the extension gives
\[
 \begin{aligned}
 g
 &=\widetilde U_t-\bar U'(\tau)-\mathscr A_{\widetilde U}\widetilde U\\
 &=\theta_E\bigl[3U_t(-x)-2U_t(-2x)-\bar U'(\tau)\bigr]
                                  -\mathscr A_{\widetilde U}\widetilde U .
 \end{aligned}
\]
On $x\le-1$, $\widetilde U=\bar U(\tau)$ and $g=0$.
The matched traces give \eqref{eq:comparison-source-split}
in distributions on the whole line.

Use $|\bar U'|\le\|\varrho\|_2\|U_t\|_{L^2(0,1)}$, differentiate
\eqref{eq:corrected-end-field} in time, and apply
\eqref{eq:st-moment-bounds} and \eqref{eq:time-energy-equivalence}.
This gives
\[
 \begin{aligned}
 |\bar U'|+\|U_t\|_{L^2(0,2)}
 &\le C_\delta\bigl(m_\omega(F)+|\alpha'|\bigr)
 \le C_\delta m_\omega(F)
 \le C_\delta\min\{\sqrt{\mathcal T},\Gamma\}.
 \end{aligned}
\]
For $|\vartheta|\le C_E S_I$, we have
$\|\widetilde U-\vartheta\|_\infty\le C_\delta S_I$.
On $x>0$ the weight $(1+x)^N\sqrt{\omega(jx)}$ belongs to
$L^2\cap L^\infty$. Thus \eqref{eq:comparison-extension-bounds},
\eqref{eq:end-source-bounds}, and Sobolev embedding for $U_x$ give,
for each fixed $N\ge0$,
\begin{align}
 &\|f_\vartheta(\tau)\|_2
 \le C_\delta S_I\bigl(\|\widetilde U_x\|_2+\|\widetilde U_{xx}\|_2\bigr)
 \le C_\delta S_I\sqrt{\mathcal T(\tau)},\notag\\
 &\|n(\tau)\|_\infty\le C_\delta\|U_x\|_\infty^2
 \le C_\delta h(\tau)^2,\qquad
 \|n(\tau)\|_2\le C_\delta\|U_x\|_2\|U_x\|_{H^1(\R_+)}
 \le C_\delta\mathcal T(\tau),\notag\\
 &\|(1+|x|)^N g(\tau)\|_1+\|(1+|x|)^N g(\tau)\|_2+|\bar U'(\tau)|\notag\\*
 &\quad\le C_{\delta,N}\left[
       \left\|\frac{DN_{q_j}(Q^j)G^j}{\sqrt{\omega}}\right\|_{L^2(\R_j)}
       +\|U_t\|_{L^2(0,2)}+\|\widetilde U_x\|_{L^2(-1,0)}
                            +\|\widetilde U_{xx}\|_{L^2(-1,0)}\right]\notag\\*
 &\quad\le C_{\delta,N}\min\{\sqrt{\mathcal T(\tau)},\Gamma(\tau)\}.
 \label{eq:comparison-source-budgets}
\end{align}
Norms of $U$ are on $\R_+$ unless a smaller set is specified.
Bounds in terms of $\mathcal T$ include the endpoints. Those in
terms of $h$ or $\Gamma$ hold almost everywhere.
By \eqref{eq:time-energy-equivalence},
\[
 \|f_\vartheta\|_{L^2(I;L^2)}
 \le C_\delta S_I\left(\int_I\mathcal T(\tau)\dd\tau\right)^{\frac{1}{2}}
 \le C_\delta S_I\beta_I .
\]

For each fixed $N\ge0$ and almost every $\tau\in I$,
\begin{equation}\label{eq:comparison-actual-source}
 \left|\left.
 \bigl[\mathsf L(\vartheta,D)E_\vartheta\bigr](\tau,x)
                  \right|_{\vartheta=\vartheta^\natural(\tau,x)}\right|
 \le C_\delta h(\tau)^2+
                      C_{\delta,N}M_I(1+x)^{-N},\qquad x>0.
\end{equation}
Let $m_\vartheta$ be the convolution kernel of $\mathsf L(\vartheta,D)$.
Integration by parts in $k$ using \eqref{eq:lift-multiplier} gives
$|m_\vartheta(z)|+|m_\vartheta'(z)|\le C_{\delta,N}(1+|z|)^{-N}$
uniformly in $\vartheta$. This bound also holds for parameter
derivatives through order two.

Fix $\tau$ and omit it in the spatial arguments. With $\vartheta$
constant, integrate the diffusion term in
\eqref{eq:comparison-source-split} by parts in $y$ to get
\[
 \begin{aligned}
 \bigl[\mathsf L(\vartheta,D)f_\vartheta\bigr](x)
 &=-j\int_{\R} m_\vartheta(x-y)
           \bigl[\mathsf A(\widetilde U(y))-\mathsf A(\vartheta)\bigr]
                         \widetilde U_y(y)\dd y\\
 &\quad+\int_{\R} m_\vartheta'(x-y)
           \bigl[\mathsf B(\widetilde U(y))-\mathsf B(\vartheta)\bigr]
                         \widetilde U_y(y)\dd y\\
 &\quad-\int_{\R} m_\vartheta(x-y)
     \bigl(D\mathsf B(\widetilde U(y))[\widetilde U_y(y)]\bigr)
                         \widetilde U_y(y)\dd y .
 \end{aligned}
\]
The traces of $\widetilde U$ and $\widetilde U_x$ match at both
interfaces, so the integration by parts holds across them. Now set
$\vartheta=\vartheta^\natural(\tau,x)$. Equations
\eqref{eq:comparison-natural-parameter} and
\eqref{eq:comparison-extension-bounds} give
\[
 \begin{aligned}
 |\widetilde U(y)-\vartheta^\natural(\tau,x)|
 &\le\int_{\R}\varphi(z)|\widetilde U(y)-\widetilde U(x-z)|\dd z\\
 &\le\|\widetilde U_x\|_\infty
                   \int_{\R}\varphi(z)(|x-y|+|z|)\dd z
 \le C_\delta h(\tau)(1+|x-y|).
 \end{aligned}
\]
Smoothness of $\mathsf A$ and $\mathsf B$ gives the same bound for
the coefficient differences. Thus
\[
 \begin{aligned}
 &\left|\left.\bigl[\mathsf L(\vartheta,D)f_\vartheta\bigr](\tau,x)
                      \right|_{\vartheta=\vartheta^\natural(\tau,x)}\right|\\
 &\quad\le C_\delta h(\tau)^2\sup_\vartheta
         \int_{\R}(1+|z|)\bigl(|m_\vartheta(z)|+|m_\vartheta'(z)|\bigr)\dd z
 \le C_\delta h(\tau)^2 .
 \end{aligned}
\]
Bound the other terms of \eqref{eq:comparison-source-split} using
$1+x\le(1+|x-y|)(1+|y|)$ for $x>0$ and
\eqref{eq:comparison-source-budgets}.
\[
 \begin{aligned}
 &\bigl|\bigl[m_\vartheta*(n+g+\bar U'(\tau)\mathbf1_{\{x<0\}})\bigr](x)\bigr|\\
 &\quad\le\|m_\vartheta\|_1\|n\|_\infty
      +C_{\delta,N}(1+x)^{-N}\|(1+|y|)^N g\|_{L^1_y}
      +|\bar U'(\tau)|\int_x^\infty |m_\vartheta(z)|\dd z\\
 &\quad\le C_\delta h(\tau)^2
     +C_{\delta,N}(1+x)^{-N}
           \bigl(\|(1+|y|)^N g\|_{L^1_y}+|\bar U'(\tau)|\bigr)\\
 &\quad\le C_\delta h(\tau)^2+C_{\delta,N}M_I(1+x)^{-N}.
 \end{aligned}
\]
These bounds give \eqref{eq:comparison-actual-source}.

\subsubsection{The parameterized correction}

At fixed $\vartheta$, the symbols in \eqref{eq:comparison-frozen-generator}
and \eqref{eq:comparison-auxiliary-generator} are
$ikjb(\tau)I_3+\mathsf P(\vartheta,ik)$ and
$ikjb(\tau)I_9+\mathsf M(\vartheta,ik)$. The translated difference kernel
from Lemma~\ref{lemma:lift-kernel} has symbol
\[
 \widehat{\mathsf K}(\vartheta,\tau,r;k)
 =e^{ikj\int_r^\tau b(\rho)\dd\rho}
    \bigl[\mathsf L e^{(\tau-r)\mathsf P}
                        -e^{(\tau-r)\mathsf M}\mathsf L\bigr],
\]
where $\mathsf L=\mathsf L(\vartheta,k)$,
$\mathsf P=\mathsf P(\vartheta,ik)$ and $\mathsf M=\mathsf M(\vartheta,ik)$.
Define
\begin{equation}\label{eq:comparison-correction}
 \begin{aligned}
 \mathcal Q(\vartheta,\tau)
 &=\mathsf K(\vartheta,\tau,s)*\widetilde U(s)
   +\int_s^\tau\mathsf K(\vartheta,\tau,r)*E_\vartheta(r)\dd r,\\
 \mathcal Y(\vartheta,\tau)
 &=\mathsf L(\vartheta,D)\widetilde U(\tau)
                                      -\mathcal Q(\vartheta,\tau).
 \end{aligned}
\end{equation}
By integrability of the kernels and \eqref{eq:comparison-source-budgets},
the convolutions are defined on the bounded extension and the step
source. Since $\widehat{\mathsf K}$ vanishes at $k=0$, its primitive
$\partial_x^{-1}\mathsf K$ is Schwartz and satisfies
\[
 \begin{aligned}
 \mathsf K(\vartheta,\tau,s)*\widetilde U(s)
 &=\bigl(\partial_x^{-1}\mathsf K\bigr)(\vartheta,\tau,s)
                                     *\widetilde U_x(s),\\
 \bigl[\mathsf K(\vartheta,\tau,r)*\mathbf1_{\{x<0\}}\bigr](x)
 &=\int_x^\infty \mathsf K(\vartheta,\tau,r;y)\dd y
   =-\bigl(\partial_x^{-1}\mathsf K\bigr)(\vartheta,\tau,r;x).
 \end{aligned}
\]

By \eqref{eq:comparison-source-split}, the Fourier transform in $x$
satisfies
\[
 \begin{aligned}
 &\frac{\dd}{\dd r}\left[
    e^{ikj\int_r^\tau b(\rho)\dd\rho}
      \mathsf L e^{(\tau-r)\mathsf P}\widehat{\widetilde U}(r,k)\right]\\
 &\quad=e^{ikj\int_r^\tau b(\rho)\dd\rho}
      \mathsf L e^{(\tau-r)\mathsf P}
      \left[\partial_r\widehat{\widetilde U}
               -(ikjb(r)I_3+\mathsf P)\widehat{\widetilde U}\right](r,k)\\
 &\quad=e^{ikj\int_r^\tau b(\rho)\dd\rho}
                   \mathsf L e^{(\tau-r)\mathsf P}\widehat E_\vartheta(r,k).
 \end{aligned}
\]
The factors with $\mathsf L$ are smooth and compactly supported
in $k$. Hence the calculation is valid for the tempered Fourier
transforms of the extended fields. Integrating from $s$ to $\tau$
and subtracting the transform of $\mathcal Q$ gives
\[
 \begin{aligned}
 \widehat{\mathcal Y}(\vartheta,\tau,k)
 &=e^{ikj\int_s^\tau b(\rho)\dd\rho}
             e^{(\tau-s)\mathsf M}\mathsf L\widehat{\widetilde U}(s,k)\\
 &\quad+\int_s^\tau e^{ikj\int_r^\tau b(\rho)\dd\rho}
             e^{(\tau-r)\mathsf M}\mathsf L\widehat E_\vartheta(r,k)\dd r .
 \end{aligned}
\]
Differentiate in $\tau$ and invert the Fourier transform to get
\begin{equation}\label{eq:comparison-frozen-output}
 (\partial_\tau-\mathscr C_\vartheta)\mathcal Y
       =\mathsf L(\vartheta,D)E_\vartheta,
 \qquad \mathcal Y(\vartheta,s)=\mathsf L(\vartheta,D)\widetilde U(s).
\end{equation}

Set $J(\vartheta)=\mathsf L(\vartheta,0)$. At $k=0$,
\eqref{eq:lift-multiplier} gives three copies in each block of the
corresponding row of $T_0(\vartheta)$ from
\eqref{eq:affine-diagonalizer}. Thus $J_{ab}=0$ for parameter
derivatives $a,b$. For $k\ne0$, set
\[
 \widehat R_a(\vartheta,k)
 =\frac{\mathsf L_a(\vartheta,k)-\chi(k)J_a(\vartheta)}{ik},
 \qquad
 \widehat R_{ab}(\vartheta,k)=\frac{\mathsf L_{ab}(\vartheta,k)}{ik}.
\]
The numerators vanish at $k=0$ because $\chi(0)=1$ and $J_{ab}=0$.
After division by $ik$, the symbols are smooth and compactly supported,
with uniformly Schwartz kernels. Hence
\begin{equation}\label{eq:comparison-affine-split}
 \mathsf L_a=\chi(D)J_a+R_a(D)\partial_x,
 \qquad \mathsf L_{ab}=R_{ab}(D)\partial_x,
 \qquad J_{ab}=0.
\end{equation}
For $\mathcal Y_a^{\rm reg}=\mathcal Y_a-\chi(D)J_a\widetilde U$,
differentiate \eqref{eq:comparison-correction} to obtain
\[
 \begin{aligned}
 \mathcal Y_a^{\rm reg}&=R_a(D)\widetilde U_x-\mathcal Q_a,\qquad
 &\mathcal Y_{xa}&=\mathsf L_a(\vartheta,D)\widetilde U_x-\mathcal Q_{xa},\\
 \mathcal Y_{ab}&=R_{ab}(D)\widetilde U_x-\mathcal Q_{ab}.
 \end{aligned}
\]
The kernel bounds and \eqref{eq:comparison-extension-bounds} give
\[
 \begin{aligned}
 \|R_a(D)\widetilde U_x\|_\infty
       +\|\mathsf L_a(\vartheta,D)\widetilde U_x\|_\infty
 &\le C_\delta\|\widetilde U_x\|_\infty\le C_\delta h(\tau),\\
 \|R_{ab}(D)\widetilde U_x\|_\infty
 &\le C_\delta\|\widetilde U_x\|_2\le C_\delta M_I .
 \end{aligned}
\]
The first two instantaneous terms are therefore bounded by
$C_\delta H_I$ in $L^2(I;L^\infty)$.

For the initial term of \eqref{eq:comparison-correction}, use the
primitive identity and \eqref{eq:lift-kernel-bound}. For
$0\le n\le2$ and $|\beta|\le2$, these give
\[
 \begin{aligned}
 &\bigl\|\partial_x^n\partial_\vartheta^\beta
       [\mathsf K(\vartheta,\tau,s)*\widetilde U(s)]\bigr\|_\infty\\
 &\quad\le\|\partial_x^{n-1}\partial_\vartheta^\beta
                   \mathsf K(\vartheta,\tau,s)\|_2
                       \|\widetilde U_x(s)\|_2\\
 &\quad\le C_\delta\mathfrak b_s
    \begin{cases}
       \tau-s,&0\le\tau-s\le1,\\
       (\tau-s)^{-\frac{5}{4}-\frac{n-|\beta|}{2}},&\tau-s\ge1.
    \end{cases}
 \end{aligned}
\]
By \eqref{eq:comparison-extension-bounds},
$\|\widetilde U_x(s)\|_2\le C_\delta\mathfrak b_s$.
For $|\beta|\le1$, the $L^2(I;L^\infty)$ bound is at most
$C_\delta\mathfrak b_s
 \bigl(\int_0^1\rho^2\dd\rho+\int_1^\infty\rho^{-\frac{3}{2}}\dd\rho\bigr)^{\frac{1}{2}}
 \le C_\delta\mathfrak b_s$. For $n=0$ and $|\beta|=2$, the factor
$(\tau-s)^{-\frac{1}{4}}$ gives the $L^\infty$ bound. The pointwise
estimate gives the $L^\infty$ bound for $n=|\beta|=0$ as well.
In each time norm take the parameter supremum at each time.

The kernel estimates give
\begin{align}
 \|\mathcal Q\|_{L^\infty(I\times\R)}
 +\sum_{n=0}^2\|\partial_x^n\mathcal Q\|_{L^2(I;L^\infty)}
 &\le C_\delta(\mathfrak b_s+Q_I+S_I\beta_I+H_I^2),
                                      \label{eq:comparison-Q-budget}\\
 \|\mathcal Y_\vartheta^{\rm reg}\|_{L^2(I;L^\infty)}
 &\le C_\delta(\mathfrak b_s+H_I+H_I^2+L_IQ_I+S_I\beta_I),
                                      \label{eq:comparison-Yq-budget}\\
 \|\mathcal Y_{x\vartheta}\|_{L^2(I;L^\infty)}
 &\le C_\delta(\mathfrak b_s+H_I+H_I^2+Q_I+S_I\beta_I),
                                      \label{eq:comparison-Yxq-budget}\\
 \|\mathcal Y_{\vartheta\vartheta}\|_{L^\infty(I\times\R)}
 &\le C_\delta(M_I+\mathfrak b_s+H_I+H_I^2
                                   +L_I^{\frac{1}{2}}Q_I+S_I\beta_I).
                                      \label{eq:comparison-Yqq-budget}
\end{align}

In \eqref{eq:comparison-budget-table}, parameter derivatives in the
source integrals of \eqref{eq:comparison-correction} act on the kernel.
The first three rows use $L^2(I;L^\infty)$ and the last uses
$L^\infty(I\times\R)$.
\begin{equation}\label{eq:comparison-budget-table}
\begin{array}{c|cccc}
 &f_\vartheta&n&g&\bar U'\mathbf1_{x<0}\\ \hline
 \partial_x^n\mathcal Q\ (0\le n\le2)
       &S_I\beta_I&H_I^2&Q_I&Q_I\\
 \mathcal Q_\vartheta
       &S_I\beta_I&H_I^2&Q_I&L_IQ_I\\
 \mathcal Q_{x\vartheta}
       &S_I\beta_I&H_I^2&Q_I&Q_I\\
 \mathcal Q_{\vartheta\vartheta}
       &S_I\beta_I&H_I^2&Q_I&L_I^{\frac{1}{2}}Q_I
\end{array}
\end{equation}
Use the spatial bounds in \eqref{eq:lift-kernel-bound} for the translated
kernel at elapsed time $\rho=\tau-r$. For $0\le m\le2$ and
$|\beta|\le1$, Young's inequality in space and time gives
\begin{align*}
 &\left\|\int_s^\tau
       \partial_x^m\partial_\vartheta^\beta\mathsf K(\vartheta,\tau,r)
       *(f_\vartheta+n+g)(r)\dd r\right\|_{L^2(I;L^\infty)}\\*
 &\quad\le C_\delta\Biggl[
   \|f_\vartheta\|_{L^2(I;L^2)}\int_0^\infty
                                     \frac{\dd\rho}{(1+\rho)^{\frac{5}{4}}}
   +\|n\|_{L^1(I;L^\infty)}
                \left(\int_0^\infty\frac{\dd\rho}{(1+\rho)^2}\right)^{\frac{1}{2}}\\*
 &\qquad\qquad+\|g\|_{L^2(I;L^1)}\int_0^\infty
                                     \frac{\dd\rho}{(1+\rho)^{\frac{3}{2}}}\Biggr]\\
 &\quad\le C_\delta(S_I\beta_I+H_I^2+Q_I).
\end{align*}
The three spatial kernel norms are $L^2_x$, $L^1_x$ and
$L^\infty_x$. Use \eqref{eq:comparison-source-budgets} for the sources.
For $|\beta|=0,2$, the same spatial norms and Cauchy--Schwarz in
time give
\[
 \begin{aligned}
 &\left\|\int_s^\tau
       \partial_\vartheta^\beta\mathsf K(\vartheta,\tau,r)
       *(f_\vartheta+n+g)(r)\dd r\right\|_{L^\infty(I\times\R)}\\
 &\quad\le C_\delta\Biggl[
   \|f_\vartheta\|_{L^2(I;L^2)}
           \left(\int_0^\infty\frac{\dd\rho}{(1+\rho)^{\frac{3}{2}}}\right)^{\frac{1}{2}}
           +\|n\|_{L^1(I;L^\infty)}\\
 &\qquad\qquad+\|g\|_{L^2(I;L^1)}
           \left(\int_0^\infty\frac{\dd\rho}{(1+\rho)^2}\right)^{\frac{1}{2}}\Biggr]\\
 &\quad\le C_\delta(S_I\beta_I+H_I^2+Q_I).
 \end{aligned}
\]

For the left step use the primitive identity after
\eqref{eq:comparison-correction} and the bound in
\eqref{eq:comparison-source-budgets}. For $0\le m\le2$ and
$|\beta|\le1$, this gives
\[
 \begin{aligned}
 &\left\|\int_s^\tau
       \partial_x^m\partial_\vartheta^\beta\mathsf K(\vartheta,\tau,r)
       *[\bar U'(r)\mathbf1_{\{x<0\}}]\dd r\right\|_{L^2(I;L^\infty)}\\
 &\quad=\left\|\int_s^\tau
       \bigl(\partial_x^{m-1}\partial_\vartheta^\beta
                      \mathsf K\bigr)(\vartheta,\tau,r)\bar U'(r)\dd r
                              \right\|_{L^2(I;L^\infty)}\\
 &\quad\le C_\delta\|\bar U'\|_{L^2(I)}
        \int_0^{t-s}\frac{\dd\rho}{(1+\rho)^{\frac{3}{2}+\frac{m-|\beta|}{2}}}
 \le C_\delta Q_I
 \begin{cases}
  L_I,&m=0,\ |\beta|=1,\\
  1,&\text{otherwise}.
 \end{cases}
 \end{aligned}
\]
For $|\beta|=0,2$, the analogous supremum estimate is
\[
 \begin{aligned}
 &\left\|\int_s^\tau
       \partial_\vartheta^\beta\mathsf K(\vartheta,\tau,r)
       *[\bar U'(r)\mathbf1_{\{x<0\}}]\dd r\right\|_{L^\infty(I\times\R)}\\
 &\quad\le C_\delta\|\bar U'\|_{L^2(I)}
         \left(\int_0^{t-s}\frac{\dd\rho}{(1+\rho)^{3-|\beta|}}\right)^{\frac{1}{2}}
 \le C_\delta Q_I
 \begin{cases}
  1,&|\beta|=0,\\
  L_I^{\frac{1}{2}},&|\beta|=2.
 \end{cases}
 \end{aligned}
\]
The factors of $L_I$ follow from
$\int_0^{t-s}\frac{\dd\rho}{1+\rho}=\log(1+t-s)\le L_I$.
With $|\beta|=0$, the bounds give the source terms in the supremum
estimate of \eqref{eq:comparison-Q-budget}.

Since $\widetilde U$ and $b$ are independent of the formal parameter,
\eqref{eq:comparison-source-split} gives
\[
 \partial_a E_\vartheta
   =j\mathsf A_a\widetilde U_x-\mathsf B_a\widetilde U_{xx},
 \qquad
 \partial_{ab}E_\vartheta
   =j\mathsf A_{ab}\widetilde U_x-\mathsf B_{ab}\widetilde U_{xx}.
\]
At fixed $\vartheta$ the matrices are constant in space, so
integration by parts gives
\begin{equation}\label{eq:comparison-density-transfer}
 \mathsf K*(\mathsf B_a\widetilde U_{xx})
        =\mathsf K_x*(\mathsf B_a\widetilde U_x),
 \qquad
 \mathsf K_b*(\mathsf B_a\widetilde U_{xx})
        =\mathsf K_{xb}*(\mathsf B_a\widetilde U_x).
\end{equation}
The first identity also holds with $\mathsf B_{ab}$ in place of
$\mathsf B_a$. Spatial derivatives act on the kernels. Evaluate
the kernel at $(\vartheta,\tau,r)$ and the source at $r$ below.

By \eqref{eq:lift-kernel-bound} in $L^1_x$ and
\eqref{eq:comparison-extension-bounds}, for $m=0,1$ we obtain
\[
 \begin{aligned}
 &\left\|\int_s^\tau\partial_x^m\mathsf K*
                   \partial_a E_\vartheta(r)\dd r\right\|_{L^2(I;L^\infty)}\\
 &\quad\le C_\delta\|\widetilde U_x\|_{L^2(I;L^\infty)}
      \int_0^\infty\left[
       \frac{1}{(1+\rho)^{\frac{3+m}{2}}}
       +\frac{1}{(1+\rho)^{\frac{4+m}{2}}}\right]\dd\rho
 \le C_\delta H_I .
 \end{aligned}
\]
For two parameter derivatives, the terms differentiated in the
source satisfy
\[
 \begin{aligned}
 &\left\|\int_s^\tau
       \bigl[\mathsf K_a*\partial_b E_\vartheta(r)
       +\mathsf K_b*\partial_a E_\vartheta(r)
       +\mathsf K*\partial_{ab}E_\vartheta(r)\bigr]\dd r
                                 \right\|_{L^\infty(I\times\R)}\\
 &\quad\le C_\delta\|\widetilde U_x\|_{L^2(I;L^\infty)}
       \left[\left(\int_0^\infty\frac{\dd\rho}{(1+\rho)^2}\right)^{\frac{1}{2}}
            +\left(\int_0^\infty\frac{\dd\rho}{(1+\rho)^3}\right)^{\frac{1}{2}}\right]
 \le C_\delta H_I .
 \end{aligned}
\]
The time Cauchy--Schwarz estimate uses the $L^1_x$ bound
$C_\delta(1+\rho)^{-1}$ for $\mathsf K,\mathsf K_a,\mathsf K_b$
and $C_\delta(1+\rho)^{-\frac{3}{2}}$ for their spatial derivatives.
Combine these bounds, \eqref{eq:comparison-budget-table} and the
initial and instantaneous estimates to obtain
\eqref{eq:comparison-Q-budget}--\eqref{eq:comparison-Yqq-budget}.

For reconstruction, use the initial bound and the pointwise
bounds in \eqref{eq:comparison-source-budgets}. The $L^2_x$ and
$L^\infty_x$ bounds for $\mathsf K$ and the $L^\infty_x$ bound for
its primitive yield
\begin{equation}\label{eq:comparison-endpoint-Q}
 \begin{aligned}
 \|\mathcal Q\|_{L^\infty(I\times\R)}
 &\le C_\delta\mathfrak b_s+
       C_\delta(S_I M_I+M_I^2)\int_0^\infty
                                  \frac{\dd\rho}{(1+\rho)^{\frac{7}{4}}}\\
 &\quad+C_\delta M_I\int_0^\infty
        \left[\frac{1}{(1+\rho)^2}+
                        \frac{1}{(1+\rho)^{\frac{3}{2}}}\right]\dd\rho\\
 &\le C_\delta\bigl[\mathfrak b_s+M_I(1+S_I+M_I)\bigr]
 \le C_\delta(\mathfrak b_s+M_I).
 \end{aligned}
\end{equation}
Here $S_I$ and $M_I$ are uniformly small.

\subsubsection{Evaluation at the natural parameter}

Set $Y(\tau,x)=\mathcal Y(\vartheta^\natural(\tau,x),\tau,x)$.
At fixed $x$, the right-hand side of the $\mathcal Y_\tau$ equation
in \eqref{eq:comparison-frozen-output} is continuous in $\vartheta$.
Its supremum over the parameter ball is time integrable by
\eqref{eq:comparison-Q-budget}, \eqref{eq:comparison-source-budgets}
and the smooth kernels of $\mathsf L$. Equation
\eqref{eq:comparison-parameter-derivatives} gives absolute continuity
of $\vartheta^\natural(\cdot,x)$, so the chain rule applies.
\[
 \begin{aligned}
 Y_\tau&=\mathcal Y_\tau+\mathcal Y_\vartheta\vartheta^\natural_t,
 &Y_x&=\mathcal Y_x+\mathcal Y_\vartheta\vartheta^\natural_x,\\
 Y_{xx}&=\mathcal Y_{xx}+2\mathcal Y_{x\vartheta}\vartheta^\natural_x
             +\mathcal Y_\vartheta\vartheta^\natural_{xx}
             +\mathcal Y_{\vartheta\vartheta}
                           [\vartheta^\natural_x,\vartheta^\natural_x].
 \end{aligned}
\]
Evaluate the fixed-parameter derivatives at
$\vartheta=\vartheta^\natural(\tau,x)$. The switching term in block
$e$ is then $\mathsf Q_e(\vartheta^\natural)Y_e$. Substitution in
\eqref{eq:comparison-auxiliary-generator} and
\eqref{eq:comparison-frozen-output} gives, for component $a$ of
block $e$,
\begin{equation}\label{eq:comparison-chain-rule}
 \begin{aligned}
 \bigl[(\partial_\tau-\mathscr C_{\vartheta^\natural})Y\bigr]_{e,a}
 &=\left.[\mathsf L(\vartheta,D)E_\vartheta]_{e,a}
                                             \right|_{\vartheta^\natural}
 +(\mathcal Y_\vartheta)_{e,a}\cdot
       (\vartheta^\natural_t-v_{e,a}\vartheta^\natural_x
                                  -D_e\vartheta^\natural_{xx})\\
 &\quad-2D_e(\mathcal Y_{x\vartheta})_{e,a}\cdot\vartheta^\natural_x
 -D_e(\mathcal Y_{\vartheta\vartheta})_{e,a}:
                    (\vartheta^\natural_x\otimes\vartheta^\natural_x).
 \end{aligned}
\end{equation}
Here $D_e$ and $v_{e,a}$ are evaluated at the natural parameter,
with $jb(\tau)$ included in $v_{e,a}$.

The convolution in \eqref{eq:comparison-natural-parameter} and
\eqref{eq:comparison-extension-bounds} give
$\|\vartheta^\natural_t\|_\infty+\|\vartheta^\natural_x\|_\infty+
\|\vartheta^\natural_{xx}\|_\infty\le C_\delta h(\tau)$, using
$\vartheta^\natural_{xx}=\varphi'*\widetilde U_x$.
In the first parameter term of \eqref{eq:comparison-chain-rule},
use $\mathcal Y_\vartheta^{\rm reg}$ from
\eqref{eq:comparison-affine-split}. By \eqref{eq:comparison-actual-source},
the sum of the resulting terms on $x>0$ is bounded by
$C_{\delta,N}M_I(1+x)^{-N}$ plus the spatially constant integrand below.
Cauchy--Schwarz and
\eqref{eq:comparison-Yq-budget}--\eqref{eq:comparison-Yqq-budget} give
\begin{equation}\label{eq:comparison-regular-chain}
 \begin{aligned}
 &C_\delta\int_I\Bigl[
    h^2\bigl(1+\|\mathcal Y_{\vartheta\vartheta}\|_\infty\bigr)
   +h\bigl(\|\mathcal Y_\vartheta^{\rm reg}\|_\infty
                          +\|\mathcal Y_{x\vartheta}\|_\infty\bigr)\Bigr]\dd\tau\\
 &\quad\le C_\delta\Bigl[
    H_I^2\bigl(1+\|\mathcal Y_{\vartheta\vartheta}\|_{L^\infty(I\times\R)}\bigr)
    +H_I\bigl(\|\mathcal Y_\vartheta^{\rm reg}\|_{L^2(I;L^\infty)}
                      +\|\mathcal Y_{x\vartheta}\|_{L^2(I;L^\infty)}\bigr)\Bigr]\\
 &\quad\le C_\delta\Bigl[
    H_I^2+H_I\bigl(\mathfrak b_s+H_I+H_I^2+L_IQ_I+S_I\beta_I\bigr)\\
 &\qquad\qquad\quad
    +H_I^2\bigl(M_I+\mathfrak b_s+H_I+H_I^2
                              +L_I^{\frac{1}{2}}Q_I+S_I\beta_I\bigr)\Bigr]\\
 &\quad\le C_\delta\bigl(H_I^2+H_I\mathfrak b_s+L_IH_IQ_I
                                      +S_I\beta_IH_I\bigr)=:\mathcal E_I .
 \end{aligned}
\end{equation}
The last line uses $L_I\ge1$ and the uniform bounds from
\eqref{eq:window-inputs} and the global $H^2$ estimate.

By \eqref{eq:affine-diagonalizer} and $J(\vartheta)=\mathsf L(\vartheta,0)$,
the rows of $J$ satisfy
\[
 (J(\vartheta)\widetilde U)_{0,a}=\widetilde U_0,\qquad
 (J(\vartheta)\widetilde U)_{e,a}
       =\widetilde U_e+k_0\vartheta_{-e}\widetilde U_0,
 \qquad e\in\{+1,-1\}.
\]
Differentiating at fixed $\widetilde U$, the contribution of
$\mathcal Y_\vartheta-\mathcal Y_\vartheta^{\rm reg}
=\chi(D)J_\vartheta\widetilde U$ to
\eqref{eq:comparison-chain-rule} is $R^0_{0,a}=0$ and
\begin{equation}\label{eq:comparison-zero-chain}
 \begin{aligned}
 R^0_{e,a}
 &=\bigl[\chi(D)J_\vartheta\widetilde U\bigr]_{e,a}
       \cdot(\vartheta^\natural_t-v_{e,a}\vartheta^\natural_x
                                  -D_e\vartheta^\natural_{xx})\\
 &=k_0[\chi(D)\widetilde U_0]
       \bigl(\vartheta^\natural_{-e,t}
       -v_{e,a}\vartheta^\natural_{-e,x}
       -D_e\vartheta^\natural_{-e,xx}\bigr),
 \qquad e\in\{+1,-1\}.
 \end{aligned}
\end{equation}

Take $P^0$ from Corollary~\ref{cor:zero-chain} with zero initial
value. Its source uses the physical half-line inputs of the two
factors in \eqref{eq:comparison-zero-chain}, with restriction before
convolution. Extend the product on $x>0$ by zero. The convolution
bounds in Corollary~\ref{cor:zero-chain} place the source in
$L^1(I;L^\infty)$. At fixed $\tau$, the Duhamel integral over times
away from $\tau$ is continuous in $x$. Unit row mass bounds the
remaining integral uniformly in $x$ by the $L^1_\tau L^\infty_x$
source norm on the terminal subinterval. This gives the trace at
$x=0$. Corollary~\ref{cor:zero-chain} now gives
\begin{equation}\label{eq:comparison-zero-corrector}
 \|P^0\|_{L^\infty(I\times\R_+)}
      +\|P^0(\cdot,0)\|_{L^\infty(I)}
 \le C_\delta\eps K_I .
\end{equation}
By the bounds used in \eqref{eq:comparison-regular-chain} and
Corollary~\ref{cor:zero-chain}, for $x>0$ and $N>1$,
\[
 \begin{aligned}
 &\left|\bigl[(\partial_\tau-\mathscr C_{\vartheta^\natural})
                                      (Y-P^0)\bigr]_{e,a}\right|\\
 &\quad\le
 \left|\bigl[(\partial_\tau-\mathscr C_{\vartheta^\natural})Y\bigr]_{e,a}
                                                 -R^0_{e,a}\right|
 +\left|R^0_{e,a}
       -\bigl[(\partial_\tau-\mathscr C_{\vartheta^\natural})
                                                P^0\bigr]_{e,a}\right|\\
 &\quad\le C_\delta\Bigl[
      h^2\bigl(1+\|\mathcal Y_{\vartheta\vartheta}\|_\infty\bigr)
     +h\bigl(\|\mathcal Y_\vartheta^{\rm reg}\|_\infty
                         +\|\mathcal Y_{x\vartheta}\|_\infty\bigr)\Bigr]
            +C_{\delta,N}(1+S_I)M_I(1+x)^{-N}.
 \end{aligned}
\]
By \eqref{eq:comparison-regular-chain}, the spatially constant term
has time integral at most $\mathcal E_I$. Since
$S_I\le C_\delta\eps$, the decaying term is bounded by
$C_{\delta,N}M_I(1+x)^{-N}$.

\subsubsection{Reconstruction and comparison at the boundary}

Since $J(\vartheta)=\mathsf L(\vartheta,0)$, the quotient
$\frac{\mathsf L(\vartheta,k)-J(\vartheta)}{ik}$ extends smoothly
to $k=0$ by the formula
$\frac{1}{i}\int_0^1(\partial_k\mathsf L)(\vartheta,rk)\dd r$.
For $|k|>k_*$ the quotient is $-\frac{J(\vartheta)}{ik}$ by
\eqref{eq:lift-multiplier}. Hence
\[
 \begin{aligned}
 \sup_\vartheta\left\|
       \frac{\mathsf L(\vartheta,k)-J(\vartheta)}{ik}\right\|_{L^2_k}^2
 &\le C_\delta\left(\int_{|k|\le k_*}\dd k
                         +\int_{|k|>k_*}\frac{\dd k}{k^2}\right)
 \le C_\delta .
 \end{aligned}
\]
For the bounded extension $\widetilde U$, the identity
\[
 \widehat{(\mathsf L(\vartheta,D)-J(\vartheta))\widetilde U}
 =\frac{\mathsf L(\vartheta,k)-J(\vartheta)}{ik}\,
                             \widehat{\widetilde U_x}
\]
holds in tempered distributions. Cauchy--Schwarz gives an $L^1_k$
bound for the right-hand side. Since $\widetilde U=U$ on $\R_+$,
Fourier inversion, Plancherel's identity,
\eqref{eq:comparison-extension-bounds} and
\eqref{eq:comparison-endpoint-Q} give
\begin{equation}\label{eq:comparison-reconstruction}
 \begin{aligned}
 &\|Y-J(\vartheta^\natural)U\|_{L^\infty(I\times\R_+)}\\
 &\quad\le C\sup_\vartheta
       \left\|\frac{\mathsf L(\vartheta,k)-J(\vartheta)}{ik}\right\|_{L^2_k}
             \|\widetilde U_x\|_{L^\infty(I;L^2)}
             +\left\|\left.\mathcal Q\right|_{\vartheta^\natural}
                         \right\|_{L^\infty(I\times\R_+)}\\
 &\quad\le C_\delta(\mathfrak b_s+M_I).
 \end{aligned}
\end{equation}
In an acoustic block $e=+1,-1$, the entries of
$(J(\vartheta^\natural)-J(0))U$ are $k_0\vartheta^\natural_{-e}U_0$.
The contact entries vanish. By
\eqref{eq:comparison-natural-parameter} and
\eqref{eq:comparison-extension-bounds},
\begin{equation}\label{eq:comparison-quadratic-reconstruction}
 \begin{aligned}
 \|J(\vartheta^\natural)U-J(0)U\|_{L^\infty(I\times\R_+)}
 &\le C_\delta\|\vartheta^\natural\|_{L^\infty(I\times\R)}
                                 \|U_0\|_{L^\infty(I\times\R_+)}
 \le C_\delta S_I^2 .
 \end{aligned}
\end{equation}
Both bounds hold at $x=0$ by spatial continuity.
Since $J(0)$ repeats each input component three times,
$J(0)^{\mathsf T}J(0)=3I_3$ and $\frac{1}{3}J(0)^{\mathsf T}$
is a left inverse. For $\tau=s$, \eqref{eq:comparison-frozen-output}
and $P^0(s)=0$ give
\[
 \|(Y-P^0)(s)\|_{L^\infty(\R_+)}
 \le \sup_\vartheta\|m_\vartheta\|_1\|\widetilde U(s)\|_\infty
 \le C_\delta S(s).
\]

In each block the operator is
\[
 \mathscr C=D(\tau,x)\partial_x^2 I_3
       +\operatorname{diag}(v_a(\tau,x))\partial_x+\mathsf Q(\tau,x).
\]
By Lemma~\ref{lemma:lift}, the system is uniformly parabolic,
$D\le D^*$, and $\mathsf Q$ has nonnegative off-diagonal entries
with $\mathsf Q\mathbf1=0$. In an incoming block
$v_a\ge g_\delta>0$ for every component, and in an outgoing block
$v_a\le-g_\delta$ for every component.
For $\gamma_\delta=\frac{g_\delta}{2D^*}$, the incoming barrier satisfies
\[
 -\bigl[\mathscr C(e^{-\gamma_\delta x}\mathbf1)\bigr]_a
 =\gamma_\delta(v_a-D\gamma_\delta)e^{-\gamma_\delta x}
 \ge \frac{g_\delta\gamma_\delta}{2}e^{-\gamma_\delta x},
 \qquad \mathscr C\mathbf1=0 .
\]
The constant vector $\mathbf1$ serves as the outgoing boundary barrier.

Fix $N>1$ and take $R_0\ge\max\{1,\frac{2ND^*}{g_\delta}\}$. For
$w(x)=(1+\frac{x}{R_0})^{-N}$ put
\[
 p(x)=\begin{cases}
       \displaystyle\int_x^\infty w(y)\,\dd y,& v_a\ge g_\delta,\\
       \displaystyle 1+\int_0^x w(y)\,\dd y,& v_a\le-g_\delta .
      \end{cases}
 \qquad 0\le p(x)\le 1+\frac{R_0}{N-1}.
\]
Since $w_x=-\frac{N}{R_0+x}w$, in both cases we have
\[
 \begin{aligned}
 -\bigl[\mathscr C(p\mathbf1)\bigr]_a
 &=\begin{cases}
     \bigl(v_a-\frac{ND}{R_0+x}\bigr)w(x),& v_a\ge g_\delta,\\[2pt]
     \bigl(\frac{ND}{R_0+x}-v_a\bigr)w(x),& v_a\le-g_\delta
   \end{cases}\\
 &\ge \frac{g_\delta}{2}w(x)\ge \frac{g_\delta}{2}(1+x)^{-N}.
 \end{aligned}
\]

Let $m_0$ and $m_1$ bound $|V(s,x)|$ and $|V(\tau,0)|$,
respectively, for a bounded block solution $V$. Suppose
$\bigl|[(\partial_\tau-\mathscr C)V]_a\bigr|
\le f(\tau)+m_2(1+x)^{-N}$, where $m_2\ge0$, $f\ge0$ and
$\int_I f(\tau)\,\dd\tau\le\mathcal E_I$. The vector
\[
 \left(m_0+\int_s^\tau f(r)\,\dd r
                 +\frac{2m_2}{g_\delta}p(x)\right)\mathbf1
       +m_1\begin{cases}
          e^{-\gamma_\delta x}\mathbf1,& v_a\ge g_\delta,\\
          \mathbf1,& v_a\le-g_\delta
        \end{cases}
\]
is an upper bound for the initial and boundary values of $\pm V$.
Its image under $\partial_\tau-\mathscr C$ is bounded below by
$(f+m_2w)\mathbf1$.

For comparison on a finite cylinder, add
$\eta e^{C(\tau-s)}(1+x^2)\mathbf1$, with $\eta>0$ and
$C\ge2D^*+\sup_{\tau,x,a}|v_a|$. Then
\[
 \begin{aligned}
 \bigl[(\partial_\tau-\mathscr C)
         (\eta e^{C(\tau-s)}(1+x^2)\mathbf1)\bigr]_a
 &=\eta e^{C(\tau-s)}\bigl[C(1+x^2)-2D-2v_a x\bigr]\\
 &\ge \eta e^{C(\tau-s)}
          \bigl(C-2D^*-\sup_{\tau,x,a}|v_a|\bigr)(1+x^2)\ge0 .
 \end{aligned}
\]
For fixed $\eta>0$, the added term bounds $|V|$ at every sufficiently
large right endpoint. On the finite cylinder use the cooperative
maximum principle for both signs and the positive-part test for
the integrable time source. Let the right endpoint tend to infinity
and then let $\eta$ tend to zero.
\begin{equation}\label{eq:comparison-halfline}
 |V(\tau,x)|\le C_\delta\left[
    m_0+\mathcal E_I+m_2+
       m_1\begin{cases}
          e^{-\gamma_\delta x},&\text{incoming},\\
          1,&\text{outgoing}.
       \end{cases}\right]
\end{equation}
The constant is independent of $I$.

Set $A_I=\sup_{\tau\in I}|\alpha(\tau)|$. For $R_1\ge0$, use
\eqref{eq:corrected-end-field}, bounded derivatives of $N_{q_j}$ and
the fundamental theorem of calculus from $y\in[-1,1]$ to $jR_1$.
This gives
\begin{equation}\label{eq:comparison-fixed-trace}
 \begin{aligned}
 |U^j(\tau,jR_1)-N_{q_j}(q_{\alpha(\tau),j})|
 &\le C_\delta|r(\tau,jR_1)|\\
 &\le C_\delta\left[\frac{1}{2}\int_{-1}^{1}|r(\tau,y)|\dd y
                    +\sqrt{1+R_1}\,\|r_\xi(\tau)\|_2\right]\\
 &\le C_\delta\left[m_\omega(r(\tau))+
                    \sqrt{1+R_1}\,\|r_\xi(\tau)\|_2\right]
 \le C_{\delta,R_1}\sqrt{\mathcal T(\tau)} .
 \end{aligned}
\end{equation}
We used the positive lower bound of $\omega$ on $[-1,1]$,
\eqref{eq:st-spatial-estimates} and \eqref{eq:time-energy-equivalence}.
Continuity of $\mathcal T$ gives $\sqrt{\mathcal T(\tau)}\le M_I$
on $I$. By \eqref{eq:st-family-bounds},
$|N_{q_j}(q_{\alpha,j})|\le C_\delta|\alpha|$. At $x=0$,
\eqref{eq:comparison-reconstruction},
\eqref{eq:comparison-quadratic-reconstruction},
\eqref{eq:comparison-zero-corrector} and
\eqref{eq:comparison-fixed-trace} give
\[
 \begin{aligned}
 \|(Y-P^0)(\cdot,0)\|_{L^\infty(I)}
 &\le C_\delta\sup_{\tau\in I}|U^j(\tau,0)|
                 +C_\delta(\mathfrak b_s+M_I+S_I^2+\eps K_I)\\
 &\le C_\delta(A_I+\mathfrak b_s+M_I+S_I^2+\eps K_I).
 \end{aligned}
\]

By Lemma~\ref{lemma:incoming}, for all sufficiently small $\alpha$,
$|\alpha|\le C_\delta\sum_{(j,e)\ {\rm incoming}}
                        |[N_{q_j}(q_{\alpha,j})]_e|$.
In \eqref{eq:comparison-halfline} use $V=Y-P^0$,
$m_0=C_\delta S(s)$ and $m_2=C_{\delta,N}M_I$, by the initial and
source bounds. Let $[(Y-P^0)(\tau,R_1)]_{j,e}$ be the $e$-block
on side $j$. Each component of $U^j$ is repeated by $J(0)$.
The reconstruction and trace bounds give
\begin{equation}\label{eq:comparison-alpha-before-absorption}
 \begin{aligned}
 A_I
 &\le C_\delta\sum_{(j,e)\ {\rm incoming}}
          \sup_{\tau\in I}|[(Y-P^0)(\tau,R_1)]_{j,e}|
          +C_{\delta,R_1}M_I+C_\delta(\mathfrak b_s+S_I^2+\eps K_I)\\
 &\le C_\delta\bigl[S(s)+\mathcal E_I+M_I+
      e^{-\gamma_\delta R_1}(A_I+\mathfrak b_s+M_I+S_I^2+\eps K_I)\bigr]\\
 &\quad+C_{\delta,R_1}M_I+C_\delta(\mathfrak b_s+S_I^2+\eps K_I)\\
 &\le C_{\delta,R_1}\bigl[S(s)+\mathfrak b_s+M_I+
                   \mathcal E_I+S_I^2+\eps K_I\bigr]
                  +C_\delta e^{-\gamma_\delta R_1}A_I .
 \end{aligned}
\end{equation}
The constant in front of $e^{-\gamma_\delta R_1}A_I$ is independent
of $R_1$. Fix $R_1$ with $C_\delta e^{-\gamma_\delta R_1}\le\frac{1}{2}$
and absorb the last term of \eqref{eq:comparison-alpha-before-absorption}.

By \eqref{eq:corrected-end-field}, \eqref{eq:st-parameter-tails} and
the local inverse of $N_{q_j}$, we have
$|Z(\tau,\xi)|\le C_\delta(|U^j(\tau,\xi)|+
|\alpha(\tau)|e^{-2\eta_0|\xi|})$ for $j\xi>0$.
Write $[Y-P^0]_j$ for the comparison vector on side $j$.
The factor multiplying $m_1$ in \eqref{eq:comparison-halfline}
is at most one in every block. Using the left inverse of $J(0)$
with \eqref{eq:comparison-zero-corrector}--\eqref{eq:comparison-quadratic-reconstruction},
the boundary bound at $x=0$, and the absorbed bound in
\eqref{eq:comparison-alpha-before-absorption}, we obtain
\begin{equation}\label{eq:comparison-before-absorption}
 \begin{aligned}
 S_I
 &\le C_\delta\left[\max_{j=\pm1}
                  \|U^j\|_{L^\infty(I\times\R_j)}+A_I\right]\\
 &\le C_\delta\left[\max_{j=\pm1}
           \|[Y-P^0]_j\|_{L^\infty(I\times\R_+)}
           +A_I+\mathfrak b_s+M_I+S_I^2+\eps K_I\right]\\
 &\le C_\delta\left[S(s)+\mathfrak b_s+M_I+\mathcal E_I+
                     S_I^2+\eps K_I+A_I\right]\\
 &\le C_\delta\left[S(s)+\mathfrak b_s+M_I+\mathcal E_I+
                     S_I^2+\eps K_I\right]\\
 &\le C_\delta\left[S(s)+\mathfrak b_s+\Theta_I+H_I^2+
          H_I\mathfrak b_s+L_IH_IQ_I+S_I\beta_IH_I+S_I^2+\eps K_I\right].
 \end{aligned}
\end{equation}
Use \eqref{eq:comparison-regular-chain} and $M_I\le\Theta_I$ in the
last line.

By \eqref{eq:global-H2} and \eqref{eq:window-inputs},
$S_I+\beta_I+H_I\le C_\delta\eps$. Substituting
$K_I=\mathfrak b_s+Q_I+S_I\beta_I$ into
\eqref{eq:comparison-before-absorption} and using $Q_I\le\Theta_I$ gives
\[
 \begin{aligned}
 S_I
 &\le C_\delta\left[S(s)+\mathfrak b_s+\Theta_I+H_I^2+
      H_I\mathfrak b_s+L_IH_IQ_I+\eps(\mathfrak b_s+Q_I)\right]\\
 &\quad+C_\delta S_I(\beta_IH_I+S_I+\eps\beta_I)\\
 &\le C_\delta\left[S(s)+\mathfrak b_s+\Theta_I+H_I^2+
      H_I\mathfrak b_s+L_IH_IQ_I\right]+C_\delta\eps S_I.
 \end{aligned}
\]
Take $C_\delta\eps\le\frac{1}{2}$ and absorb the final term to obtain
\eqref{eq:comparison-final}.

For a strong solution, $U^j$ is continuous in $H^2(\R_j)$ and
$U^j_t$ in $L^2(\R_j)$. The extension adds a bounded time-dependent
constant and terms supported in a fixed strip. Its convolution with
$\varphi$ is smooth in space and continuously differentiable in
time. The multipliers in \eqref{eq:comparison-correction} have fixed
frequency support. The bounds \eqref{eq:comparison-source-budgets}
and \eqref{eq:comparison-Q-budget}--\eqref{eq:comparison-Yqq-budget}
justify differentiation of the Duhamel formulas and the chain rule.
The flux regularity in \eqref{eq:time-potential-input} allows the
integrated cumulative-kernel estimate for the time differentiated
conserved variables. These are the regularity bounds of
Propositions~\ref{prop:global}, \ref{prop:time} and \ref{prop:spatial}.

\section{Proof of Theorem~\ref{thm:main}}\label{sec:completion}

We apply \eqref{eq:window-inputs} and \eqref{eq:comparison-final}
on subintervals of $[\frac{t}{2},t]$.

\begin{lemma}[Decay from finite-window estimates]\label{lemma:window-decay}
Let $S\ge0$ be continuous on $[0,\infty)$, and let $B,\mathcal T\ge0$
be measurable. Assume that for some $A\ge1$ and $\eps\ge0$,
\begin{equation}\label{eq:abstract-basic}
 \begin{aligned}
 \sup_{t\ge0}S(t)&\le A\eps,\qquad
          \int_0^\infty B(\tau)\dd\tau\le A\eps^2,\\
 \mathcal T(t)&\le AB(t),\qquad
          S(t)^4\le A\eps^2B(t)\quad\text{for a.e. }t\ge0 .
 \end{aligned}
\end{equation}
For $I=[s,t]$ use $S_I,\beta_I,\mathfrak b_s$ from
\eqref{eq:time-window-notation} and $M_I,L_I$ from
\eqref{eq:window-notation}. Suppose that for each $t>0$ and almost
every $s\in[0,t)$ there are $Q_I,H_I\ge0$ such that, with
$\Theta_I=M_I+Q_I$,
\begin{align}
 \Theta_I&\le A(\mathfrak b_s+\beta_IS_I),\notag\\
 H_I^2&\le A\beta_IQ_I,\notag\\
 S_I&\le A\bigl[S(s)+\mathfrak b_s+\Theta_I
                   +H_I^2+H_I\mathfrak b_s+L_IH_IQ_I\bigr].
 \label{eq:abstract-window}
\end{align}
If $\eps\le\eps_*(A)$, then
\[
 \sup_{t\ge0}(1+t)^{\frac{1}{4}}S(t)\le C(A)\eps,\qquad
 (1+t)S(t)^4\le C(A)\eps^2\int_{\frac{t}{2}}^t B(\tau)\dd\tau,
 \quad t\ge1 .
\]
\end{lemma}
\begin{proof}
The bound $\sup_{t\ge0}S(t)\le A\eps$ covers $\eps=0$.
Assume $0<\eps\le1$. For $T\ge1$, suppose that
\[
 \mathfrak M_T=\sup_{0\le\tau\le T}(1+\tau)^{\frac{1}{4}}S(\tau)\le1 .
\]
Fix $1\le t\le T$ and set
$\beta^2=\int_{\frac{t}{2}}^t B(\tau)\dd\tau$.
If $\beta=0$, \eqref{eq:abstract-basic} and continuity give
$S(t)=0$. For $\beta>0$, \eqref{eq:abstract-basic} gives
\[
 \frac{4}{t}\int_{\frac{t}{2}}^{\frac{3t}{4}}
       \left[\mathcal T(\tau)+\frac{S(\tau)^4}{\eps^2}\right]\dd\tau
 \le \frac{8A}{t}\int_{\frac{t}{2}}^t B(\tau)\dd\tau
 =\frac{8A\beta^2}{t}.
\]
The integrand is at most $\frac{16A\beta^2}{t}$
on a set of measure at least $\frac{t}{8}$. Choose $s$ from this set
where \eqref{eq:abstract-window} holds. Then
\begin{equation}\label{eq:common-good-time}
 \mathfrak b_s^2+\frac{S(s)^4}{\eps^2}\le\frac{16A\beta^2}{t},
 \qquad
 \mathfrak b_s\le C\frac{\beta}{\sqrt t},\qquad
 S(s)\le C\frac{\sqrt{\eps\beta}}{t^{\frac{1}{4}}}.
\end{equation}
In this proof, $C$ depends only on $A$. Since $s\ge\frac{t}{2}$,
the bound on $\mathfrak M_T$ gives
$S_I\le (1+s)^{-\frac{1}{4}}\le C t^{-\frac{1}{4}}$ on $I=[s,t]$.
With $\beta_I\le\beta$, \eqref{eq:abstract-window} gives
\[
 \begin{aligned}
 Q_I\le\Theta_I
 &\le A(\mathfrak b_s+\beta_I S_I)
 \le C\beta\left(\frac{1}{\sqrt t}+\frac{1}{t^{\frac{1}{4}}}\right)
 \le C\frac{\beta}{t^{\frac{1}{4}}},\\
 H_I^2&\le A\beta_I Q_I\le C\frac{\beta^2}{t^{\frac{1}{4}}}.
 \end{aligned}
\]
The bounds for $H_I,Q_I$, the inequality $L_I\le 1+\log(2+t)$
and \eqref{eq:common-good-time} give
\[
 \begin{aligned}
 H_I^2+H_I\mathfrak b_s+L_IH_IQ_I
 &\le C\beta^2\left[\frac{1}{t^{\frac{1}{4}}}
         +\frac{1}{t^{\frac{5}{8}}}
         +\frac{1+\log(2+t)}{t^{\frac{3}{8}}}\right]
 \le C\frac{\beta^2}{t^{\frac{1}{4}}},
 \end{aligned}
\]
where $1+\log(2+t)\le C t^{\frac{1}{8}}$ for $t\ge1$.
Since $\beta\le\sqrt{A}\eps$ by \eqref{eq:abstract-basic},
\eqref{eq:abstract-window} and \eqref{eq:common-good-time} give
\[
 \begin{aligned}
 S(t)&\le S_I\le A\bigl[S(s)+\mathfrak b_s+\Theta_I+H_I^2
                            +H_I\mathfrak b_s+L_IH_IQ_I\bigr]\\
 &\le C\left[\frac{\sqrt{\eps\beta}}{t^{\frac{1}{4}}}
       +\frac{\beta}{\sqrt t}+\frac{\beta+\beta^2}{t^{\frac{1}{4}}}\right]
 \le C\left[1+(2+\beta)\sqrt{\frac{\beta}{\eps}}\right]
      \frac{\sqrt{\eps\beta}}{t^{\frac{1}{4}}}
 \le C\frac{\sqrt{\eps\beta}}{t^{\frac{1}{4}}}.
 \end{aligned}
\]
For $1\le t\le T$, this gives
\begin{equation}\label{eq:abstract-tail-proof}
 (1+t)S(t)^4\le2tS(t)^4\le C\eps^2\beta^2
             = C\eps^2\int_{\frac{t}{2}}^t B(\tau)\dd\tau.
\end{equation}

Under $\mathfrak M_T\le1$, \eqref{eq:abstract-basic} on $[0,1]$
and \eqref{eq:abstract-tail-proof} on $[1,T]$ give
\[
 \mathfrak M_T
 \le \max\left\{2^{\frac{1}{4}}A\eps,
       \left(C\eps^2\int_0^\infty B(\tau)\dd\tau\right)^{\frac{1}{4}}\right\}
 \le C\eps .
\]
Choose $\eps_*(A)\le1$ so that
$\max\{2^{\frac{1}{4}}A,C\}\eps_*(A)\le\frac{1}{2}$.
Then $\mathfrak M_1\le\frac{1}{2}$ by \eqref{eq:abstract-basic}.
The map $T\mapsto\mathfrak M_T$ is continuous. At a first time with
$\mathfrak M_T=1$, the preceding bound would give
$\mathfrak M_T\le\frac{1}{2}$. Hence
$\mathfrak M_T\le C\eps$ for every $T\ge1$ and
\eqref{eq:abstract-tail-proof} holds for every $t\ge1$.
\end{proof}

Proposition~\ref{prop:global} gives the unique global solution.
Set $S(t)=\|Z(t)\|_\infty$ and use $B$ from \eqref{eq:def-B}.
Since $Z\in C([0,\infty);H^2)$, the Sobolev embedding gives
continuity of $S$. Equations \eqref{eq:global-H2},
\eqref{eq:time-B-integral}, \eqref{eq:time-energy-equivalence}
and \eqref{eq:st-nash} give
\[
 \begin{aligned}
 \sup_{t\ge0}S(t)+\left(\int_0^\infty B(\tau)\dd\tau\right)^{\frac{1}{2}}
 &\le C_\delta\left(\sup_{t\ge0}\|Z(t)\|_{H^2}+\eps\right)
 \le C_\delta\eps,\\
 \mathcal T(t)&\le C_\delta B(t),\qquad
 S(t)^4\le C_\delta\|Z(t)\|_2^2\|F(t)\|_2^2
       \le C_\delta\eps^2\mathcal T(t)\le C_\delta\eps^2 B(t).
 \end{aligned}
\]
This proves \eqref{eq:abstract-basic}. Equations
\eqref{eq:window-inputs} and \eqref{eq:comparison-final} give
\eqref{eq:abstract-window} for every $t>0$ and almost every $s\in[0,t)$.
Choose $A$ large enough for
\eqref{eq:abstract-basic} and \eqref{eq:abstract-window}, and take
$\eps_\delta\le\eps_*(A)$.

By \eqref{eq:feedbackbounds}, $|X'(t)|\le C_\delta S(t)$.
Lemma~\ref{lemma:window-decay} and \eqref{eq:abstract-tail-proof} give
\[
 \begin{aligned}
 \sup_{t\ge0}(1+t)^{\frac{1}{4}}\bigl(S(t)+|X'(t)|\bigr)
 &\le C_\delta\sup_{t\ge0}(1+t)^{\frac{1}{4}}S(t)\le C_\delta\eps,\\
 (1+t)\bigl(S(t)+|X'(t)|\bigr)^4
 &\le C_\delta(1+t)S(t)^4
 \le C_\delta\eps^2\int_{\frac{t}{2}}^t B(\tau)\dd\tau\\
 &\le C_\delta\eps^2\int_{\frac{t}{2}}^\infty B(\tau)\dd\tau,
 \qquad t\ge1 .
 \end{aligned}
\]
This proves \eqref{eq:main-uniform} and \eqref{eq:main-tail}.
Since $B\in L^1(0,\infty)$, the tail integral vanishes as $t$ tends
to infinity. In the original coordinate, \eqref{eq:moving-frame}
gives $S(t)=\|q(t,\cdot)-\bar q(\cdot-\sigma t-X(t))\|_\infty$.
Taking the fourth root gives \eqref{eq:main-little-o}
and proves Theorem~\ref{thm:main}.

\subsection{Sharpness of the uniform decay rate}

\begin{proposition}\label{prop:uniform-sharpness}
For the fixed shock with the shift in \eqref{eq:moving-frame},
there are $c_\delta>0$ and $t_\delta\ge1$ such that
\begin{equation}\label{eq:uniform-rate-lower}
 \sup_{0<\|Z_0\|_{H^2}\le\eps_\delta}
 \frac{\|Z(t)\|_\infty}{\|Z_0\|_{H^2}}
 \ge c_\delta t^{-\frac{1}{4}},\qquad t\ge t_\delta .
\end{equation}
\end{proposition}

\begin{proof}
Fix $t>0$ and a nonzero real Schwartz function $f:\R\to\R^3$.
Take all limits in $a$ and $L$ at this time. For $L>0$ and $0<a\|f\|_{H^2}\le\eps_\delta$, let $Z^{a,L}$ be
the perturbation with initial value $af(\cdot-L)$. Since
$B\le C_\delta\|Z\|_{H^2}^2$ by \eqref{eq:def-B} and
\eqref{eq:feedbackbounds}, \eqref{eq:global-H2} and
\eqref{eq:time-F-controlled-by-B} give
\[
 \begin{aligned}
 &\frac{1}{a}\sup_{0\le\tau\le t}
       \left(\|Z^{a,L}(\tau)\|_{H^2}
                          +\|\partial_\tau Z^{a,L}(\tau)\|_2\right)
 \le \frac{C_\delta}{a}\sup_{0\le\tau\le t}
                     \|Z^{a,L}(\tau)\|_{H^2}\le C_\delta\|f\|_{H^2},\\
 &\frac{1}{a}\|Z^{a,L}(\tau)-Z^{a,L}(\rho)\|_{H^1}\\
 &\quad\le \frac{C}{a}
      \|Z^{a,L}(\tau)-Z^{a,L}(\rho)\|_2^{\frac{1}{2}}
      \|Z^{a,L}(\tau)-Z^{a,L}(\rho)\|_{H^2}^{\frac{1}{2}}
 \le C_\delta\|f\|_{H^2}\sqrt{\tau-\rho},
 \qquad 0\le\rho\le\tau\le t .
 \end{aligned}
\]
Taylor's formula in \eqref{eq:st-spatial-operator}, with
$\mathscr N_\sigma(\bar q)=0$ and $X'=\ell(Z^{a,L})$, gives on each
bounded interval $K$
\[
 \begin{aligned}
 &\frac{1}{a}\left\|\partial_\tau Z^{a,L}
       -D_q\mathscr N_\sigma(q_+)Z^{a,L}\right\|_{L^2(K+L)}\\
 &\quad\le \frac{1}{a}
   \left\|\bigl[D_q\mathscr N_\sigma(\bar q)
        -D_q\mathscr N_\sigma(q_+)\bigr]Z^{a,L}\right\|_{L^2(K+L)}\\
 &\qquad+\frac{1}{a}\int_0^1(1-\rho)
   \left\|D_q^2\mathscr N_\sigma(\bar q+\rho Z^{a,L})
              [Z^{a,L},Z^{a,L}]\right\|_{L^2(K+L)}\dd\rho\\
 &\qquad+\frac{|\ell(Z^{a,L})|}{a}
          \left(\|\partial_\xi Z^{a,L}\|_{L^2(K+L)}
                        +\|\bar q'\|_{L^2(K+L)}\right)\\
 &\quad\le C_{\delta,K}\|f\|_{H^2}
       \left(a\|f\|_{H^2}+\|\bar q-q_+\|_{W^{2,\infty}(K+L)}\right)
 \le C_{\delta,K}\|f\|_{H^2}\bigl(a\|f\|_{H^2}+e^{-c\delta L}\bigr)
 \end{aligned}
\]
for large $L$. The product rule in \eqref{eq:st-spatial-operator}
and $H^2\hookrightarrow W^{1,\infty}$ bound the Taylor integrand by $C_\delta\|Z^{a,L}\|_{H^2}^2$. The last two bounds
use \eqref{eq:feedbackbounds} and \eqref{eq:profilebounds}.

Let $a_n$ tend to zero and $L_n$ tend to infinity. The time estimate
and the compact embedding of $H^2(K)$ in $H^1(K)$ give a subsequence
along which $\frac{Z^{a_n,L_n}(\tau,\cdot+L_n)}{a_n}$ converges to $w$
in $C([0,t];H^1(K))$ for every bounded interval $K$. The global bounds give
$w\in L^\infty(0,t;H^2)\cap W^{1,\infty}(0,t;L^2)$.
After translation, the residual estimate on $K+L_n$ gives the
equation on $K$ in distributions. The initial value follows from
$\frac{Z^{a_n,L_n}(0,y+L_n)}{a_n}=f(y)$. Thus
\begin{equation}\label{eq:right-end-linear}
 \begin{aligned}
 w_\tau&=(\sigma I-A_+)w_y+D_+w_{yy},\qquad w(0)=f,\\
 A_+&=\begin{pmatrix}
       0&-1&0\\
       -\frac{R\theta_+}{v_+^2}&0&\frac{R}{v_+}\\
       0&\frac{R\theta_+}{c_vv_+}&0
      \end{pmatrix},
 \qquad D_+=\operatorname{diag}\left(0,\frac{\mu}{v_+},
                                      \frac{\kappa}{c_vv_+}\right).
 \end{aligned}
\end{equation}
The diagonal matrix with entries $\frac{R\theta_+}{v_+^2}$, $1$ and
$\frac{c_v}{\theta_+}$ symmetrizes $A_+$. Multiplying
\eqref{eq:right-end-linear} by the corresponding weighted $w$
and integrating by parts gives
\[
 \begin{aligned}
 &\frac{1}{2}\frac{\mathrm d}{\mathrm d\tau}
   \int_\R\left(\frac{R\theta_+}{v_+^2}w_1^2+w_2^2
                           +\frac{c_v}{\theta_+}w_3^2\right)\dd y\\
 &\qquad+\frac{\mu}{v_+}\|\partial_y w_2\|_2^2
              +\frac{\kappa}{v_+\theta_+}\|\partial_y w_3\|_2^2=0 .
 \end{aligned}
\]
The same identity for a difference with zero initial data gives
uniqueness. The embedding $H^1(K)\hookrightarrow C(\overline K)$
gives pointwise convergence along the subsequence. For each fixed
$y\in\R$, we obtain
\begin{equation}\label{eq:uniform-rate-limit}
 \begin{aligned}
 \sup_{0<\|Z_0\|_{H^2}\le\eps_\delta}
       \frac{\|Z(t)\|_\infty}{\|Z_0\|_{H^2}}
 &\ge \lim_{n\to\infty}
       \frac{|Z^{a_n,L_n}(t,y+L_n)|}{a_n\|f\|_{H^2}}
 =\frac{|w(t,y)|}{\|f\|_{H^2}} .
 \end{aligned}
\end{equation}

The zero eigenvalue of $-A_+$ is simple. For $k$ near zero,
the corresponding right eigenvector $r_0(k)$ of $-A_++ikD_+$
is analytic and has first component one. Multiplying its
eigenvalue by $ik$ gives the contact eigenvalue $\lambda_0(k)$
of $-ikA_+-k^2D_+$.
By \eqref{eq:right-end-linear},
$r_0(0)=(1,0,\frac{\theta_+}{v_+})^{\mathsf T}$ and
$(\frac{R}{v_+},0,\frac{c_v}{\theta_+})A_+=0$. Multiplying the
eigenvector equation by this row gives
\[
 \begin{aligned}
 \lambda_0(k)
 &=-k^2\frac{(\frac{R}{v_+},0,\frac{c_v}{\theta_+})D_+r_0(k)}
                {(\frac{R}{v_+},0,\frac{c_v}{\theta_+})r_0(k)}\\
 &=-k^2\frac{\frac{\kappa}{v_+^2}}{\frac{R+c_v}{v_+}}+O(|k|^3)
   =-\frac{\kappa}{\gamma c_vv_+}k^2+O(|k|^3).
 \end{aligned}
\]

Fix a nonnegative, even, nonzero $g\in C_c^\infty((-1,1))$.
For $t$ sufficiently large, take $f_t$ with Fourier transform
$t^{\frac{1}{4}}g(\sqrt t\,k)r_0(k)$ supported in
$|k|<\frac{1}{\sqrt t}$. Conjugate symmetry of $-A_++ikD_+$
and the normalization give $r_0(-k)=\overline{r_0(k)}$ for real $k$,
so $f_t$ is real and Schwartz. Plancherel's identity gives
\[
 \begin{aligned}
 \|f_t\|_{H^2}^2
 &\le \frac{t^{\frac{1}{2}}}{2\pi}
       \int_{-\frac{1}{\sqrt t}}^{\frac{1}{\sqrt t}}
           (1+k^2)^2|g(\sqrt t\,k)|^2|r_0(k)|^2\dd k\\
 &=\frac{1}{2\pi}\int_{-1}^1
        \left(1+\frac{k^2}{t}\right)^2|g(k)|^2
             \left|r_0\left(\frac{k}{\sqrt t}\right)\right|^2\dd k
 \le C_\delta\|g\|_2^2.
 \end{aligned}
\]
On $|k|\le1$ the expansion gives
$t\lambda_0(\frac{k}{\sqrt t})=-\frac{\kappa k^2}{\gamma c_vv_+}
+O(t^{-\frac{1}{2}})$ uniformly. For sufficiently large $t$,
$\operatorname{Re}e^{t\lambda_0(\frac{k}{\sqrt t})}
\ge\frac{1}{2}\exp(-\frac{\kappa k^2}{\gamma c_vv_+})$.
At $y=-\sigma t$, the Fourier inversion phase cancels the drift
factor $e^{i\sigma kt}$. Apply \eqref{eq:uniform-rate-limit} with
$f=f_t$ and use the first component of $w$ to obtain
\[
 \begin{aligned}
 t^{\frac{1}{4}}\!
 \sup_{0<\|Z_0\|_{H^2}\le\eps_\delta}\frac{\|Z(t)\|_\infty}{\|Z_0\|_{H^2}}
 &\ge\frac{t^{\frac{1}{4}}|w(t,-\sigma t)|}{\|f_t\|_{H^2}}\\
 &\ge\frac{1}{2\pi\|f_t\|_{H^2}}
     \int_{-1}^1 g(k)\operatorname{Re}
             e^{t\lambda_0(\frac{k}{\sqrt t})}\dd k\\
 &\ge\frac{1}{4\pi C_\delta\|g\|_2}
     \int_{-1}^1 g(k)\exp\left(-\frac{\kappa k^2}{\gamma c_vv_+}\right)\dd k
 >0 .
 \end{aligned}
\]
This proves \eqref{eq:uniform-rate-lower}.
\end{proof}


\begin{thebibliography}{99}
\bibitem{MN}
A.~Matsumura and K.~Nishihara,
On the stability of travelling wave solutions of a one-dimensional
model system for compressible viscous gas,
\emph{Japan J. Appl. Math.} \textbf{2} (1985), 17--25.
\href{https://doi.org/10.1007/BF03167036}{doi:10.1007/BF03167036}.

\bibitem{KM}
S.~Kawashima and A.~Matsumura,
Asymptotic stability of traveling wave solutions of systems for
one-dimensional gas motion,
\emph{Comm. Math. Phys.} \textbf{101} (1985), 97--127.
\href{https://doi.org/10.1007/BF01212358}{doi:10.1007/BF01212358}.

\bibitem{Liu}
T.-P.~Liu,
Shock waves for compressible Navier--Stokes equations are stable,
\emph{Comm. Pure Appl. Math.} \textbf{39} (1986), 565--594.
\href{https://doi.org/10.1002/cpa.3160390502}{doi:10.1002/cpa.3160390502}.

\bibitem{TZ}
S.~Tang and L.~Zhang,
Nonlinear stability of viscous shock waves for one-dimensional
nonisentropic compressible Navier--Stokes equations with a class
of large initial perturbation,
\emph{Acta Math. Sci. Ser. B Engl. Ed.} \textbf{38} (2018),
973--1000.
\href{https://doi.org/10.1016/S0252-9602(18)30797-5}{doi:10.1016/S0252-9602(18)30797-5}.

\bibitem{ZH}
K.~Zumbrun and P.~Howard,
Pointwise semigroup methods and stability of viscous shock waves,
\emph{Indiana Univ. Math. J.} \textbf{47} (1998), 741--872.
\href{https://doi.org/10.1512/iumj.1998.47.1604}{doi:10.1512/iumj.1998.47.1604}.

\bibitem{SK}
Y.~Shizuta and S.~Kawashima,
Systems of equations of hyperbolic-parabolic type with applications
to the discrete Boltzmann equation,
\emph{Hokkaido Math. J.} \textbf{14} (1985), 249--275.
\href{https://doi.org/10.14492/hokmj/1381757663}{doi:10.14492/hokmj/1381757663}.

\bibitem{MZG}
C.~Mascia and K.~Zumbrun,
Pointwise Green function bounds for shock profiles of systems with
real viscosity,
\emph{Arch. Ration. Mech. Anal.} \textbf{169} (2003), 177--263.
\href{https://doi.org/10.1007/s00205-003-0258-5}{doi:10.1007/s00205-003-0258-5}.

\bibitem{MZ}
C.~Mascia and K.~Zumbrun,
Stability of small-amplitude shock profiles of symmetric
hyperbolic-parabolic systems,
\emph{Comm. Pure Appl. Math.} \textbf{57} (2004), 841--876.
\href{https://doi.org/10.1002/cpa.20023}{doi:10.1002/cpa.20023}.

\bibitem{Raoofi}
M.~Raoofi,
$L^p$ asymptotic behavior of perturbed viscous shock profiles,
\emph{J. Hyperbolic Differ. Equ.} \textbf{2} (2005), 595--644.
\href{https://arxiv.org/abs/math/0408227}{arXiv:math/0408227}.

\bibitem{HRZ}
P.~Howard, M.~Raoofi, and K.~Zumbrun,
Sharp pointwise bounds for perturbed viscous shock waves,
\emph{J. Hyperbolic Differ. Equ.} \textbf{3} (2006), 297--373.
\href{https://doi.org/10.1142/S0219891606000793}{doi:10.1142/S0219891606000793}.

\bibitem{LZ}
T.-P.~Liu and Y.~Zeng,
Time-asymptotic behavior of wave propagation around a viscous shock
profile,
\emph{Comm. Math. Phys.} \textbf{290} (2009), 23--82.
\href{https://doi.org/10.1007/s00220-009-0820-6}{doi:10.1007/s00220-009-0820-6}.

\bibitem{KV}
M.-J.~Kang and A.~F.~Vasseur,
Contraction property for large perturbations of shocks of the
barotropic Navier--Stokes system,
\emph{J. Eur. Math. Soc.} \textbf{23} (2021), 585--638.
\href{https://doi.org/10.4171/JEMS/1018}{doi:10.4171/JEMS/1018}.

\bibitem{KVWbar}
M.-J.~Kang, A.~F.~Vasseur, and Y.~Wang,
Time-asymptotic stability of composite waves of viscous shock and
rarefaction for barotropic Navier--Stokes equations,
\emph{Adv. Math.} \textbf{419} (2023), Paper No.~108963.
\href{https://doi.org/10.1016/j.aim.2023.108963}{doi:10.1016/j.aim.2023.108963}.

\bibitem{KVW}
M.-J.~Kang, A.~F.~Vasseur, and Y.~Wang,
Time-asymptotic stability of generic Riemann solutions for compressible
Navier--Stokes--Fourier equations,
\emph{Arch. Ration. Mech. Anal.} \textbf{249} (2025), Paper No.~42.
\href{https://doi.org/10.1007/s00205-025-02116-w}{doi:10.1007/s00205-025-02116-w}.
The numbering used for this reference is that of
\href{https://arxiv.org/abs/2306.05604v1}{arXiv:2306.05604v1}.

\bibitem{HM}
F.~Huang and A.~Matsumura,
Stability of a composite wave of two viscous shock waves for
the full compressible Navier--Stokes equation,
\emph{Comm. Math. Phys.} \textbf{289} (2009), 841--861.
\href{https://doi.org/10.1007/s00220-009-0843-z}{doi:10.1007/s00220-009-0843-z}.

\bibitem{HW}
F.~Huang and T.~Wang,
Vanishing dissipation limit for non-isentropic Navier--Stokes equations
with shock data,
\emph{Sci. China Math.} \textbf{69} (2026), 2975--3006.
\href{https://doi.org/10.1007/s11425-025-2479-6}{doi:10.1007/s11425-025-2479-6}.

\bibitem{WW}
T.~Wang and Y.~Wang,
Nonlinear stability of planar viscous shock wave to three-dimensional
compressible Navier--Stokes equations,
\emph{J. Eur. Math. Soc.} \textbf{28} (2026), 505--557.
\href{https://doi.org/10.4171/JEMS/1486}{doi:10.4171/JEMS/1486}.

\bibitem{HX}
F.~Huang and L.~Xu,
Nonlinear stability of planar shock waves for the
Navier--Stokes--Fourier equations in $\R^3$,
\href{https://arxiv.org/abs/2611.11383}{arXiv:2611.11383}, 2026.

\bibitem{HHLX}
M.~He, F.~Huang, L.~Liu, and L.~Xu,
Time-decay rates toward rarefaction waves for the one-dimensional
isentropic compressible Navier--Stokes equations,
\href{https://arxiv.org/abs/2611.01956v1}{arXiv:2611.01956v1}, 2026.

\bibitem{LWX}
L.~Liu, S.~Wang, and L.~Xu,
Optimal decay rates to the contact wave for 1-D compressible
Navier--Stokes equations,
\emph{J. Math. Pures Appl.} \textbf{214} (2026), Paper No.~103948.
\href{https://arxiv.org/abs/2310.12747}{arXiv:2310.12747}.
\end{thebibliography}
\end{document}